\documentclass[11pt]{amsart}
\usepackage{hyperref}
\usepackage{amsfonts, amsmath, amssymb, amscd, amsthm, graphicx, setspace, enumitem, color,cleveref,comment}
\usepackage{hyperref}
\usepackage{pinlabel}
\usepackage{tikz-cd}
\usepackage{tikz}

\usepackage[alphabetic]{amsrefs}

\makeatletter
\def\blfootnote{\xdef\@thefnmark{}\@footnotetext}
\makeatother

\newtheorem{thm}{Theorem}[section]
\newtheorem{cor}[thm]{Corollary}
\newtheorem{lemma}[thm]{Lemma}
\newtheorem{prop}[thm]{Proposition}

\theoremstyle{definition}
\newtheorem{df}[thm]{Definition}
\theoremstyle{remark}
\newtheorem{rem}[thm]{Remark}

\newtheorem{Ex}[thm]{Example}

\newcommand{\Z}{{\mathbb Z}}
\newcommand{\R}{\mathbb{R}}
\newcommand{\Q}{\mathbb{Q}}
\newcommand{\C}{\mathbb{C}}

\newcommand{\A}{\mathbb{A}}
\newcommand{\Proj}{\mathbb{P}}
\newcommand{\PD}{\mathcal{P}\mathcal{D}}

\newcommand{\reg}{\text{reg}}
\newcommand{\Pic}{\operatorname{Pic}}
\newcommand{\Cl}{\operatorname{Cl}}
\newcommand{\Mod}{\operatorname{Mod}}
\newcommand{\Flat}{\textrm{Flat}}
\newcommand{\Art}{\textrm{Art}}
\newcommand{\Spec}{\operatorname{Spec}}

\newcommand{\Hom}{\operatorname{Hom}}
\newcommand{\Aut}{\textrm{Aut}}
\newcommand{\NE}{\operatorname{NE}}
\newcommand{\Ess}{\operatorname{Ess}}
\newcommand{\HP}{\text{HP}}

\newcommand{\Ab}{\textrm{Ab}}

\newcommand{\pd}{\operatorname{PD}}
\newcommand{\Def}{\textrm{Def}}
\newcommand{\HC}{\mathsf{HC}}
\newcommand{\Spr}{\mathsf{Spr}}
\newcommand{\Perv}{\operatorname{Perv}}

\newcommand{\End}{\operatorname{End}}
\newcommand{\sympl}{\operatorname{sympl}}
\newcommand{\Loc}{\operatorname{Loc}}
\newcommand{\IC}{\operatorname{IC}}
\newcommand{\BM}{\text{BM}}

\newcommand{\Mov}{\operatorname{Mov}}
\newcommand{\Amp}{\operatorname{Amp}}
\newcommand{\Nef}{\operatorname{Nef}}
\newcommand{\diff}{\operatorname{d}}
\newcommand{\PDef}{\operatorname{PDef}}
\newcommand{\Spf}{\operatorname{Spf}}

\newcommand{\limit}{\operatorname{lim}}
\newcommand{\Sym}{\operatorname{Sym}}
\newcommand{\Hilb}{\operatorname{Hilb}}
\newcommand{\PROJ}{\operatorname{Proj}}
\newcommand{\Cok}{\operatorname{Coker}}
\newcommand{\Stab}{\operatorname{Stab}}

\newcommand{\cH}{\mathcal{H}}
\newcommand{\cN}{\mathcal{N}}
\newcommand{\cM}{\mathcal{M}}
\newcommand{\cS}{\mathcal{S}}
\newcommand{\cO}{\mathcal{O}}
\newcommand{\cX}{\mathcal{X}}
\newcommand{\cU}{\mathcal{U}}
\newcommand{\cY}{\mathcal{Y}}
\newcommand{\cV}{\mathcal{V}}
\newcommand{\cL}{\mathcal{L}}
\newcommand{\cT}{\mathcal{T}}
\newcommand{\cF}{\mathcal{F}}

\newcommand{\cZ}{\mathcal{Z}}

\newcommand{\ol}[1]{\overline{{#1}}}

\begin{document}
\title[Partial Resolutions of Affine Symplectic Singularities]{Partial Resolutions of Affine Symplectic Singularities}

\author{Alberto San Miguel Malaney}

\date{\today}

\begin{abstract}

We explore the relationship between the Poisson deformation theory, birational geometry, and Springer theory of partial resolutions of affine symplectic singularities. Let $\rho: X' \rightarrow X$ be a crepant partial resolution of a conical affine symplectic singularity $X$. We show that the Poisson deformation functor of $X'$ is prorepresentable and unobstructed. Additionally, we define a version of the Namikawa Weyl group for these crepant partial resolutions. In particular, the Namikawa Weyl group of $X'$ is a parabolic subgroup of the Namikawa Weyl group of $X$ that is determined by the birational geometry of $X'$. If $\pi: Y \rightarrow X$ is a $\Q$-factorial terminalization of $X$ that covers $X'$, we show there is a natural functor from Poisson deformations of $Y$ to those of $X'$. Building on work of Namikawa in \cite{NaPD2}, we show that this morphism is a Galois covering and that the Galois group is the Namikawa Weyl group of $X'$. Finally, we put these partial resolutions and their universal deformations into the context of recent work of McGerty and Nevins in \cite{MN}, obtaining some preliminary results concerning their Springer theory. In particular, if the universal deformation of $X'$ is rationally smooth, we compute the cohomology of the fibers of $\rho$ in terms of the cohomology of the fibers of $\pi$ and the Namikawa Weyl group of $X'$.

\end{abstract}

\maketitle

\tableofcontents

\section{Introduction} \label[section]{Intro}

In this paper we study the deformation theory, birational geometry, and Springer theory of partial resolutions of conical affine symplectic singularities.

\subsection{Motivation}

The Springer resolution $\pi: \tilde{\cN} \rightarrow \cN$ of the nilpotent cone of a semisimple Lie algebra $\mathfrak{g}$ is a central object in geometric representation theory. In particular, it is the unique symplectic resolution of the nilpotent cone $\cN \subset \mathfrak{g}$. More generally, for every conjugacy class $\mathcal{P}$ of parabolic subalgebras of $\mathfrak{g}$, we obtain the following partial resolution of $\cN$: 
\begin{equation*}
\tilde{\cN}^{\mathcal{P}} := \{(x, \mathfrak{p}) | x \in \mathfrak{p} \cap \cN, \mathfrak{p} \in \mathcal{P} \}.
\end{equation*}
If we let $\mathcal{P}$ be the set of Borel subalgebras we recover the Springer resolution. Furthermore, each Borel subalgebra $\mathfrak{b} \subset \mathfrak{g}$ is contained in a unique element of each conjugacy class $\mathcal{P}$ of parabolic subalgebras, so we obtain birational and projective maps 
\begin{equation*}
\tilde{\cN} \xrightarrow{\eta} \tilde{\cN}^{\mathcal{P}} \xrightarrow{\xi} \cN.
\end{equation*}

Let $W$ be the Weyl group of $\mathfrak{g}$. Then we can take the Weyl group of the Levi subalgebra of any $\mathfrak{p} \in \mathcal{P}$ to obtain a parabolic subgroup $W_{\mathcal{P}} \subset W$ of the Weyl group. This gives us a bijection between conjugacy classes of parabolic subalgebras of $\mathfrak{g}$ and parabolic subgroups of $W$, both of which can also be classified by subsets of the simple roots of $\mathfrak{g}$. 

The relationship between the geometry of these partial resolutions and the combinatorics of the Weyl group was studied extensively in \cite{BM}. For example, \cite[2.8, Corollary]{BM} states that for any $x \in \cN$,
\begin{equation} \label{FibersInv}
H^*(\xi^{-1}(x), \Q) \cong H^*(\pi^{-1}(x), \Q)^{W_{\mathcal{P}}}.
\end{equation}
In proving results such as these it is essential that we can drop the nilpotent condition to obtain a Grothendieck-Springer resolution as below, where $\tilde{\eta}$ is generically a $W_{\mathcal{P}}$-cover and $\tilde{\pi} := \tilde{\xi} \circ \tilde{\eta}$ is generically a $W$-cover.
\[\begin{tikzcd}[row sep=scriptsize]
	{\tilde{\cN}} && {\tilde{\cN}^{\mathcal{P}}} && \cN \\
	\\
	{\tilde{\mathfrak{g}}} && {\tilde{\mathfrak{g}}^{\mathcal{P}}} && {\mathfrak{g}}
	\arrow["{\tilde{\eta}}", from=3-1, to=3-3]
	\arrow["{\tilde{\xi}}", from=3-3, to=3-5]
	\arrow[hook, from=1-1, to=3-1]
	\arrow["\eta", from=1-1, to=1-3]
	\arrow[hook, from=1-3, to=3-3]
	\arrow["\xi", from=1-3, to=1-5]
	\arrow[hook, from=1-5, to=3-5]
\end{tikzcd}\]

For a general conical affine symplectic singularity $X$, with a $\Q$-factorial terminalization $\pi: Y \rightarrow X$, the generalization of the Grothendieck-Springer resolution is the universal graded deformation. These have been studied extensively by Namikawa, who in a series of papers (primarily \cite{NaFlops, NaPD1, NaPD2}) showed that the Poisson deformation functors $\pd_X$ and $\pd_Y$ of $X$ and $Y$ are prorepresentable and unobstructed, and that there is a natural map 
\begin{equation*}
\pi_*: \pd_Y \rightarrow \pd_X.
\end{equation*}
In particular, this implies that $X$ and $Y$ have universal formal deformations over the formal neighborhood of zero in their tangent spaces. Using the conical action on $X$, which also lifts to $Y$, these formal deformations can be algebraicized to obtain universal graded deformations $\cX$ and $\cY$ (see \cite[Proposition 2.6, Proposition 2.12]{Lo}), which fit into the following diagram:
\begin{equation} \label{IntroUniv}
\begin{tikzcd} [row sep=scriptsize]
	Y && \cY & {\A_Y} \\
	\\
	X && \cX & {\A_X \cong \A_Y/W_X}
	\arrow[hook, from=1-1, to=1-3]
	\arrow["\pi", from=1-1, to=3-1]
	\arrow[from=1-3, to=1-4]
	\arrow["{\tilde{\pi}}", from=1-3, to=3-3]
	\arrow["{\pi_*}", from=1-4, to=3-4]
	\arrow[hook, from=3-1, to=3-3]
	\arrow[from=3-3, to=3-4]
\end{tikzcd}
\end{equation}
In the above diagram, $X$ and $Y$ embed in $\cX$ and $\cY$ as the central fibers, and $\pi_*$ is the quotient by $W_X$. Additionally, $\A_Y$ and $\A_X$ are both affine spaces produced by giving $\pd_{Y}(\C[\epsilon])$ and $\pd_{X}(\C[\epsilon])$ the natural scheme structure. 

In the case where $\pi$ is a symplectic resolution, the viability of symplectic Springer theory was further established by McGerty and Nevins in \cite{MN} using the above constructions. In particular, \cite{MN} employs a nearby cycles construction over the universal graded deformation to show that the symplectic Springer sheaf $\Spr := \pi_*\C_Y$ is constructible with respect to the stratification by symplectic leaves on $X$. Additionally, letting $\HC := \tilde{\pi}_*{\C_{\cY}}$ be the symplectic Harish-Chandra sheaf, there is a restriction map 
\begin{equation*}
\End(\HC) \cong \C[W_X] \rightarrow \End(\Spr),
\end{equation*}
which generalizes the classical isomorphism $\C[W] \cong \End(\Spr)$ in Springer theory. It also produces an action of $W_X$ on the cohomology of the fibers of $\pi$.

\subsection{Outline of Results}

In this paper we introduce partial resolutions into the above story, treating them as the appropriate analogue of parabolic subalgebras in the symplectic setting. We introduce a natural definition of a symplectic partial resolution, which we go on to show is equivalent to being crepant and to being covered by a $\Q$-factorial terminalization. The main construction of the paper is a generalization of the Namikawa Weyl group to crepant partial resolutions $\rho: X' \rightarrow X$ of conical affine symplectic singularities, such that the Namikawa Weyl group $W_{X'}$ of a partial resolution is a parabolic subgroup of $W_X$.

In order to proceed we prove that the Poisson deformation functor $\pd_{X'}$ of $X'$ is also prorepresentable and unobstructed and that $X'$ has a universal graded deformation over 
$$
\A_{X'} := \Spec \Sym(\pd_{X'}(\C[\epsilon])^*).
$$ 
Given a sequence of maps
\begin{equation*}
Y \xrightarrow{\pi'} X' \xrightarrow{\rho} X
\end{equation*}
where the composition $\pi: Y \rightarrow X$ is a $\Q$-factorial terminalization of $X$, we also define and study the pushforward functors 
\begin{equation*}
\pd_{Y} \xrightarrow{\pi'_*} \pd_{X'} \xrightarrow{\rho_*} \pd_X.
\end{equation*}
This culminates in \Cref{MainThm}, the main theorem of the paper, which states that there is a diagram of universal graded deformations
\[\begin{tikzcd}[row sep=scriptsize]
	Y && \cY & {\A_Y} \\
	\\
	{X'} && {\cX'} & {\A_{X'} \cong \A_{Y}/W_{X'}} \\
	\\
	X && \cX & {\A_{X} \cong \A_{Y}/W_{X}}
	\arrow["{i_{Y}}", hook, from=1-1, to=1-3]
	\arrow["{\pi'}"', from=1-1, to=3-1]
	\arrow[from=1-3, to=1-4]
	\arrow["{\tilde{\pi}'}"', from=1-3, to=3-3]
	\arrow["{\pi'_*}"', from=1-4, to=3-4]
	\arrow["{i_{X'}}", hook, from=3-1, to=3-3]
	\arrow["\rho"', from=3-1, to=5-1]
	\arrow[from=3-3, to=3-4]
	\arrow["{\tilde{\rho}}"', from=3-3, to=5-3]
	\arrow["{\rho_*}"', from=3-4, to=5-4]
	\arrow["{i_X}", hook, from=5-1, to=5-3]
	\arrow[from=5-3, to=5-4]
\end{tikzcd}\]
extending diagram \eqref{IntroUniv}, where $\pi'_*$ is the quotient by $W_{X'}$. The diagram above can be seen as a justification for calling $W_{X'}$ the Namikawa Weyl group. 

One key result and motivation for this paper is Namikawa's result in \cite{NaMori} that a $\Q$-factorial terminalization $\pi: Y \rightarrow X$ of a conical affine symplectic singularity is also a relative Mori Dream Space. We explain that this implies that crepant partial resolutions of $X$ are classified by faces of the Mori fan of $\Mov(Y/X)$. This provides a combinatorial description of the crepant partial resolutions of $X$, and it is also the setting in which we first define the Namikawa Weyl group of a partial resolution. 

In particular, to every simple generator $s$ of the Namikawa Weyl group $W_X$, we assign a hyperplane $H_s \subset H^2(Y^{\reg}, \C)$, which bounds $\Mov(Y/X)$. Given a face $F$ of $\Mov(Y/X)$, we define the Namikawa Weyl group of the corresponding crepant partial resolution to be the parabolic subgroup of $W_X$ generated by all the simple generators $s$ such that $F \subset H_s$.

The nilpotent cone also has the special property that it has a unique $\Q$-factorial terminalization. Whenever this is the case, each face of the Mori fan will be given by intersecting $\Mov(Y/X)$ with a unique subset of the $H_s$ above. Thus, we can view this as a way of expressing the bijection between conjugacy classes of parabolic subalgebras, parabolic subgroups of the Weyl group, and subsets of the simple roots. It also shows that these are in natural bijection with crepant partial resolutions. 

These bijections generalize to a bijection $[\rho: X' \rightarrow X] \mapsto W_{X'}$ between crepant partial resolutions and parabolic subgroups of $W_X$ anytime $X$ has a unique $\Q$-factorial terminalization. As proved in \Cref{AllSubgroups}, this map is not a bijection for every conical affine symplectic singularity, but it is in general surjective. This is because if there is more than one $\Q$-factorial terminalization, there will be faces in the interior of $\Mov(Y/X)$ that correspond to small contractions that allow us to flop between $\Q$-factorial terminalizations. As a corollary of this theorem we also prove in \Cref{NilpCovers} that every crepant partial resolution of a nilpotent cone $\cN$ is of the form $\tilde{\cN}^{\mathcal{P}} \rightarrow \cN$ described at the beginning of the introduction.

In the final section of the paper we begin to study Springer theory of crepant partial resolutions as described above, with the added assumption that $\pi: Y \to X$ is a symplectic resolution. We briefly highlight two instances of interest. The first is a generalization of equation \eqref{FibersInv}. There are simple examples of conical affine symplectic singularities with trivial Namikawa Weyl group but nontrivial crepant partial resolutions, which shows the analogous equation could not be true in general. However, we show in \Cref{SymplFibersInv} that the analogous statement 
$$
H^*(\rho^{-1}(x), \C) \cong H^*(\pi^{-1}(x), \C)^{W_{X'}}
$$ 
is true as long as the universal graded deformation $\cX'$ of $X$ is rationally smooth.

Secondly, we introduce the algebra $\End(\Spr') \cong H^{\BM}_{top}(Z')$ where $\Spr' := \pi'_*(\C_Y)$ and $Z' := Y \times_{X'} Y$. As explained in diagram \eqref{End}, this algebra has a $W_{X'}$-action that is compatible with the natural map $\rho_*: \End(\Spr') \rightarrow \End(\Spr)$, which we show is an injection, and with the $W_X$-action on $\Spr$. In particular, this algebra provides a finer invariant than $W_{X'}$ that we can attach to a crepant partial resolution. 

\subsection{Summary}

In \Cref{Fold} we review group actions on simply-laced Dynkin diagrams and fix a formalism for the folded Dynkin diagrams we obtain by taking a quotient. We note that unlike in some formalisms for folded Dynkin diagrams, we allow the group action to permute nodes that have edges connecting them. For example we will have $A_4/(\Z/2\Z) \cong B_2$. In \Cref{Background} we review some basics about symplectic singularities. We recall some facts about $\Q$-factorial terminalizations as well as the definition of the Namikawa Weyl group. 

In \Cref{Bir} we begin to study crepant partial resolutions of conical affine symplectic singularities. In particular, we define them and prove a few equivalent characterizations of them. We review the main theorem of \cite{NaMori}, which roughly states that if $\pi: Y \rightarrow X$ is a $\Q$-factorial terminalization of a conical affine symplectic singularity, then it is also a relative Mori Dream space. Using this we obtain a combinatorial description of crepant partial resolutions by faces of the Mori fan of $\Mov(Y/X)$, and we use this description to define the Namikawa Weyl group of a crepant partial resolution. In \Cref{AllSubgroups} we prove the key result that the map from crepant partial resolutions to parabolic subgroups of $W_X$ is surjective, and bijective if and only if $X$ has a unique $\Q$-factorial terminalization.

In \Cref{Def} we begin to study the deformation theory of crepant partial resolutions. We use many of the techniques in \cite{NaFlops} and \cite{NaPD1} to show that $\pd_{X'}$ is prorepresentable and unobstructed for any crepant partial resolution $\rho: X' \rightarrow X$. Additionally, we define the pushforward functors $\pd_Y \rightarrow \pd_{X'} \rightarrow \pd_X$ on deformations, where $\pi: Y \rightarrow X$ is a $\Q$-factorial terminalization that covers $\rho$. In \Cref{Klein} we continue by studying the deformation theory of crepant partial resolutions of Kleinian singularities. The main strategy is that these partial resolutions have finitely many singular points, with singularities that are also Kleinian and correspond to a subset of the Dynkin diagram corresponding to the original Kleinian singularity. Using this we are able to prove the main theorem for partial resolutions of Kleinian singularities.

In \Cref{Main} we adapt the argument in \cite[Theorem 1.1]{NaPD2} to prove the main theorem. Given $Y \rightarrow X' \rightarrow X$ as above, we proceed by looking at the partial resolution and $\Q$-factorial terminalization locally over a small neighborhood of a point in each codimension two symplectic leaf of $X$. Thus, we obtain full and partial resolutions of the Kleinian singularity given by taking a transverse slice to the symplectic leaf. We are able to use the results from the previous section to describe the pushforward maps $\pd_Y \rightarrow \pd_{X'} \rightarrow \pd_X$ in the general case. This provides an explicit description of the universal graded deformation of $\pd_{X'}$, proving the main theorem, \Cref{MainThm}.

Finally, in \Cref{Springer} we use the main theorem to study the symplectic Springer theory of crepant partial resolutions, using many of the techniques from \cite{MN}. For this section we assume that $\pi: Y \to X$ is a symplectic resolution. We explain how many of the results in that paper that establish the viability of symplectic Springer theory carry over for crepant partial resolutions. Thus, we are able to define the partial symplectic Springer and Harish-Chandra sheaves, the former of which is also constructible with respect to the stratification by symplectic leaves. We end with a generalization of \cite[2.8, Corollary]{BM}, showing that if the universal graded deformation $\cX'$ is rationally smooth, then $H^*(\rho^{-1}(x), \C) \cong H^*(\pi^{-1}(x), \C)^{W_{X'}}$.

\subsection{Notation and Conventions}

All varieties and schemes will be assumed to be over $\C$. The results of this paper primarily involve algebraic varieties, but if we use the word ``neighborhood" it will usually mean an open set in the complex analytic topology. Germs will also be taken in the complex analytic topology. Finally, $\epsilon$ will be assumed to square to zero, so that $\C[\epsilon]$ will always denote the dual numbers. 

Throughout the paper, $X$ will always be a conical affine symplectic singularity, $\rho: X' \rightarrow X$ will always be a crepant partial resolution, and $\pi: Y \rightarrow X$ will always be a $\Q$-factorial terminalization of $X$, except in \Cref{Springer}, where $\pi$ will always be a symplectic resolution. If $\pi$ covers $\rho$, then $\pi': Y \rightarrow X'$ will be the map making $Y \rightarrow X' \rightarrow X$ commute. Adding a prime symbol will denote the ``partial" version. For example, the partial Springer and Harish-Chandra sheaves will be denoted $\Spr'$ and $\HC'$. 

Whenever there is a group action on a Dynkin diagram, we will denote the unfolded Dynkin diagram and all of its data (such as the Cartan and Weyl group) with a hat symbol. The folded Dynkin diagram and its data will have no hat. In \Cref{Klein} there are no group actions so we will drop the hats entirely. Additionally, hats are also used to denote formal completions, but the notation $\hat{\mathfrak{h}}$ will only appear to denote the Cartan of a simply laced Dynkin diagram.

Finally, we call a map of normal varieties $f: Y \rightarrow X$ a Galois cover if it is finite, surjective, and the degree of the map is $|\Aut(Y/X)|$. This last condition is equivalent to the extension $[\C(Y), \C(X)]$ of function fields being Galois. In the literature this is often called ``generically Galois." We will also refer to $\Aut(Y/X)$ as the Galois group. Given a Galois cover $f: Y \rightarrow X$ and an open set $U \subset X$ in the complex analytic topology, we will also refer to the map $f^{-1}(U) \rightarrow U$ as a Galois cover.

\subsection{Acknowledgments}

I want to thank David Ben-Zvi for his continuous help and support, without which this paper would not exist. I want to thank Rok Gregoric and Joakim F\ae rgeman for helping me work through a number of details. I want to thank Liam Riordan for pointing out that I could relax the condition in \Cref{SymplFibersInv} from smoothness to rational smoothness. I also want to thank Ivan Losev, Travis Schedler, Gwyn Bellamy, Chiara Damiolini, Paul Levy, Alastair Craw, Lewis Topley, JiWoong Park, Yaochen Wu, Austin Hubbard, and Jose Guzman for helpful discussions and suggestions. Finally, I want to thank the referees at \textit{Adv. Math.}, whose thorough comments greatly improved the paper.

\section{Lie Theory and Folding} \label[section]{Fold}

Let $\hat{D}$ be a simply laced Dynkin diagram and let $\hat{\mathfrak{h}}$ and $\hat{W}$ be the corresponding Cartan subalgebra and Weyl group, with $\langle \ , \rangle$ the natural inner product on $\hat{\mathfrak{h}}$ given by the Killing form. Furthermore, let $I$ be the set of nodes of $\hat{D}$, and let $\alpha_i \in \hat{\mathfrak{h}}$ and $s_i \in \hat{W}$ be the corresponding simple root and simple reflection, for each $i \in I$. Let $G$ be a finite group which acts on $\hat{D}$ by graph automorphism. We are going to use some results from \cite[Section 1]{St} to obtain a folded Dynkin diagram $D := \hat{D}/G$ with Cartan $\mathfrak{h} := \hat{\mathfrak{h}}^G$ and Weyl group $W := \hat{W}^G$.

Let $J := I/G$ be the set of orbits of the action of $G$ on $I$. For each $j \in J$, let $\beta_j = \sum_{i \in j} \alpha_i$. Since the $\alpha_i$ form a basis for $\hat{\mathfrak{h}}$, the permutation action of $G$ on the $\alpha_i$ induces an action of $G$ on $\hat{\mathfrak{h}}$, and $\mathfrak{h} := \hat{\mathfrak{h}}^G$ has a basis given by the $\beta_j$. 

\begin{rem}

The $\beta_j$ we define below have a slightly different normalization than the one given by the $a_j$ in \cite[p. 14]{St}, which will result in the short and long roots of the folded Dynkin diagrams being flipped. 

\end{rem}

In \cite[Section 12.2]{Hum} (see also \cite[Section 4.1]{NaPD1}) graph automorphisms of Dynkin diagrams are classified. The automorphism groups are $S_3$ for $D_4$ and $\Z/2\Z$ for $A_n, \ n \geq 2$, $D_n, \ n \geq 5$, and $E_6$. It follows from this classification that all non-trivial faithful group actions on simply laced Dynkin diagrams are of the following forms:
\begin{equation*}
\Z/2\Z \ \rotatebox[origin=c]{-90}{$\circlearrowright$}\  A_n \text{ for } n \geq 2, \ \Z/2\Z \ \rotatebox[origin=c]{-90}{$\circlearrowright$}\  D_n \text{ for } n \geq 3, \ \Z/3\Z \ \rotatebox[origin=c]{-90}{$\circlearrowright$}\  D_4, \ S_3 \ \rotatebox[origin=c]{-90}{$\circlearrowright$}\  D_4, \ \Z/2\Z \ \rotatebox[origin=c]{-90}{$\circlearrowright$}\  E_6.
\end{equation*}
In all of the above cases, the matrix $[c_{jk}]_{j, k \in J}$ given by 
\begin{equation} \label{Cartan}
c_{jk} = \frac{2 \langle \beta_j, \beta_k \rangle}{\langle \beta_j, \beta_j \rangle}
\end{equation}
is the Cartan matrix of a simple Lie algebra. We define $D = \hat{D}/G$ to be the Dynkin diagram corresponding to this Lie algebra. In particular we have:
\begin{gather*}
A_n/(\Z/2\Z) \cong B_{\lceil n/2 \rceil} \text{ for } n \geq 2, \ D_n/(\Z/2\Z) \cong C_{n - 1} \text{ for } n \geq 3 \\
D_4/(\Z/3\Z) \cong G_2, \ D_4/S_3 \cong G_2, \ E_6/(\Z/2\Z) \cong F_4.
\end{gather*}

\begin{Ex}

We shall compute $A_6/(\Z/2\Z) \cong B_3$ explicitly. We have 
\begin{equation*}
\beta_1 = \alpha_1 + \alpha_6, \hspace{5pt} \beta_2 = \alpha_2 + \alpha_5, \hspace{5pt} \beta_3 = \alpha_3 + \alpha_4.
\end{equation*}
Thus using the inner product from $A_6$ we obtain 
\begin{equation*}
\langle \beta_1, \beta_1 \rangle = \langle \beta_2, \beta_2 \rangle = 4, \hspace{5pt} \langle \beta_3, \beta_3 \rangle = 2, \hspace{5pt} \langle \beta_1, \beta_2 \rangle = \langle \beta_2, \beta_3 \rangle = -2, \hspace{5pt} \langle \beta_1, \beta_3 \rangle = 0. 
\end{equation*}
Using equation \eqref{Cartan} above we obtain the Cartan matrix for $B_3$. In particular, $\beta_1$ and $\beta_2$ are the long roots, while $\beta_3$ is the short root. 

\end{Ex}

Since the action of $G$ on $\hat{\mathfrak{h}}$ simply permutes the simple roots, conjugation by any element of $G$ permutes the simple reflections $s_i$. Thus $G$ normalizes $\hat{W}$, and we let $G$ act on $\hat{W}$ by conjugation. 

\begin{lemma}[{\cite[Lemma 1.2]{NaPD2}}]

The subset of $\hat{W}$ which preserves $\mathfrak{h}$ is $W := \hat{W}^G$.

\end{lemma}

Given $j \in J$, let $\hat{W}_j \subset \hat{W}$ be the parabolic subgroup generated by $s_i$ for all $i \in j$. Furthermore, let $w^0_j$ be the longest element of $\hat{W}_j$.

\begin{thm}[{\cite[1.30, 1.32]{St}}] \label[thm]{Stein}

For each $j \in J$, $w^0_j \in W$ and the restriction of $w^0_j$ to $\mathfrak{h}$ is a simple reflection corresponding to $\beta_j$. Additionally, the $w^0_j$ generate $W$. 

\end{thm}

In particular, \Cref{Stein} tells us that $W$ is simply the Weyl group corresponding to the Dynkin diagram $D$. For the remainder of the paper, whenever we have a Weyl group $W$ of a folded Dynkin diagram $D$ as above, we define $W$ to be a Coxeter group with the $w^0_j$ as the simple generators.

Finally, we want to compare the parabolic subgroups of $W$ with those of $\hat{W}$. Thus, let $J' \subset J$ and let $I' \subset I$ be such that $i \in I'$ if and only if $i \in j$ for some $j \in J'$. We further let $\hat{W}'$ be the parabolic subgroup of $\hat{W}$ generated by $s_i$ for each $i \in I'$.

\begin{prop} \label[prop]{FoldingResult}

The intersection $\hat{W}' \cap W$ is the parabolic subgroup $W' \subset W$ generated by the set $\{w^0_j \ | \ j \in J'\}$.

\end{prop}

\begin{proof}

By definition $w^0_j$ is a product of simple generators $s_i \in \hat{W}$ satisfying $i \in j$. Thus $w^0_j \in \hat{W}'$ for each $j \in J'$. Furthermore, by \Cref{Stein} we have $w^0_j \in W$ so $W' \subset \hat{W}' \cap W$.

We show the other inclusion by contradiction. Assume there is some $w \in \hat{W}' \cap W$ such that $w \notin W'$. We have an expansion $w = w^0_{j_1} \dots w^0_{j_n}$ for $j_1, \dots, j_n \in J$, and some $j_l$ is not in $J'$. We can take $w$ to have minimal length (as a word in $W$) out of elements of $\hat{W}' \cap W$ which are not in $W'$. Thus, $j_n \notin J'$, because otherwise multiplying by $w^0_{j_n}$ on the right would create a shorter word in $\hat{W}' \cap W$ but not in $W'$. 

If we let $\Delta = \{\beta_j\}_{j \in J}$ and $\Phi = W\Delta$ be the simple roots and the roots respectively, then $(\mathfrak{h}, \Delta, \Phi, W|_{\mathfrak{h}})$ satisfies conditions \cite[1.1-1.6]{St}. Clearly as a word in $W$, $l(ww^0_{j_n}) = l(w) - 1$ where $l$ is the word length function. Thus, by \cite[1.8]{St}, $ww^0_{j_n}$ makes one less positive root of $\Phi$ negative than $w$.

Since $w \in \hat{W}'$, \cite[1.15]{St} tells us that as an element of $\hat{W}$ it permutes all the positive roots of $\hat{\mathfrak{h}}$ whose support is not contained in $I'$. In particular it sends $\alpha_i$ to a positive root for each $i \in j_n$. Thus since $\beta_{j_n}$ is a linear combination of these $\alpha_i$, it must map $\beta_{j_n}$ to the positive cone of $\hat{\mathfrak{h}}$. Since the positive cone of $\mathfrak{h}$ is contained in the positive cone of $\hat{\mathfrak{h}}$ and $\beta_{j_n} \in \Phi$, $w$ must send $\beta_{j_n}$ to a positive root of $\Phi$. By \cite[1.7]{St}, $w^0_{j_n}$ permutes all the positive roots of $\Phi$ except for $\beta_{j_n}$, which it negates. Thus if $\{r_1, \dots, r_k\} \subset \Phi$ is the set of positive roots that $w$ sends to negative roots, the analogous set for $ww^0_{j_n}$ is $\{w^0_{j_n}(r_1), \dots, w^0_{j_n}(r_k), \beta_{j_n}\}$. However this contradicts the previous paragraph, so $\hat{W}' \cap W \subset W'$ and we are done. 
\end{proof}

We will now briefly recall the relationship between Dynkin diagrams and subregular Springer fibers. Given a Lie algebra $\mathfrak{g}$ with Springer resolution $\pi: \tilde{\cN} \rightarrow \cN \subset \mathfrak{g}$, the subregular orbit $\mathbb{O}_{\text{subreg}} \subset \cN$ is the unique nilpotent orbit which is open in the complement of the regular nilpotent orbit. A number of equivalent conditions to being subregular are stated in \cite[3.10, Theorem 2]{St2}. These conditions also imply that the subregular nilpotent orbit can be characterized as the Richardson orbit of any minimal parabolic. 

Let $D$ be the Dynkin diagram of $\mathfrak{g}$ and let $\hat{D}$ be the unique simply laced Dynkin diagram such that $D \cong \hat{D}/\Gamma$, where $\Gamma \ \rotatebox[origin=c]{-90}{$\circlearrowright$} \ \hat{D}$ is a faithful group action which does not permute nodes that are connected by an edge. Additionally, let $G$ be a simple Lie algebra of adjoint type with Lie algebra $\mathfrak{g}$. Let $x \in \cN$ be a subregular nilpotent element and let $\beta$ be a simple root of $\mathfrak{g}$. Furthermore, let $\mathcal{P}_{\beta}$ be the corresponding conjugacy class of minimal parabolic subalgebras, which is in bijection with $G/P$, for any parabolic subgroup $P$ of type $\beta$. By \cite[3.10, Theorem 2]{St2}, $x$ is in the nilradical of $|\beta|^2/|\beta_{\text{min}}|^2$ parabolic subalgebras $\mathfrak{p} \in \mathcal{P}_{\beta}$. Note that $|\beta|^2/|\beta_{\text{min}}|^2$ is the same as the number of nodes in $\hat{D}$ that map to the node in $D$ corresponding to $\beta$. Each of the these parabolic subalgebras $\mathfrak{p}$ produces a projective line in $\pi^{-1}(x)$ given by the set of Borel subalgebras contained in $\mathfrak{p}$, which is isomorphic to $P/B \cong \Proj^1$. Thus, $\pi^{-1}(x)$ consists of a union of projective lines that are in bijection with the nodes of $\hat{D}$. By \cite[3.10, Theorem 2]{St2} or \cite[Section 6.3]{Slo}, two of these projective lines intersect at a point if the corresponding nodes of $\hat{D}$ are connected and do not intersect otherwise. Furthermore, by \cite[7.5, Proposition]{Slo}, when $D$ is not simply laced the action of $G$ on $\hat{D}$ is isomorphic to the action of $Z_{G}(x)/Z_{G}(x)^{\circ}$, the component group of the centralizer of $x$, on the projective lines comprising $\pi^{-1}(x)$. Note that if $D$ is simply laced then this action is trivial since the projective lines correspond to minimal parabolic subalgebras in different conjugacy classes. Finally, two projective lines are permuted by $Z_{G}(x)/Z_{G}(x)^{\circ}$ if and only if they are permuted by the monodromy action of $\pi_1(\mathbb{O}_{\text{subreg}}, x)$. The results described in this paragraph produce the following result:

\begin{cor} \label[cor]{SubregFiber}

Given any $x \in \mathbb{O}_{\text{subreg}}$, there is a canonical bijection between the simple roots of $\mathfrak{g}$ and the $\pi_1(\mathbb{O}_{\text{subreg}}, x)$-orbits of the irreducible components of $\pi^{-1}(x)$. This bijection is given by sending each simple root $\beta$ to the set of projective lines given by $\{(x, \mathfrak{b}) \ | \ \mathfrak{b} \subset \mathfrak{p}\}$, for every $\mathfrak{p} \in \mathcal{P}_{\beta}$ with $x$ in the nilradical of $\mathfrak{p}$. 

\end{cor}

\section{Preliminaries and Background} \label[section]{Background}

\begin{df}

A symplectic singularity (or symplectic variety) is a normal algebraic variety $X$ with a symplectic form $\omega$ on $X_{\reg}$ such that for any (equivalently, some) resolution of singularities $\pi: \tilde{X} \rightarrow X$, $\pi^*(\omega)$ extends to a regular 2-form on $\tilde{X}$. If it extends to a symplectic form then $\tilde{X}$ is called a symplectic resolution. 

\end{df}

The symplectic form induces a Poisson bracket on $\cO(X_{\reg})$, and since $X$ is normal, this bracket extends to all of $X$, making $X$ a Poisson scheme. Throughout the paper we shall study the universal Poisson deformations of affine symplectic singularities and their (partial) symplectic resolutions. In order to ensure that universal formal deformations algebraicize to deformations over affine spaces, we need the following condition as well.

\begin{df} \label[df]{ConicDef}

A symplectic singularity $X$ is conical if it has an algebraic $\C^*$-action such that $\cO(X)$ has non-negative weights,  $\cO(X)^{\C^*} = \C$, and $\omega_X$ is positively weighted with respect to $\C^*$.

\end{df}

\begin{rem}

The conditions that the $\C^*$-action has non-negative weights and $\cO(X)^{\C^*} = \C$ are equivalent to the $\C^*$-action contracting $X$ to a unique cone point. 

\end{rem}

The following theorem is one of the most important facts about symplectic singularities:

\begin{thm}[{\cite[Theorem 2.3]{Kal}}]

Every symplectic singularity has a finite stratification by symplectic leaves.

\end{thm}

Throughout this paper we will be particularly interested in the codimension two symplectic leaves. Thus we recall (see \cite[Corollary 1]{NaNote}) that a symplectic singularity has terminal singularities if and only if it has no codimension two symplectic leaves. 

In general not every symplectic singularity has a symplectic resolution, but they always have a $\Q$-factorial terminalization, which is defined below:

\begin{df}

Given a symplectic singularity $X$, a $\Q$-factorial terminalization (or simply $\Q$-terminalization) $\pi: Y \rightarrow X$ is a crepant birational projective morphism such that $Y$ is a normal $\Q$-factorial variety with terminal singularities. 

\end{df}

We recall some essential facts about $\Q$-factorial terminalizations below:

\begin{prop}[{[\cite[Proposition 2.1]{LoSRA}}] \label[thm]{QFac}

Let $X$ be a symplectic singularity:

\begin{enumerate}
\item $X$ has a $\Q$-factorial terminalization, $\pi: Y \rightarrow X$.
\item $Y$ is a symplectic singularity and $\pi$ is a Poisson map.
\item If $X$ is conical, there is a $\C^*$-action on $Y$ such that $\pi$ is $\C^*$-equivariant. 
\end{enumerate}

\end{prop}

\begin{rem}

The proof of part (1) is based primarily on \cite[Corollary 1.4.3]{BCHM}. Additionally, part (3) was first proven in \cite[Proposition A.7]{NaFlops}.

\end{rem}

From now on, let $X$ be a conical affine symplectic singularity, and $\pi: Y \rightarrow X$ a $\Q$-factorial terminalization. Let $U \subset X$ be the complement of the symplectic leaves of codimension four or greater. This inherits the structure of a Poisson scheme from $X$, and since $\Q$-factorial terminalizations do not have any codimension two leaves, $\tilde{U}:= \pi^{-1}(U)$ is a symplectic resolution of $U$. In particular, is it the unique minimal resolution of $U$ (see \cite[p. 52]{NaPD1}, \cite[Section 3.7]{MMY25}).

Given a codimension two symplectic leaf $\mathcal{L}_i \subset U$ and a point $x_i \in \mathcal{L}_i$, we can take a transverse slice $S_i$ to $x_i$, which will be a Kleinian singularity.  Let $\hat{D}_i$ be the Dynkin diagram corresponding to this singularity, with Cartan $\hat{\mathfrak{h}_i}$ and Weyl group $\hat{W_i}$. Thus $\pi^{-1}(x_i)$ is a union of projective lines in the shape of $\hat{D}_i$, and we obtain a monodromy action of $\pi_1(\mathcal{L}_i, x_i)$ on $\hat{D}_i$ by graph automorphism.  This induces an action of $\pi_1(\mathcal{L}_i, x_i)$ on $\hat{\mathfrak{h}_i}$, and we define 
\begin{equation*}
\mathfrak{h}_i := \hat{\mathfrak{h}}_i^{\pi_1(\mathcal{L}_i, x_i)}  \hspace{20pt} \text{and} \hspace{20pt} W_i = \{w \in \hat{W_i} \ | \ w(\mathfrak{h}_i) = \mathfrak{h}_i\}.
\end{equation*} 
In particular, these are the Cartan and Weyl group of $D_i := \hat{D}_i/{\pi_1(\mathcal{L}_i, x_i)}$, as defined in \Cref{Fold}. 

\begin{rem}

Although it is not always possible to take an algebraic transverse slice to each codimension two symplectic leaf, we can always take either a formal slice or a transverse slice in the complex analytic category. 

\end{rem}

\begin{df}

Let $\{\mathcal{L}_i\}_{i \in I}$ be the set of codimension two symplectic leaves of $X$. Then the Namikawa Weyl group of $X$ is 
\begin{equation*}
W_X:= \prod_{i \in I} W_i
\end{equation*}
\end{df}

This group appears naturally when studying the Poisson deformations of $X$ and $Y$, as explained in \Cref{Intro}.

\section{Partial Symplectic Resolutions and Birational Geometry} \label[section]{Bir}

\subsection{Birational Geometry Preliminaries} \label[subsection]{BirPrelim}

In this subsection we will recall some statements about birational geometry and (relative) Mori Dream Spaces that will be relevant for the rest of the section. For more details, see \cite{HK} and \cite{Ohta}.

For the rest of this subsection let $\pi: Y \to X$ be a projective map of quasi-projective normal varieties over $\C$. Then we have the Abelian groups:
$$
\Pic(Y/X) := \Pic(Y) / \pi^*(\Pic(X)) 
$$
$$
Z_1(Y/X) := \Z[\{ \text{Integral curves } C \subset Y \text{ contracted by } \pi \}].
$$

We have a pairing
$$
\Pic(Y/X) \times Z_1(Y/X) \rightarrow \Z
$$
given by $L \cdot C := \deg(f^*L)$ where $f: \ol{C} \to C \to Y$ is the normalization of $C$. 

This pairing naturally extends to a pairing between $\Pic(Y/X)_{\Q} := \Pic(Y/X) \otimes \Q$ and $Z_1(Y/X)_{\Q} := Z_1(Y/X) \otimes \Q$. We say two elements $L_1, L_2 \in \Pic(Y/X)_{\Q}$ are numerically equivalent if
$$
L_1 \cdot C = L_2 \cdot C \hspace{10pt} \forall C \in Z_1(Y/X)_{\Q}
$$
and similarly two elements $C_1, C_2 \in Z_1(Y/X)_{\Q}$ are numerically equivalent if
$$
C_1 \cdot L = C_2 \cdot L \hspace{10pt} \forall L \in \Pic(Y/X)_{\Q}.
$$

Thus, we define 
$$
\cN^1(Y/X) := \{\Pic(Y/X)_{\Q}/\equiv\}, \hspace{10pt} \cN_1(Y/X) := \{Z_1(Y/X)_{\Q}/\equiv\},
$$
and there is a natural pairing $\cN^1(Y/X) \times \cN_1(Y/X) \to \Q$. We now introduce the definition of a partial resolution that we will be using for this paper.

\begin{df}

A map $\pi: Y \rightarrow X$ of normal varieties is a partial resolution if $\pi$ is projective and birational and the restriction of $\pi$ to $\pi^{-1}(X_{\textrm{reg}})$ is an isomorphism. 

\end{df}

\begin{rem}

Traditionally a (partial) resolution of singularities is only required to be proper, rather than projective. However, for our purposes we will only be interested in projective (partial) resolutions.

\end{rem}

Before proceeding we recall that Zariski's main theorem (see \cite[Corollary III.11.4]{Ha}) implies that for any partial resolution $\pi: Y \rightarrow X$, we have $\pi_*\cO_Y = \cO_X$, i.e. $\pi$ is a contraction. Given any contraction $f: Y \to Z$ over $X$, we have natural pushforward and pullback maps:
$$
f^*: \cN^1(Z/X) \to \cN^1(Y/X), \hspace{10pt} f_*: \cN_1(Y/X) \to \cN_1(Z/X)
$$

We now introduce some cones of interest in $\cN^1(Y/X)$

\begin{df}

The nef cone $\Nef(Y/X) \subset \cN^1(Y/X)$ is the cone of classes $[L] \in \cN^1(Y/X)$ that have nonnegative pairing with every element of $\cN_1(Y/X)$. A line bundle $L \in \Pic(Y)$ is $\pi$-movable if the base locus $\Cok[\pi^* \pi_* L \to L]$ has codimension at least 2. The movable cone $\Mov(Y/X) \subset \cN^1(Y/X)$ is the closure of the convex cone generated by classes $[L] \in \cN^1(Y/X)$ such that $L$ is $\pi$-movable. 

\end{df}

We now introduce the definition of a relative Mori Dream Space. Mori Dream Spaces were first introduced in \cite{HK}, and the definition is extended to the relative setting in \cite{Ohta}.

\begin{df}[{\cite[Definition 3.1]{Ohta}}]

Let $\pi: Y \to X$ be a projective map of quasi-projective normal varieties over $\C$. Then $\pi$ is a relative Mori Dream Space if 

\begin{enumerate}
\item $Y$ is $\Q$-factorial.
\item The natural surjection $\Pic(Y/X)_{\Q} \to \cN^1(Y/X)$ is an isomorphism.
\item $\Nef(Y/X)$ is a polyhedral cone generated by finitely many $\pi$-semiample line bundles. 
\item There exist finitely many small $\Q$-factorial modifications (SQMs) $f_i: Y \dashrightarrow Y_i$ over $X$ such that every $Y_i$ satisfies the previous condition, and 
$$
\Mov(Y/X) = \bigcup_{i} f_i^*(\Nef(Y_i/X)).
$$
\end{enumerate}

\end{df}

For the remainder of the paper, we will treat $\Mov(Y/X)$ as a fan, where the faces are given the faces of the polyhedral cones $f_i^*(\Nef(Y_i/X))$ described in the definition above. We will sometimes refer to this fan as the Mori fan. The faces of this fan are described using birational geometry in the following theorem, which is the natural generalization of \cite[Proposition 1.11 (3)]{HK}:

\begin{thm}[{\cite[Theorem 4.4]{Ohta}}] \label[thm]{Ohta4.4}

There is a natural bijection between the set of faces of the fan $\Mov(Y/X)$, and the set of rational contractions $f: Y \dashrightarrow Z$, where $Z$ is projective over $X$. This bijection is given by 
$$
f \mapsto f^*(\Nef(Z/X)).
$$

\end{thm}

Note that a rational contraction is a contraction morphism $f: Y \to Z$ if and only if it corresponds to a face of $\Nef(Y/X)$ in the above theorem. 

\subsection{Symplectic Partial Resolutions} \label{BirSymp}

The following definition is a natural extension of the definition of a symplectic resolution.

\begin{df}

A partial resolution $\pi: Y \rightarrow X$ is symplectic if the symplectic form $\omega_{X_{\textrm{reg}}}$ extends to a symplectic form on $Y_{\textrm{reg}}$.

\end{df}

\begin{lemma} \label[lemma]{PartialSympl}

Any symplectic partial resolution $Y$ of a symplectic singularity $X$ is itself a symplectic singularity.

\end{lemma}

\begin{proof}

By definition $Y$ is normal and has a symplectic form $\omega_{Y_{\textrm{reg}}}$ on its smooth locus that extends $\omega_{X_{\textrm{reg}}}$.  All that is left to check is that given a resolution $\tilde{\pi}: Z \rightarrow Y$, $\tilde{\pi}^{*}(\omega_{Y_{\textrm{reg}}})$ extends to a regular 2-form on $Z$.  Since $\pi \circ \tilde{\pi}$ is a resolution of $X$, and $X$ is a symplectic singularity, $\omega_{X_{\textrm{reg}}}$ extends to a regular 2-form $\omega_Z$ on $Z$. Note that since $Z$ is smooth, the sheaf of 2-forms is a vector bundle, and thus restriction to a Zariski open set is injective. Since we can treat $X_{\textrm{reg}}$ and $Y_{\textrm{reg}}$ as Zariski open sets of $Z$, this means that the extensions $\omega_{Y_{\textrm{reg}}}$ and $\omega_Z$ of $\omega_{X_{\textrm{reg}}}$ are unique.  In particular, $\omega_Z$ must extend $\omega_{Y_{\textrm{reg}}}$.  
\end{proof}

The following argument is identical for full resolutions (see \cite[Proposition 1.1]{Fu}), but we include it for completeness.

\begin{lemma} 

A partial resolution is symplectic if and only if it is crepant.

\end{lemma}

\begin{proof}

Let $X$ be a symplectic singularity with $\dim X = 2n$ and let $\omega_X$ be the symplectic form on $X_{\textrm{reg}}$. Then $\wedge^n \omega_X$ provides a trivialization of $K_X$ over $X_{\reg}$ which extends to a trivialization of $K_X$ by normality. Thus, every symplectic singularity has trivial canonical sheaf. In particular, given a partial symplectic resolution $\pi: Y \rightarrow X$, we have $\pi^*(K_X) = \pi^*0 = 0 = K_Y$, so the map is crepant. 

Conversely, say $\pi: Y \rightarrow X$ is crepant. Since $X$ has symplectic singularities, $\pi^{*}\omega_X$ extends to a regular 2-form $\omega_Y$ on $Y_{\textrm{reg}}$. Additionally, $\wedge^n \omega_X$ extends to a non-vanishing global section $f$ of $K_X \cong \cO_X$. Thus, $f$ pulls back to a non-vanishing section of $K_Y \cong \cO_Y$ that extends $\pi^*(\wedge^n \omega_X) = \wedge^n \pi^*\omega_X$ and thus $\wedge^n \omega_Y$. Therefore, $\wedge^n \omega_Y$ is a non-vanishing section of $K_Y|_{Y_{\reg}} \cong \Omega^n_{Y_{\reg}}$, which implies that $\omega_Y$ is non-degenerate. Finally, since $\pi^*$ commutes with differentials, $\diff(\omega_Y)$ extends $\pi^*\diff(\omega_X) = 0$ on $\pi^{-1}(X_{\reg})$, which extends uniquely to 0 on $Y_{\reg}$. Thus, $\omega_Y$ is closed, so $\pi$ is symplectic. 
\end{proof}

The above theorem implies that if $\pi: Y \rightarrow X$ is a $\Q$-factorial terminalization of a symplectic singularity $X$, then $Y$ is a symplectic partial resolution of $X$. In fact, $\Q$-factorial terminalizations can be described as the maximal elements of the poset of symplectic partial resolutions, as shown below.

\begin{lemma} \label[lemma]{Cover}

A partial resolution $\rho: X' \rightarrow X$ of a symplectic singularity $X$ is symplectic if and only if it is covered by a $\Q$-factorial terminalization of $X$.

\end{lemma}

\begin{proof}

If $\rho: X' \rightarrow X$ is symplectic then by \Cref{PartialSympl} it has symplectic singularities and thus by \Cref{QFac} there is some $\Q$-factorial terminalization $\pi': Y \rightarrow X'$. Thus, $\pi := \rho \circ \pi'$ is a $\Q$-factorial terminalization of $X$ that covers $\rho$.

Conversely, let $\pi: Y \rightarrow X$ be a $\Q$-factorial terminalization of $X$ that covers $\rho$, with covering map $\pi': Y \rightarrow X'$. Then the symplectic form $\omega_{X_{\textrm{reg}}}$ on $X_{\reg}$ extends to a symplectic form on $Y_{\reg}$. In particular, $\pi = \rho \circ \pi'$ implies that $X'_{\reg} \subset Y_{\reg}$, so $\omega_{X_{\textrm{reg}}}$ extends to a symplectic form on $X'_{\reg}$. Thus, $\rho$ is symplectic.
\end{proof}

For the remainder of the paper we will only be interested in symplectic partial resolutions, so we will simply call them partial resolutions and they will be assumed to be symplectic. 

Additionally, for the remainder of the paper let $X$ be a conical affine symplectic singularity, and let $\rho: X' \rightarrow X$ be a partial resolution that is covered by a $\Q$-factorial terminalization $\pi: Y \rightarrow X$. Since $X$ is affine we have $H^i(X, \cO_X) = 0, \forall i > 0$. Furthermore, by \cite[Proposition 1.3]{Beau}, symplectic singularities have rational singularities. Thus $X$, $X'$, and $Y$ all have rational singularities, so we have quasi-isomorphisms $\pi_*\cO_Y \cong \cO_X$ and $\rho_*\cO_{X'} \cong \cO_X$ in the derived category, which implies that 
$$
H^i(Y, \cO_{Y}) = H^i(X', \cO_{X'}) = 0, \forall i > 0.
$$

The following theorem will allow us pass back and forth between the settings of birational geometry and Poisson deformation theory, as will become clearer in the following section.

\begin{lemma}[{\cite[Lemma 4.23]{LMM}}]

The first Chern class map provides an isomorphism
\begin{equation*}
c_1: \Pic(Y^{\reg}) \otimes_{\Z} \Q \rightarrow H^2(Y^{\reg}, \Q).
\end{equation*}

\end{lemma}

\begin{rem}

The analogue of the above theorem for the analytic Picard group can be deduced via a standard computation using the exponential short exact sequence. However, the above statement applies to the algebraic Picard group, which requires an argument using the GAGA principle.

\end{rem}

As detailed in \cite{NaPD2}, $W_X$ acts on $H^2(Y^{\reg}, \C)$ as a finite reflection group, and the action is generically free. The following celebrated theorem of Namikawa will allow us to classify partial resolutions of conical affine symplectic singularities.

\begin{thm}[{\cite[Main Theorem, Remark 10]{NaMori}}] \label[thm]{NamikawaMDS}

If $\pi: Y \rightarrow X$ is a $\Q$-factorial terminalization, then there are finitely many hyperplanes $\{H_i\} \subset H^2(Y^{\reg}, \R)$ such that

\begin{enumerate}
\item The union $\mathcal{D}$ of the $H_i$ is invariant under the action of $W_X$.
\item $\mathcal{D}$ is the locus over which the universal deformation map $\mathcal{Y} \rightarrow \mathcal{X}$ is not a fiberwise isomorphism, as described in diagram \eqref{IntroUniv}.
\item For every $\Q$-factorial terminalization $\pi': Y' \rightarrow X$, $\ol{\Amp(\pi')} = \Nef(\pi')$ is a unique chamber cut out by the hyperplanes.
\item The union of $\Nef(\pi')$ over every $\Q$-factorial terminalization $\pi': Y' \rightarrow X$ is $\Mov(Y/X)$.
\end{enumerate}
In particular, $\pi: Y \rightarrow X$ is a relative Mori Dream Space.

\end{thm}

Note that traditionally the cones in \Cref{NamikawaMDS} below would be described as subsets of 
\begin{equation*}
\Pic(Y/X) \otimes \R := (\Pic(Y)/\pi^*(\Pic(X))) \otimes \R. 
\end{equation*}
However, by \cite[Lemma 7.1]{LMM}, the Picard group of every normal affine cone is trivial, and thus $\Pic(X) = 0$. Additionally, by \cite[Proposition III.6.5]{Ha}, $\Cl(Y) \cong \Cl(Y^{\reg})$, where $\Cl$ denotes the divisor class group, since $Y - Y^{\reg}$ has codimension at least four. Finally, since $Y$ and $Y^{\reg}$ are $\Q$-factorial, this implies $\Pic(Y) \otimes \Q = \Pic(Y^{\reg}) \otimes \Q$. So indeed
\begin{equation*}
\Pic(Y/X) \otimes \R \cong H^2(Y^{\reg}, \R).
\end{equation*}

\begin{rem}

As emphasized by Namikawa in \cite[p. 341]{NaMori}, a crepant resolution $Y \to X$ is a relative Mori Dream space whenever $X$ has rational Gorenstein singularities, but in the case of a conical affine symplectic singularity and a $\Q$-factorial terminalization we have the added feature that the Mori fan is given by full hyperplanes. 

\end{rem}

\begin{cor} \label[cor]{face=res}

Partial resolutions of $X$ are in bijection with faces of the Mori fan of $\Mov(Y/X)$.

\end{cor}

\begin{proof}

By \Cref{Ohta4.4}, there is a bijection between faces of $\Mov(Y/X)$ and birational contractions $Y \dashrightarrow X'$, where $X'$ is projective over $X$. Each of these birational contractions is a contraction of some $\Q$-factorial terminalization of $X$, and thus by \Cref{Cover} it is a partial resolution of $X$. Conversely, given a partial resolution $\rho: X' \rightarrow X$ of $X$, it is covered by some $\Q$-factorial terminalization of $X$ via a map $\pi': Y' \rightarrow X'$. This map is a contraction, so there is a face of $\Mov(Y/X)$ that corresponds to $\rho: X' \rightarrow X$.
\end{proof}

Because $\pi: Y \rightarrow X$ is a relative Mori Dream Space, $\Nef(Y/X)$ is a closed polyhedral cone generated by semiample line bundles. The cone of effective curves 
$$
\NE_1(Y/X) \subset \cN_1(Y/X)
$$
is dual to $\Nef(Y/X)$, and thus it is a polyhedral cone generated by finitely many curve classes $[C_k]$.  Thus every codimension one face $F$ of $\Nef(Y/X)$ is of the form $\Nef(Y/X) \cap [C_k]^{\perp}$ for a unique extremal ray, $\langle [C_k] \rangle \subset \NE_1(Y/X)$. We can recover this extremal ray as the kernel of $\pi'_*: \cN_1(Y/X) \rightarrow \cN_1(X'/X)$, where $\pi':Y \rightarrow X'$ is the contraction corresponding to $F$.

We follow \cite[Definition 3.4]{AW} with the following definition

\begin{df}

Given a $\Q$-factorial terminalization $Y \rightarrow X$, the cone of essential curves, $\Ess(Y/X)$ is the cone generated by curve classes $[C_i] \in \mathcal{N}_1(Y/X)$ of contracted curves $C_i$ whose image intersects a codimension two symplectic leaf of $X$.

\end{df}

As explained in \Cref{Background}, given any codimension two leaf $\mathcal{L}_i \subset X$ and $x_i \in \mathcal{L}_i$, $\pi^{-1}(x_i)$ is given by a union of projective lines $C_k$ in the shape on the Dynkin diagram $\hat{D}_i$ corresponding to the transverse slice to $\cL_i$. Since these curve classes are all contracted by $\pi$, they also represent rays of $\NE_1(Y/X)$, and in particular they generate $\Ess(Y/X) \subset \NE_1(Y/X)$. 

Given a simple generator $s_k$ of $W_X$ corresponding to a node of $D_i := \hat{D}_i/\pi_1(\mathcal{L}_i, x_i)$, we can take a projective line $C_k$ that maps to $x_i \in \cL_i$ and that is in the orbit of the monodromy action of $\pi_1(\mathcal{L}_i, x_i)$ that corresponds to the same node of $D_i$ as $s_k$. It follows from \cite[Proposition 4.4]{Wie} that two of the above projective lines $C_k$ and $C_l$ are in the same curve class if they map to the same symplectic leaf and they are in the same orbit of the $\pi_1(\mathcal{L}_i, x_i)$-action. Thus $s_k \mapsto [C_k]$ is a well-defined map
\begin{equation*}
\Phi: \{\text{Simple generators of }W_X\} \rightarrow \{\text{Generators of }\Ess(Y/X)\}.
\end{equation*}

This map is surjective by construction since for every projective line $C_k$ that maps to a codimension two leaf $\cL_i$, we can take the node of $D_i$ corresponding to its orbit under the monodromy action of $\pi_1(\mathcal{L}_i, x_i)$, and the corresponding simple generator $s_k$ will map to $[C_k]$. The fact that it is a bijection will follow from the work below. 

Given a simple generator $s_k$ of the Namikawa Weyl group $W_X$, let 
$$
H_{s_k} \subset H^2(Y^{\reg}, \R) \cong \Pic(Y/X) \otimes \R
$$
be the hyperplane $\Phi(s_k)^{\perp}$. 

\begin{lemma} \label[lemma]{Hgen}

The hyperplanes $H_{s_i}$ are in general position in $H^2(Y^{\reg}, \R)$. 

\end{lemma}

\begin{proof}
This follows from \cite[Proposition 2.17]{BPW}. Although in that paper it is assumed that $Y$ is a symplectic resolution, the same argument applies, which we shall briefly explain below. 

We will show the natural restriction map
\begin{equation} \label{BPWmap}
H^2(Y^{\reg}, \R) \rightarrow \bigoplus_{i \in I} H^2(\pi^{-1}(x_i), \R)^{\pi_1(\mathcal{L}_i, x_i)}
\end{equation}
is surjective. It follows from \cite[p. 81]{NaPD1} that the restriction map $H^2(Y^{\reg}, \R) \rightarrow H^2(\pi^{-1}(U), \R)$ is an isomorphism, where $U$ is the union of the open and codimension two symplectic leaves of $X$. To show that the restriction map 
\begin{equation} \label{BPWmap2}
H^2(\pi^{-1}(U), \R) \rightarrow \bigoplus_{i \in I} H^2(\pi^{-1}(x_i), \R)^{\pi_1(\mathcal{L}_i, x_i)}
\end{equation}
is surjective we can apply \cite[Section 4]{NaPD1}, which explains that the irreducible components of the exceptional divisor of $\pi^{-1}(U) \rightarrow U$ are in bijection with the nodes of the $D_i$ and that two projective lines in a fiber $\pi^{-1}(x_i)$ are contained in the same irreducible component if and only if they are in the same $\pi_1(\mathcal{L}_i, x_i)$-orbit. We can see that the classes in $H^2(\pi^{-1}(U), \R)$ of these divisors are linearly independent because their classes form a basis for $H^0(U, R^2\pi_*\C)$, as noted in the first sentence of \cite[p. 78 (ii)]{NaPD1}. Thus, the subspace of $H^2(\pi^{-1}(U), \R)$ generated by these irreducible components is mapped isomorphically into the image by \eqref{BPWmap2}.

Let $s$ be a simple generator of $W_X$ corresponding to a ${\pi_1(\mathcal{L}_i, x_i)}$-orbit of curves in $\pi^{-1}(x_i)$. Then $H_s$ pulls back under \eqref{BPWmap} to the hyperplane in $H^2(\pi^{-1}(x_i), \R)^{\pi_1(\mathcal{L}_i, x_i)}$ perpendicular to the sum of the fundamental classes of the curves in the ${\pi_1(\mathcal{L}_i, x_i)}$-orbit corresponding to $s$. This implies that the intersection of all the $H_s$ is contained in the kernel of \eqref{BPWmap}. Since the number of hyperplanes is equal to the dimension of the codomain of \eqref{BPWmap}, the codimension of the intersection of all the hyperplanes is the same as the number of hyperplanes, so they are in general position. 
\end{proof}

It follows from the above lemma that the $[\Phi(s_k)]$ are not only distinct but linearly independent in $\cN_1(Y/X)$. Additionally, by \cite[Theorem 3.5]{AW}, $\Ess(Y/X)$ is dual to $\Mov(Y/X)$ (see also the argument in \cite[Proposition 2.17]{BPW}). Since $[\Phi(s_k)]$ is a generator of $\Ess(Y/X)$ and $\Ess(Y/X)$ is dual to $\Mov(Y/X)$, $H_{s_k}$ bounds $\Mov(Y/X)$ and thus is one of the hyperplanes in \Cref{NamikawaMDS}. Additionally, if $\pi': Y \rightarrow X'$ is a contraction corresponding to a face $F \subset \Nef(Y/X)$, then by definition $\pi'$ contracts $[\Phi(s_k)]$ if and only if $F \subset H_{s_k}$.

\begin{rem}

As explained in the proof of \cite[Proposition 2.17]{BPW}, it follows from the above work that $\Mov(Y/X)$ is the preimage under the map \eqref{BPWmap} of the chamber that has non-negative pairing with the fundamental classes of all the irreducible components of each $\pi^{-1}(x_i)$. In particular, this implies that $\Mov(Y/X)$ is always a direct sum of a Weyl chamber of a semisimple Lie algebra with a product of real lines. The relationship between Kleinian singularities and Lie theory will be further explored in \Cref{Klein}, where we will also show that each $H_s$ is the subspace of $H^2(Y^{\reg}, \R)$ fixed by $s$.

\end{rem}

\begin{df} \label[df]{PartWeyl}

Let $\rho: X' \rightarrow X$ be the partial resolution of $X$ corresponding to a face $F \subset \Mov(Y/X)$. We define the Namikawa Weyl group $W_{X'}$ of this partial resolution to be the parabolic subgroup $W_{X'} \subset W_X$ generated by all the $s_k$ such that $F \subset H_{s_k}$.

\end{df}

\begin{rem}

The Namikawa Weyl group has previously only been defined for conical affine symplectic singularities, although the above definition has similar motivation. We shall show later that it satisfies a similar condition as the Namikawa Weyl group $W_X$ stemming from Poisson deformation theory.

\end{rem}

\begin{lemma}

The above construction does not depend on a choice of $\Q$-factorial terminalization. 

\end{lemma}

\begin{proof}

Let $U \subset X$ be the union of the open and codimension two symplectic leaves. Then any two $\Q$-factorial terminalizations $\pi_1: Y_1 \rightarrow X$ and $\pi_2: Y_2 \rightarrow X$ restrict to the minimal resolution $\tilde{U} \to U$ over $U$ (see \cite[p. 52]{NaPD1}, \cite[Section 3.7]{MMY25}). Thus, the essential cones $\Ess(Y_1/X)$ and $\Ess(Y_2/X)$ are canonically identified.
\end{proof}

\begin{rem}

A different way of seeing that every $\Q$-factorial terminalization of $X$ is isomorphic over $U$ is to note that any two $\Q$-factorial terminalizations are related by a finite sequence of flops. In particular, given one such flop $Y_1 \rightarrow X' \leftarrow Y_2$, the contractions $Y_1 \rightarrow X'$ and $Y_2 \rightarrow X'$ are both small. Thus, they cannot contract any of the curve classes that maps to codimension two leaves in $X$, which would constitute a divisorial contraction. So the flop $Y_1 \dashrightarrow Y_2$ restricts to an isomorphism $\pi_1^{-1}(U) \xrightarrow{\sim} \pi_2^{-1}(U)$ over $U$.

\end{rem} 

\begin{thm} \label[thm]{AllSubgroups}

The map 
\begin{equation*}
\Psi: \{\text{Partial resolutions of }X\} \rightarrow \{\text{Parabolic subgroups of }W_X\}
\end{equation*}
given by $[\rho: X' \rightarrow X] \mapsto W_{X'}$ is surjective, and it is bijective if and only if $X$ has a unique $\Q$-factorial terminalization. 

\end{thm}

\begin{proof}

Let $W' \subset W_X$ be the parabolic subgroup of $W_X$ generated by $s_{1}, \dots , s_{l}$. By \Cref{Hgen}, the hyperplanes $H_s$ are in general position, and in particular $H_{s_{1}} \cap \dots \cap H_{s_{l}}$ has codimension $l$. It also bounds $\Mov(Y/X)$, so we can pick some codimension $l$ face $F'$ of $\Mov(Y/X)$ that is open in $H_{s_{1}} \cap \dots \cap H_{s_{l}}$. Since $F'$ is codimension $l$, it cannot be contained in any of the other hyperplanes $H_s$. By \Cref{face=res} and \Cref{PartWeyl}, the partial resolution corresponding to $F'$ has Namikawa Weyl group $W'$.

If $X$ has a unique $\Q$-factorial terminalization, then the Mori fan of $\Mov(Y/X)$ is simply given by the faces that bound $\Mov(Y/X) = \Nef(Y/X)$. The hyperplanes that bound $\Mov(Y/X)$ are all hyperplanes $H_s$, and thus by \Cref{Hgen} they are in general position. Thus the faces of the Mori fan are each of the form $\Nef(Y/X) \cap H_{s_{1}} \cap \dots \cap H_{s_{i}}$, for some unique subset $\{s_{1}, \dots , s_{i}\}$ of the simple roots, so $\Psi$ is a bijection. However, if there is more than one $\Q$-factorial terminalization, they both have trivial Namikawa Weyl group, so $\Psi$ is not a bijection.
\end{proof}

The following two corollaries are not strictly needed but they explain how the constructions in this section can be viewed as generalizations of Springer theory. 

\begin{cor} \label[cor]{NilpCovers}

Every partial resolution of the nilpotent cone $\cN$ of a semisimple Lie algebra $\mathfrak{g}$ is of the form $\rho: \tilde{\cN}^{\mathcal{P}} \rightarrow \cN$ with
\begin{equation*}
\tilde{\cN}^{\mathcal{P}} := \{(x, \mathfrak{p}) | x \in \mathfrak{p} \cap \cN, \mathfrak{p} \in \mathcal{P} \}, \hspace{10pt} \rho: (x, \mathfrak{p}) \mapsto x
\end{equation*}
for some conjugacy class ${\mathcal{P}}$ of parabolic subalgebras. 
\end{cor}

\begin{proof}

Each of these partial resolutions is covered by the Springer resolution $\mu: \tilde{\cN} \rightarrow \cN$ because every Borel subalgebra is contained in a unique element of ${\mathcal{P}}$. Thus, by \Cref{Cover}, all of these partial resolutions are symplectic. Additionally, they are in bijection with subsets of the simple roots of $\mathfrak{g}$, which are also in bijection with parabolic subgroups of the Weyl group of $\mathfrak{g}$. Furthermore, the Weyl group of $\mathfrak{g}$ is the Namikawa Weyl group of $\cN$. Since the Springer resolution is the unique symplectic resolution (or $\Q$-factorial terminalization) of $\cN$, by \Cref{AllSubgroups}, these must be all the partial resolutions.
\end{proof}

Recall from \eqref{FibersInv} we can identify a canonical parabolic subgroup $W^\mathcal{P} \subset W$ by taking the Weyl group of Levi part of any $\mathfrak{p} \in \mathcal{P}$. 

\begin{cor} \label[cor]{WeylLevi}

The Namikawa Weyl group of the partial resolution $\rho: \tilde{\cN}^{\mathcal{P}} \rightarrow \cN$ is $W^\mathcal{P}$.

\end{cor}

\begin{proof} 

Recall that $\cN$ has a unique codimension two symplectic leaf given by the subregular orbit. Furthermore, by \Cref{SubregFiber} the essential curve classes of the Springer resolution $\pi: \tilde{\cN} \rightarrow \cN$ are in canonical bijection with the nodes of the Dynkin diagram $D$ of $\mathfrak{g}$. In particular, if $\alpha$ is a root of $\mathfrak{g}$, then the corresponding essential curve $[C_{\alpha}]$ of $\tilde{\cN}$ is represented by $\{(x, \mathfrak{b}) \ | \ \mathfrak{b} \subset \mathfrak{p}\}$, where $\mathfrak{p}$ is any minimal parabolic of type $\alpha$ that contains $x$ in its nilradical. The map $\pi': \tilde{\cN} \rightarrow \tilde{\cN}^{\mathcal{P}}$ is given by $(x, \mathfrak{b}) \mapsto (x, \mathfrak{p})$, where $\mathfrak{p}$ is the unique element of $\mathcal{P}$ containing $\mathfrak{b}$. Thus, if $\mathcal{P}$ corresponds to a subset $I$ of the simple roots, $\pi'$ contracts $[C_{\alpha}]$ if and only if $\alpha \in I$. Therefore, the Namikawa Weyl group of $\rho: \tilde{\cN}^{\mathcal{P}} \rightarrow \cN$ is the parabolic subgroup of $W$ corresponding to $I$, which is precisely $W^\mathcal{P}$.
\end{proof}

\begin{Ex}

The following example is taken from the introduction of arXiv v7 of \cite{NaBir}. 

Let $X$ be the closure of the subregular orbit of the nilpotent cone of $\mathfrak{sl}_4$. Note that by \cite[Remark 1]{He}, $X$ is normal. Furthermore, let $\{\alpha_1, \alpha_2, \alpha_3\}$ be the simple roots of $\mathfrak{sl}_4$, so that $\alpha_2$ corresponds to the middle node of the Dynkin diagram. Then $X$ has three symplectic resolutions, given by the Springer maps $T^*(G/P_{\alpha_i}) \rightarrow X$ for $P_{\alpha_1}$, $P_{\alpha_2}$, and $P_{\alpha_3}$ elements of the three conjugacy classes of minimal parabolics of $SL_4$. The Mori fan of $\Mov(T^*(G/P_{\alpha_i})/X)$ is below. 
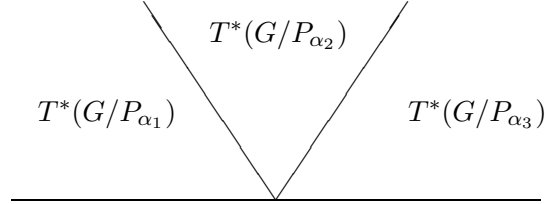
\begin{equation*}
\begin{picture}(300,100)(0,0) 
\put(150,0){\line(2,3){50}}\put(150,0){\line(1,0){100}}
\put(150,0){\line(-1,0){100}}\put(150,0){\line(-2,3){50}}
\put(125,60){$T^*(G/P_{\alpha_2})$}\put(200,30){$T^*(G/P_{\alpha_3})$}
\put(60,30){$T^*(G/P_{\alpha_1})$}
\end{picture} 
\end{equation*}

The Namikawa Weyl group is $\Z/2\Z$ and it acts by reflection across the bottom hyperplane. This hyperplane is also $H_s$, for $s$ the unique simple generator of the Namikawa Weyl group. Thus, the three partial resolutions corresponding to faces which are contained in $H_s$ have Namikawa Weyl group $\Z/2\Z$, while the other five have trivial Namikawa Weyl group.

\end{Ex}

\begin{Ex}

Another interesting family of examples is given by the symmetric powers $\Sym^n(S)$ of Kleinian singularities $S$, which each have a symplectic resolution given by $\Hilb^n(\tilde{S})$, where $\tilde{S}$ is the unique symplectic resolution of $S$. These can be constructed as quiver varieties, and are studied extensively from this perspective in \cite{BC}. In particular, \cite[Proposition 2.2]{BC} tells us that they have Namikawa Weyl group $\Z/2\Z \times W_S$, where $W_S$ is the Namikawa Weyl group corresponding to the Kleinian singularity $S$.

Treating $\C^2$ as the trivial Kleinian singularity, $\pi: \Hilb^n(\C^2) \rightarrow \Sym^n(\C^2)$ is the unique symplectic resolution of $\Sym^n(\C^2)$. Additionally, $H^2(\Hilb^n(\C^2), \R)$ is one dimensional. Thus the Namikawa Weyl group $\Z/2\Z$ acts by reflection across 0, and $\Nef(\Hilb^n(\C^2)/\Sym^n(\C^2)) = \Mov(\Hilb^n(\C^2)/\Sym^n(\C^2))$ is one of the two rays. This implies that $\pi$ and the identity are the only partial resolutions of $\Sym^n(\C^2)$. In particular, this is a family of examples of arbitrarily large dimension and with arbitrarily many symplectic leaves (given by partitions of $n$) whose only partial resolutions are the identity and a single symplectic resolution each. 

On the other hand, for any non-trivial Kleinian singularity $S$, the number of symplectic resolutions of $\Sym^n(S)$ grows with $n$. For example, the examples below, which appear in \cite[Example 2.6]{CGGS} and \cite[Example 2.6]{BC} respectively, are cross sections of the three-dimensional Mori fans for $\Sym^3(\C^2/(\Z/3\Z))$ and $\Sym^4(\C^2/(\Z/3\Z))$.
\begin{equation*}
\begin{tikzpicture}[baseline={(0,-0.3)},xscale=1.75,yscale=1.5]
			\tikzset{>=latex}
            \draw (4,4)--(2,0)--(0,4)
            --(4,4)
            --(10/11,24/11)
            --(34/11,24/11)
            --(0,4)
            --(2.75,1.5)
            --(1.25,1.5)
            --(4,4)
            --(0,4)
            -- cycle;
               \node at (2,3.7) {$\Hilb^4(\C^2/(\Z/3\Z))$};
   \node at (2,4.2) {$H_0$};
  \node at (0.7,2) {$H_1$};
 \node at (3.3,2) {$H_2$};
\end{tikzpicture}
\hspace{20pt}
\begin{tikzpicture}[baseline={(0,-0.3)},xscale=1.75,yscale=1.5]
			\tikzset{>=latex}
            \draw (2,0) -- (4,4);
			\draw (0,4) -- (4,4);
			\draw (0,4) -- (2,0);
            \draw (1.25,1.5) -- (4,4);
            \draw (2.75,1.5) -- (0,4);
            \draw (1.25,1.5) -- (2.75,1.5);
            \draw (10/11,24/11) -- (34/11,24/11);
            \draw (10/11,24/11) -- (4,4);
            \draw (0,4) -- (34/11,24/11);
            \draw (5/7,72/28) -- (23/7,72/28);
            \draw (5/7,72/28) -- (4,4);
            \draw (0,4) -- (23/7,72/28);
  \node at (2,3.7) {$\Hilb^4(\C^2/(\Z/3\Z))$};
   \node at (2,4.2) {$H_0$};
  \node at (0.7,2) {$H_1$};
 \node at (3.3,2) {$H_2$};
\end{tikzpicture}
\end{equation*}

In both of these cases $H_0$ is the hyperplane corresponding to the generator $s_0$ of the $\Z/2\Z$ summand of the Namikawa Weyl group $\Z/2\Z \times S_3$, while $H_1$ and $H_2$ correspond to the two simple generators $s_1$ and $s_2$ of $S_3$. Thus, in the first example, $\Sym^3(\C^2/(\Z/3\Z))$, there are 12 symplectic resolutions, 31 partial resolutions (including the 12 symplectic resolutions) with trivial Namikawa Weyl group, 5 each with Namikawa Weyl groups $\langle s_1 \rangle$ and $\langle s_2 \rangle$, and 1 each with Namikawa Weyl groups $\langle s_0 \rangle$, $\langle s_0, s_1 \rangle$, $\langle s_0, s_2 \rangle$, $\langle s_1, s_2 \rangle$, and the full subgroup $\Z/2\Z \times S_3$. Note that the identity is the unique partial resolution with Namikawa Weyl group $\Z/2\Z \times S_3$, and that the corresponding face is the origin, which is not pictured in the above image because the image is a two dimensional cross section. 

\end{Ex}

\section{Deformation Theory and Prorepresentability} \label[section]{Def}
In this section we study the Poisson deformation theory of partial resolutions of conical affine symplectic singularities.  Let $\Art_{\C}$ be the category of Artinian local rings with residue field $\C$. Given a Poisson scheme (or complex analytic space) $X$, we let 
$$
\PD_X: \Art_{\C} \rightarrow \text{Groupoids},
$$ 
be the standard the Poisson deformation $2$-functor. In particular, given any $A \in \Art_{\C}$, $\PD_X(A)$ is the groupoid where objects are flat Poisson schemes $X_A$ over $\Spec(A)$ such that $X_A \times_{\Spec(A)} \Spec(k) \cong X$ and morphisms are given by isomorphisms that restrict to the identity on $X$ over $\Spec(k)$. Finally, functors are given by base change. We will use the same notation to refer to the category cofibered in groupoids associated to this $2$-functor (see \cite[\href{https://stacks.math.columbia.edu/tag/0048}{Section 0048}]{Stacks}).

If we take the underlying functor to Sets, $\pd_X: \Art_{\C} \rightarrow \textrm{Sets}$ given by $\pd_X(A) = \textrm{Ob}(\PD_X(A))/\cong$, we recover the usual Poisson deformation functor, studied extensively by Namikawa in \cite{NaFlops, NaPD1, NaPD2, NaMori}. In particular, Namikawa proves in these papers that for a $\Q$-factorial terminalization $\pi: Y \rightarrow X$ of a conical affine symplectic singularity $X$, $\pd_X$ and $\pd_Y$ are both prorepresentable and unobstructed. Throughout this section we shall adapt many of the techniques and arguments in those papers to partial resolutions.

\subsection{Prorepresentability} \label[subsection]{ProDef}

In this subsection we will show that given any symplectic partial resolution $\rho: X' \to X$ of a conical affine symplectic singularity, $\pd_{X'}$ is prorepresentable. We shall do this by showing that $\pd_{X'}$ satisfies Schlessinger's criterion (\cite[Theorem 2.11]{Sch}). We recall this criterion below:

\begin{thm}[{\cite[Theorem 2.11]{Sch}}] \label[thm]{SchThm}

Let $F: \Art_{\C} \to \text{Sets}$ be a functor such that $F(\C)$ is a point. Given any two maps $A' \to A$ and $A'' \to A$ in $\Art_{\C}$, there is a natural map
$$
\text{res}: F(A' \times_{A} A'') \to F(A') \times_{F(A)} F(A'').
$$ 
The functor $F$ has a hull if it satisfies the first three conditions below, and it is prorepresentable if it satisfies all four:
\begin{enumerate}
\item \text{res} is a surjection whenever $A'' \to A$ is a small extension.
\item \text{res} is a bijection when $A = \C$ and $A'' = \C[\epsilon]$.
\item $F(\C[\epsilon])$ is a finite dimensional vector space.
\item \text{res} is a bijection when $A' \to A$ and $A'' \to A$ are the same small extension.
\end{enumerate}

\end{thm}

As explained in \cite[2.13]{Sch}, when condition (2) of \Cref{SchThm} is satisfied, the set $F(\C[\epsilon])$ has a natural vector space structure. This explains condition (3), and we will also treat $\pd_Z(\C[\epsilon])$ as a vector space for any symplectic singularity $Z$ after establishing the first two conditions of Schlessinger's criterion. In order to show that $\pd_{X'}$ satisfies the first two conditions of Schlessinger's criterion, we will show that the 2-functor $\PD_{X'}$ satisfies the Rim-Schlessinger condition. This condition was first introduced in \cite{Rim}, and the definition can be found in \cite[\href{https://stacks.math.columbia.edu/tag/06J2}{Definition 06J2}]{Stacks}. By \cite[\href{https://stacks.math.columbia.edu/tag/06J5}{Lemma 06J5}]{Stacks}, the Rim-Schlessinger condition for a 2-functor $F: \Art_{\C} \rightarrow \text{Groupoids}$ is satisfied if and only if the natural map of groupoids
$$
\text{res}: F(A' \times_{A} A'') \to F(A') \times_{F(A)} F(A'')
$$
is an equivalence for any maps $A' \to A$ and $A'' \to A$ in $\Art_{\C}$, where the first is a surjection. For more details on the Rim-Schlessinger condition, see \cite{TV}, especially Section 2.3.

The proof of \Cref{Rim} below is very similar to that of other deformations of algebraic objects, and we assume it is known to experts. However, we include it for the benefit of the reader, since as far as we are aware it is absent from other papers on symplectic singularities. Before the proof we recall the following key lemma:

\begin{lemma}[{\cite[\href{https://stacks.math.columbia.edu/tag/0D2I}{Lemma 0D2I}]{Stacks}}] \label[lemma]{FlatRim}

Assume we have the diagram of rings
\[\begin{tikzcd}
	{B'} & {A'} \\
	\
	{B} & {A}
	\arrow[from=1-1, to=1-2]
	\arrow[from=1-1, to=2-1]
	\arrow[from=1-2, to=2-2]
	\arrow[from=2-1, to=2-2]
\end{tikzcd}\]
with $A' \rightarrow A$ surjective and $B' = B \times_{A} A'$. Additionally, for a ring $R$, let $\Flat_R$ be the full subcategory of $\Mod_R$ consisting of flat modules. Then there is an equivalence of categories 
\begin{equation*}
\Flat_{B'} \rightarrow \Flat_{B} \times_{\Flat_{A}} \Flat_{A'}
\end{equation*}
given by $M \mapsto (M \otimes_{B'} B) \times_{(M \otimes_{B'} A)} (M \otimes_{B'} A')$.

\end{lemma}

The following proposition establishes that $\PD_X$ satisfies the Rim-Schlessinger condition for any symplectic singularity $X$.

\begin{lemma} \label[lemma]{Rim}

Let $X$ be a symplectic singularity and let $A' \to A$ and $B \to A$ be maps in $\Art_{\C}$, where the first map is surjective. Then the map 
$$
\text{res}: \PD_X(A' \times_{A} A'') \to \PD_X(A') \times_{\PD_X(A)} \PD_X(A'')
$$
is an equivalence.
\end{lemma}

\begin{proof}

Since any $A \in \Art_{\C}$ is a thickening of $k$, any $X_A \in \PD_X(A)$ is a thickening of $X$ and thus their underlying topological spaces $|X|$ and $|X_A|$ are the same. In particular, we can think of $X_A$ as a sheaf of Poisson $A$-algebras on $X$. Thus given an object $(X_{B}, X_{A'}, \psi) \in \PD_X(B) \times_{\PD_X(A)} \PD_X(A')$, we can explicitly construct an Poisson scheme $X_{B'} \in \PD_X(B')$ by setting 
\begin{equation*}
X_{B'}(U) := X_{B}(U) \times_{\psi, X_{A}(U)} X_{A'}(U).
\end{equation*}
Since all the data commutes with localization, this indeed forms a sheaf of Poisson $B'$-algebras on $X$ (in fact, the pushout $X_{B} \sqcup_{\psi,  X_{A}} X_{A'}$), and thus an element of $\PD_X(B')$. Furthermore, by \Cref{FlatRim} the resulting sheaf is flat. This shows that the functor is essentially surjective. 

To show that the functor is fully faithful,  let $X_1, X_2 \in \PD_X(B')$, and let $U$ be an open set in the underlying topological space $|X| = |X_1| = |X_2|$. If we consider $X_1(U)$ only as a flat $B'$-module, \Cref{FlatRim} states that it is a fiber product of flat modules over $B$ and $A'$. Furthermore, its algebra and Poisson structures are both defined by homomorphisms in $\Flat_{B'}$. By \Cref{FlatRim} again, these morphisms between fiber products are determined by their restrictions to each component, which tells us that $X_1(U)$ and $X_2(U)$ are fiber products as a Poisson $B'$-algebras. Thus $\Hom(X_1(U), X_2(U))$ is isomorphic to 
\begin{equation*}
\Hom(X_1|_{B}(U), X_2|_{B}(U)) \times_{\Hom(X_1|_{A}(U), X_2|_{A}(U))} \Hom(X_1|_{A'}(U), X_2|_{A'}(U))
\end{equation*}
This bijection commutes with localization and thus holds at a level of Poisson schemes, showing the map is fully faithful. 
\end{proof}

We recall some constructions and theorems from \cite{NaPD1}. In particular, computations regarding the Poisson deformation functor $\pd_X$ are usually made using Poisson cohomology.  If $X$ is smooth, and $\cX \in \pd_X(S)$ is a deformation of $X$ over $T := \Spec(S)$, $\HP^*(\cX/T)$ is given by the cohomology of the Poisson-Lichnerowicz complex $(\wedge^{\geq 1}\Theta_{\cX/T}, \delta)$, where $\delta$ is defined in \cite[p. 271-272]{NaFlops} and $\Theta_{\cX/T}$ is the sheaf of relative differentials.  The following key theorem allows us to compute Poisson cohomology using the de Rham complex.

\begin{prop}[{\cite[Proposition 9]{NaFlops}}] \label[prop]{smoothHP}

For a smooth symplectic variety ($\cX, \omega$) over $T$, the symplectic form induces a quasi-isomorphism between $(\wedge^{\geq 1}\Theta_{\cX/T}, \delta)$ and the truncated de Rham complex $(\Omega_{\cX/T}^{\geq 1}, d)$. 

\end{prop}

We define $PG_{\cX/T}$ to be the sheaf which associates to any open set $U \subset X$ the set of Poisson automorphisms of $\cX|_U$ which restrict to the identity over the central fiber of $T$. This sheaf is essential to showing that $\pd_X$ is prorepresentable, and as indicated above, it can be computed using Poisson cohomology. In particular, let $P\Theta_{\cX/T}$ be the first cohomology sheaf of $(\wedge^{\geq 1}\Theta_{\cX/T}, \delta)$, and let $P\Theta^0_{\cX/T}$ be the subsheaf of $P\Theta_{\cX/T}$ that restricts to zero over the central fiber. 

\begin{prop}[{\cite[Proposition 2.2]{NaPD1}}] \label[prop]{PGPO}

For a smooth $T$-scheme $\cX$, the sheaf $PG_{\cX/T}$ is isomorphic to $P\Theta^0_{\cX/T}$.

\end{prop}

For the remainder of this section we fix a partial resolution $\rho: X' \rightarrow X$ of a conical affine symplectic singularity $X$, and let $\pi: Y \rightarrow X$ be a $\Q$-factorial terminalization which covers $\rho$ via the map $\pi': Y \rightarrow X'$. We are interested in showing that the Poisson deformation functor $\pd_{X'}$ is prorepresentable. As a tool, we will also be considering other Zariski open subsets $V \subset X'$ satisfying $\rho^{-1}(X_{\reg}) \subset V$, but in all the cases we will be interested in we will have $\pd_V \cong \pd_{X'}$.

The following lemma will be essential for proving that $X'$ and various $V$ satisfy the final condition in Schlessinger's criterion.

\begin{lemma} \label[lemma]{pi1}

Let $V \subset X'$ be a Zariski open subset such that $\rho^{-1}(X_{\reg}) \subset V$. Then $H^1(V_{\reg}, \C) = 0$.

\end{lemma}

\begin{proof}

Since $X$ is a symplectic singularity, it has canonical singularities and thus $(X, 0)$ has klt singularities. By \cite[Theorem 1]{Bra}, this implies that $\pi_1(X_{\reg})$ is finite. Since $\rho$ is a partial resolution, $X_{\reg}$ is a Zariski open subset of $X'_{\reg}$ and thus of $V_{\reg}$. By \cite[Proposition 2.10]{Ko}, this means that the map $i_*: \pi_1(X_{\reg}) \rightarrow \pi_1(V_{\reg})$ induced by the inclusion is surjective. Thus $\pi_1(V_{\reg})$ is finite, so $H_1(V_{\reg}, \Z) \cong \Ab(\pi_1(V_{\reg}))$ is finite as well. Thus $H^1(V_{\reg}, \C) = 0$.
\end{proof}

In the following lemma we apply the reasoning of \cite[Corollary 2.5]{NaPD1} to partial resolutions:

\begin{lemma} \label[lemma]{AutSurj}

Let $V \subset X'$ be a Zariski open subset with $\rho^{-1}(X_{\reg}) \subset V$ and $H^1(V, \cO_V) = 0$. Then for every surjection $S \rightarrow \overline{S}$ in $\Art_{\C}$ and every $\cV \in \PD_{V}(S)$, the induced map 
$$
H^0(V, PG_{\cV /T}) \rightarrow H^0(V , PG_{\overline{\cV}/\overline{T}})
$$
is surjective, where $T := \Spec{S}$, $\overline{T} := \Spec({\overline{S}})$, and $\overline{\cV}$ is the restriction of $\cX$ to $\overline{T}$.

\end{lemma}

\begin{proof}

Let $j: \cU \to \cV$ be the inclusion of the smooth locus $\cU$ of $\cV$ over $T$, and similarly let $\ol{j}: \ol{\cU} \to \ol{\cV}$ be the inclusion of the smooth locus $\ol{\cU}$ of $\ol{\cV}$ over $\ol{T}$. As stated in the proof of \cite[Corollary 2.5]{NaPD1}, we have 
\begin{equation} \label{resSmooth}
j_*\cO_{\cU} = \cO_{\cV} \hspace{5pt} \text{and} \hspace{5pt} \ol{j}_*\cO_{\ol{\cU}} = \cO_{\ol{\cV}}.
\end{equation}
In particular, every Poisson automorphism of $\cV$ is determined by its restriction to $\cU$ and every Poisson automorphism of $\overline{\cV}$ is determined by its restriction to $\overline{\cU}$. Thus, it is enough to show that the restriction map $H^0(\cU, PG_{\cU/T}) \rightarrow H^0(\ol{\cU}, PG_{\overline{\cU}/\overline{T}})$ is surjective. By \Cref{PGPO}, this is isomorphic to the restriction map $H^0(\cU, P\Theta^0_{\cU/T}) \rightarrow H^0(\ol{\cU}, P\Theta^0_{\overline{\cU}/\overline{T}})$. In particular, this means it is sufficient to show the map $H^0(U, P\Theta_{\cU/T}) \rightarrow H^0(U, P\Theta_{\overline{\cU}/\overline{T}})$ is surjective. 

Since $\cU$ and $\overline{\cU}$ are smooth over $T$ and $\overline{T}$ respectively, we can apply \Cref{smoothHP} and compute Poisson cohomology as the cohomology of the truncated de Rham complex. In particular, we have $P\Theta_{\cU/T} \cong H^1(U, \Omega_{\cU/T}^{\geq 1})$, and similarly for $\overline{\cU}$. To compute this we use the following distinguished triangle 
\begin{equation} \label{eq:ses}
0 \rightarrow \Omega_{\cU/T}^{\geq 1} \rightarrow \Omega_{\cU/T} \rightarrow \cO_{\cU} \rightarrow 0.
\end{equation}

Note that by \cite[p. 276]{NaFlops}, $H^i(U, \Omega_{\cU/T}) \cong H^i(U, S)$. This distinguished triangle is functorial in $\Art_{\C}$, so taking the long exact sequence in cohomology for both $\cU$ and $\overline{\cU}$ we obtain the diagram 
\begin{equation} \label{Sch4}
\begin{tikzcd}
	{H^0(U, S)} && {H^0(\cU, \cO_{\cU})} && {\HP^1(\cU/T)} && {H^1(U, S)} \\
	\\
	{H^0(U, S)} && {H^0(\overline{\cU}, \cO_{\overline{\cU}})} && {\HP^1(\overline{\cU}/\overline{T})} && {H^1(U, \overline{S})}
	\arrow[from=1-1, to=1-3]
	\arrow[from=1-1, to=3-1]
	\arrow[from=3-1, to=3-3]
	\arrow[from=1-3, to=3-3]
	\arrow[from=1-3, to=1-5]
	\arrow[from=3-3, to=3-5]
	\arrow[from=1-5, to=3-5]
	\arrow[from=1-5, to=1-7]
	\arrow[from=3-5, to=3-7]
	\arrow[from=1-7, to=3-7]
\end{tikzcd}
\end{equation}
By \Cref{pi1}, $H^1(U, S) = H^1(U, \C) \otimes S = 0$ and similarly $H^1(U,  \overline{S}) = 0$. Finally, $H^1(V, \cO_V) = 0$ implies that $H^0(\cV, \cO_{\cV}) \rightarrow H^0(\ol{\cV}, \cO_{\ol{\cV}})$ is surjective. By \eqref{resSmooth}, this map is isomorphic to $H^0(\cU, \cO_{\cU}) \rightarrow H^0(\overline{\cU}, \cO_{\overline{\cU}})$, so the latter is also surjective. Thus, \eqref{Sch4} implies that $\HP^1(\cU/T) \rightarrow \HP^1(\overline{\cU}/\overline{T})$ is surjective as well. Using a spectral sequence argument, we have $\HP^1(\cU/T) \cong H^0(X, P\Theta_{\cU/T})$ and $\HP^1(\overline{\cU}/\overline{T}) \cong H^0(X, P\Theta_{\overline{\cU}/\overline{T}})$, so we are done.
\end{proof}

Finally, will now switch our focus to the third condition in Schlessinger's criterion. Given $A \in \Art_{\C}$, we have a forgetful map $\pd_X(A) \rightarrow \Def_X(A)$ from Poisson deformations to flat deformations. Note that Poisson and flat deformations both form presheaves in the complex analytic topology. Following Namikawa, we let $PT^1_X$ and $T^1_X$ be the sheafification of the presheaves given by first order Poisson and flat deformations respectively. Thus, for any open set $U \subset X$, there are canonical maps $\pd_U(\C[\epsilon]) \rightarrow PT^1_X(U)$ and $\Def_U(\C[\epsilon]) \rightarrow T^1_X(U)$. Finally, we also have a forgetful map of sheaves $PT^1_X \rightarrow T^1_X$, and for any $U \subset X$ all of the above maps fit into the commutative diagram:

\[\begin{tikzcd}
	{\pd_U(\C[\epsilon])} && {\Def_U(\C[\epsilon])} \\
	\\
	{PT^1_X(U)} && {T^1_X(U)}
	\arrow[from=1-1, to=1-3]
	\arrow[from=1-1, to=3-1]
	\arrow[from=3-1, to=3-3]
	\arrow[from=1-3, to=3-3]
\end{tikzcd}\]

The kernel of the map composition $\pd_U(\C[\epsilon]) \rightarrow T^1_X(U)$ is the set of first order Poisson deformations that are locally trivial as flat deformations. We denote this $\pd_{lt, U}(\C[\epsilon])$. Additionally, let $\cH$ be the subsheaf of $T^1_X$ given by the image of $PT^1_X$. Thus we have the exact sequence 
\begin{equation} \label{ses:lt}
0 \rightarrow \pd_{lt, U}(\C[\epsilon]) \rightarrow \pd_U(\C[\epsilon]) \rightarrow \Gamma(U, \cH) 
\end{equation}

\begin{prop} \label[prop]{ltpd}

Let $X$ be a symplectic singularity such that every symplectic leaf has codimension 0 or 2 and $H^i(X, \cO_X) = 0$ for $i = 1, 2$.  Then $\pd_{lt, X}(\C[\epsilon]) = H^2(X, \C)$. 

\end{prop}

\begin{proof}

The proof is identical \cite[Proposition 1.11]{NaPD1}.
\end{proof}

This theorem has the following key corollary:

\begin{cor} \label[cor]{Partlt}

Let $X'$ be a partial resolution of an affine symplectic singularity and let $V \subset X'$ be the union of the open stratum and the codimension two strata. Then $\pd_{lt, V}(\C[\epsilon]) = H^2(V, \C)$. 

\end{cor}

\begin{proof} 

We have that $H^i(X', \cO_{X'}) = 0$ for all $i > 0$. Additionally since $X'$ is Cohen-Macaulay and $X' - V$ has codimension four, by a depth argument we have $H^i(V, \cO_{V}) = 0$ for $i = 1, 2$. Thus, by \Cref{ltpd}, we are done. 
\end{proof}

Finally we are in a position to prove the main result of this subsection.

\begin{prop} \label[prop]{ProRep}

Given a partial resolution $\rho: X' \rightarrow X$ of a conical affine symplectic singularity $X$, $\pd_{X'}$ is prorepresentable.

\end{prop}

\begin{proof}

We shall show that $\pd_{X'}$ satisfies the four conditions of Schlessinger's criterion, \Cref{SchThm}. By \cite[\href{https://stacks.math.columbia.edu/tag/06J5}{Lemma 06J5}]{Stacks}, \Cref{Rim} shows that the $2$-functor $\PD_{X'}$ satisfies the Rim-Schlessinger condition. This implies that $\pd_{X'}$ satisfies the first two conditions of Schlessinger's criterion, as we shall explain below. 

For $2$-functors, there are conditions $(S_1)$ and $(S_2)$ that naturally generalize the first two conditions in Schlessinger's criterion. These can be found in \cite[\href{https://stacks.math.columbia.edu/tag/06HW}{Definition 06HW}]{Stacks}, and it is explained in \cite[\href{https://stacks.math.columbia.edu/tag/06HY}{Remark 06HY}]{Stacks} that these agree with the first two conditions in \Cref{SchThm} for functors valued in sets. Furthermore, it is shown in \cite[\href{https://stacks.math.columbia.edu/tag/06J7}{Lemma 06J7}]{Stacks} that the Rim-Schlessinger condition implies $(S_1)$ and $(S_2)$, so $\PD_{X'}$ satisfies $(S_1)$ and $(S_2)$. Finally, it is shown in \cite[\href{https://stacks.math.columbia.edu/tag/06I0}{Lemma 06I0}]{Stacks} that if a $2$-functor $\cF$ satisfies $(S_1)$ and $(S_2)$ and outputs the trivial groupoid when the input is the ground field, then the underlying functor $\ol{\cF}$ to sets also satisfies $(S_1)$ and $(S_2)$. Thus, since $\PD_{X'}(\C)$ is a point, and $\pd_{X'}$ is its underlying functor to sets of $\PD_{X'}$, $\pd_{X'}$ satisfies the first two conditions of Schlessinger's criterion. 

Additionally, by \cite[\href{https://stacks.math.columbia.edu/tag/06J8}{Lemma 06J8}]{Stacks}, the final condition is equivalent to the condition that for every (small) surjection $S \rightarrow \overline{S}$ in $\Art_{\C}$, and every $\cX \in \PD_X(S)$, the induced map $\Aut(\cX) \rightarrow \Aut(\overline{\cX})$ is surjective, where $\overline{\cX}$ is the restriction of $\cX$ to $\overline{S}$. This is precisely the content of \Cref{AutSurj}.

The only statement left to prove is that $\pd_{X'}(\C[\epsilon])$ is finite dimensional. By \cite[p. 62]{NaPD1} we have that $\pd_{X'}(\C[\epsilon]) \cong \pd_{V}(\C[\epsilon])$ where $V \subset X'$ is the union of the smooth and codimension two leaves of $X'$. We apply diagram \eqref{ses:lt} to $V$. By \cite[Lemma 1.5]{NaPD1}, $\cH$ is a finite dimensional local system on the singular locus of $V$, so $\Gamma(V, \cH)$ is finite dimensional. By \Cref{Partlt}, $\pd_{lt, V}(\C[\epsilon])$ is finite dimensional as well. So by diagram \eqref{ses:lt}, $\pd_{X'}(\C[\epsilon])$ is finite dimensional.
\end{proof}

Next we relate $\pd_{X'}$ to the Poisson deformation functors of $X$ and $Y$. 

\begin{lemma}

We have natural pushforward maps $\pd_Y \xrightarrow{\pi'_*} \pd_{X'} \xrightarrow{\rho_*} \pd_X$. 

\end{lemma}

\begin{proof}

We construct the map $\pi'_*$, as the construction of $\rho_*$ is identical. By \cite[Theorem 1.4, Remark 1.4.1]{Wa}, if $\pi'_*\cO_Y \cong \cO_{X'}$ and $R^1\pi'_*\cO_Y = 0$, then there is a pushforward map $\ol{\pi}'_*: \Def_Y \rightarrow \Def_{X'}$ on flat deformations. For another construction of this map see the argument in \cite[\href{https://stacks.math.columbia.edu/tag/0E3X}{Lemma 0E3X}]{Stacks}, which can be used to construct the map $\ol{\pi}'_*: \Def_Y \rightarrow \Def_{X'}$ by taking $\cY \in \Def_Y(A)$ as a sheaf of flat $A$-algebras on $Y$ and pushing it forward to a sheaf of $A$-algebras on $X'$.

We take a resolution of singularities $\nu: Z \rightarrow Y$, so $\pi' \circ \nu$ is a resolution of singularities of $X'$. Furthermore, since $Y$ and $X'$ have rational singularities, $\nu_*\cO_Z$ is quasi-isomorphic to $\cO_Y$ in the derived category and $(\pi' \circ \nu)_*(\cO_Z) \cong (\pi'_* \circ \nu_*)(\cO_Z)$ is also quasi-isomorphic to $\cO_{X'}$ in the derived category. Thus, $\pi'_*\cO_Y$ is quasi-isomorphic to $\cO_{X'}$, which indeed implies that $\pi'_*\cO_Y \cong \cO_{X'}$ and $R^1\pi'_*\cO_Y = 0$. Thus, we have a pushforward map of flat deformations $\ol{\pi}'_*: \Def_Y \rightarrow \Def_{X'}$.

To extend this to a pushforward map $\pi'_*: \pd_Y \rightarrow \pd_{X'}$, we observe that since $\pi'$ is an isomorphism over $X'_{\reg}$, a Poisson deformation of $Y$ restricts to a Poisson deformation of $X'_{\reg}$. Let $\cY \in \pd_Y(A)$, for some $A \in \Art_{\C}$. Since the above pushforward map of flat deformations is the identity over $X'_{\reg}$, we can simply extend this Poisson structure on $\cY|_{X'_{\reg}} \cong \ol{\pi}'_*(\cY)|_{X'_{\reg}}$ to $\ol{\pi}'_*(\cY)$ by normality. 
\end{proof}

\begin{rem}

The above pushforward maps compose to standard pushforward map $\pi_*: \pd_Y \rightarrow \pd_X$ studied by Namikawa. This map (as well as the map $\rho_*$) can be described more succinctly by sending some $\cY \in \pd_Y(A)$ to $\Spec\Gamma(Y, \cO_{\cY})$, which is a deformation of $X = \Spec\Gamma(Y, \cO_Y)$.

\end{rem}

\subsection{Unobstructedness} \label[subsection]{AllSmoothSec}

The primary goal of this section is to prove that $\pd_{X'}$ is unobstructed. In the course of the proof we will need to consider $\tilde{V} := \pi'^{-1}(V)$ where $V$ is the union of the open and codimension two leaves of $X'$. We will also be considering the union $U$ of the open and codimension two leaves of $X$, as well as $U' := \rho^{-1}(U)$ and $\tilde{U} := \pi^{-1}(U)$. Considering all the spaces mentioned above, we obtain the following commutative diagram:
\begin{equation} \label{AllSubs}
\begin{tikzcd}
	{\pd_{Y}} && {\pd_{Y_{\reg}}} && {\pd_{\tilde{V}}} && {\pd_{\tilde{U}}} \\
	\\
	{\pd_{X'}} &&&& {\pd_{V}} && {\pd_{U'}} \\
	\\
	{\pd_{X}} &&&&&& {\pd_U}
	\arrow["{r_{\tilde{V}}}", from=1-3, to=1-5]
	\arrow["{(\pi'|_{\tilde{V}})_*}"', from=1-5, to=3-5]
	\arrow["{r_{\tilde{U}}}", from=1-5, to=1-7]
	\arrow["{(\pi'|_{\tilde{U}})_*}"', from=1-7, to=3-7]
	\arrow["{r_{U'}}", from=3-5, to=3-7]
	\arrow["{(\rho|_{U'})_*}"', from=3-7, to=5-7]
	\arrow["{r_{Y_{\reg}}}", from=1-1, to=1-3]
	\arrow["{\pi'_*}"', from=1-1, to=3-1]
	\arrow["{\rho_*}"', from=3-1, to=5-1]
	\arrow["{r_V}", from=3-1, to=3-5]
	\arrow["{r_U}", from=5-1, to=5-7]
\end{tikzcd}
\end{equation}

Throughout this section we prove that all the functors in the diagram below are unobstructed and that all the horizontal maps are equivalences. We note that many of these functors were already proven to be unobstructed in \cite[Corollary 13]{NaFlops} and \cite[Theorem 5.1]{NaPD1}, and many of the arguments in this subsection for the other spaces are adapted from these.

Before beginning the argument, we briefly outline the main techniques that we will use. We show that the restriction maps $r_{Y_{\reg}}$, $r_V$, and $r_U$ are equivalences because the complements $X - U$, $X' - V$, and $Y - Y^{\reg}$ all have codimension at least 4. Additionally, since $Y_{\reg}$, $\tilde{V}$, and $\tilde{U}$ are smooth, we can describe their first order Poisson deformations in terms of their second de Rham cohomologies, which are all isomorphic. Furthermore, we can show they are unobstructed using the $T^1$-lifting criterion, which is described below in the proof of \Cref{TopRow}. Next we note that $V$, $U'$, and $U$ have only Kleinian singularities, are these are resolved by $\tilde{V}$ and $\tilde{U}$. Thus, the techniques in \cite[Section 4]{NaPD1} and \cite[Theorem 5.1]{NaPD1} can be used to deduce the unobstructedness of $\pd_{V}$, $\pd_{U'}$, and $\pd_{U}$ from the unobstructedness of $\pd_{\tilde{V}}$ and $\pd_{\tilde{U}}$. Using the pushforward map $\pd_{\tilde{V}} \to \pd_{V}$ as an example, the key idea is that the induced map $\Spec(R_{\tilde{V}}) \rightarrow \Spec(R_V)$ on prorepresenting rings has closed finite fiber. This fact along with the regularity of $R_{\tilde{V}}$ and a dimension count on both tangent spaces implies the regularity of $R_V$. Finally, we show that $r_{U'}: \pd_{V} \to \pd_{U'}$ is an isomorphism by showing it induces an isomorphism on first order deformations, which we do by considering the sequence \eqref{ses:lt} functorially for both $V$ and $U'$.

\begin{lemma}

All of the functors in diagram \eqref{AllSubs} are prorepresentable.

\end{lemma}

\begin{proof}

This follows from Schlessinger's criterion, \Cref{SchThm}. \Cref{Rim} implies that the first two conditions are satisfied for the deformation functor of any symplectic singularity, as explained in the proof of \Cref{ProRep}. The third criterion is satisfied by applying the short exact sequence \eqref{ses:lt} to each space. Finally, the last criterion follows from \Cref{AutSurj}, since all the spaces considered satisfy the hypotheses of that lemma.
\end{proof}

\begin{lemma} \label[lemma]{TopRow}

The functors $\pd_{Y_{\reg}}$, $\pd_{\tilde{V}}$, and $\pd_{\tilde{U}}$ in the top row of diagram \eqref{AllSubs} are all unobstructed, and the restriction maps between them are equivalences.

\end{lemma}

\begin{proof}

Because $X$, $X'$, and $Y$ are all Cohen-Macaulay and $X - U$, $X' - V$, and $Y - Y^{\reg}$ all have codimension at least 4, by a depth argument we have 
$$
H^i(U, \cO_U) = H^i(V, \cO_V) = H^i(Y^{\reg}, \cO_{Y^{\reg}}) = 0 \hspace{10pt} \text{for} \hspace{10pt} i = 1, 2.
$$ 
Furthermore, since $U$ and $V$ have rational singularities, we have 
$$
H^i(\tilde{U}, \cO_{\tilde{U}}) = H^i(\tilde{V}, \cO_{\tilde{V}}) = 0 \hspace{10pt} \text{for} \hspace{10pt} i = 1, 2.
$$
Since $\tilde{U}$ and $\tilde{V}$ are Zariski open subsets of $Y^{\reg}$, we can use the functoriality of the short exact sequence of sheaves
\begin{equation*}
0 \rightarrow \Omega_{Y^{\reg}}^{\geq 1} \rightarrow \Omega_{Y^{\reg}} \rightarrow \cO_{Y^{\reg}} \rightarrow 0
\end{equation*}
under restriction. Additionally, we recall that since $Y^{\reg}$ is smooth, by \cite[Proposition 8]{NaFlops}, $\pd_{Y^{\reg}}(\C[\epsilon]) \cong H^2(Y^{\reg}, \Omega_{Y^{\reg}}^{\geq 1})$ and the analogous statements are true for $\tilde{U}$ and $\tilde{V}$. Thus taking long exact sequences in cohomology, we have the diagram

\[\begin{tikzcd} [row sep=scriptsize]
	{H^1(Y^{\reg}, \cO_{Y^{\reg}})} & {\pd_{Y^{\reg}}(\C[\epsilon])} & {H^2(Y^{\reg}, \C)} & {H^2(Y^{\reg}, \cO_{Y^{\reg}})} \\
	\\
	{H^1(\tilde{V}, \cO_{\tilde{V}})} & {\pd_{\tilde{V}}(\C[\epsilon])} & {H^2(\tilde{V}, \C)} & {H^2(\tilde{V}, \cO_{\tilde{V}})} \\
	\\
	{H^1(\tilde{U}, \cO_{\tilde{U}})} & {\pd_{\tilde{U}}(\C[\epsilon])} & {H^2(\tilde{U}, \C)} & {H^2(\tilde{U}, \cO_{\tilde{U}})}
	\arrow[from=1-2, to=3-2]
	\arrow[from=1-1, to=1-2]
	\arrow[from=1-1, to=3-1]
	\arrow[from=3-1, to=3-2]
	\arrow[from=1-2, to=1-3]
	\arrow[from=3-2, to=3-3]
	\arrow[from=1-3, to=3-3]
	\arrow[from=1-3, to=1-4]
	\arrow[from=1-4, to=3-4]
	\arrow[from=3-3, to=3-4]
	\arrow[from=3-1, to=5-1]
	\arrow[from=3-2, to=5-2]
	\arrow[from=5-1, to=5-2]
	\arrow[from=3-3, to=5-3]
	\arrow[from=5-2, to=5-3]
	\arrow[from=3-4, to=5-4]
	\arrow[from=5-3, to=5-4]
\end{tikzcd}\]

Since all the cohomology groups of the structure sheaves are zero, this implies that $r_{\tilde{V}}(\C[\epsilon])$ and $r_{\tilde{U}}(\C[\epsilon])$ in diagram \eqref{AllSubs} are given by restriction in cohomology. Furthermore, it follows from the proof of \cite[Proposition 2]{Na1} that these maps are both isomorphisms. 

Next we show these functors are unobstructed. We begin by briefly recalling the $T^1$-lifting criterion, which can be proved using the argument in \cite[Corollary 13]{NaFlops}. Let $S_n := k[x]/(x^n)$ and $T_n := \Spec(S_n)$. For any deformation functor $D$ and $X_n \in D(S_n)$, we let $T^1(X_n/T_n)$ be the set of elements of $D(S_n[\epsilon])$ that restrict to $X_n$ over $T_n$. The $T^1$-lifting criterion says that if $D$ is a prorepresentable functor, then it is unobstructed if and only if the natural map 
\begin{equation*}
T^1(X_n/T_n) \rightarrow T^1(X_{n - 1}/T_{n - 1})
\end{equation*}
induced by $S_n[\epsilon] \rightarrow S_{n - 1}[\epsilon]$ is a surjection for all positive $n$ and any (equivalently, some) $X_n \in D(S_n)$. If $D$ is the Poisson deformation functor $\pd_Z$ of some symplectic singularity $Z$, we denote $T^1(X_n/T_n)$ by $\pd_Z(X_n/T_n, T_n[\epsilon])$, following \cite[p. 272]{NaFlops}.

We have shown that all of the functors in diagram \eqref{AllSubs} are prorepresentable, so we will use the $T^1$-lifting criterion. We now consider the deformation functors of the smooth symplectic manifolds $Y_{\reg}$, $\tilde{V}$, and $\tilde{U}$. We prove the $T^1$-lifting criterion for $\tilde{V}$ since the other two proofs are identical. Applying \cite[Proposition 8, Corollary 10]{NaFlops} and the fact that $H^i(\tilde{V}, \cO_{\tilde{V}}) = 0$ for $i = 1, 2$, we obtain for any $\tilde{V}_n \in \pd_{\tilde{V}}(S_n)$
\begin{equation*}
\pd_{\tilde{V}}(\tilde{V}_n/T_n, T_n[\epsilon]) \cong HP^2(\tilde{V}_n/T_n) \cong H^2(\tilde{V}, S_n) \cong H^2(\tilde{V}, \C) \otimes S_n.
\end{equation*}
Additionally, the above isomorphisms are functorial in $\Art_{\C}$. Thus, since $S_n \rightarrow S_{n - 1}$ is surjective, $\pd_{\tilde{V}}$ is unobstructed, and the same argument applies to the other two functors. Finally, the maps $r_{\tilde{V}}$ and $r_{\tilde{U}}$ of diagram \eqref{AllSubs} are maps between prorepresentable and unobstructed deformations functors that induce isomorphisms on first order deformations, so they are equivalences. 
\end{proof}

\begin{prop} \label[prop]{AllSmooth}

All the functors in diagram \eqref{AllSubs} are unobstructed, and they all have tangent spaces of the same dimension.

\end{prop}

\begin{proof}

We begin by observing that \cite[Section 4]{NaPD1} applies to any morphism $Y \rightarrow X$ where $X$ is a symplectic singularity whose singularities are all locally isomorphic to $(S, 0) \times (\C^{2n -2}, 0)$, for $S$ Kleinian, and $Y$ is a symplectic resolution. In particular, it applies to $\pi'|_{\tilde{V}}: \tilde{V} \rightarrow V$, $\pi'|_{\tilde{U}}: \tilde{U} \rightarrow U'$, and $\pi|_{\tilde{U}}: \tilde{U} \rightarrow U$.

The sheaf $R^2\pi'_*\C_{\tilde{V}}$ is given by the sheafification of the presheaf
\begin{equation*}
W \mapsto \pd_{\pi'^{-1}(W)}(\C[\epsilon]).
\end{equation*}
Thus, we have an exact sequence
\begin{equation} \label{tildV:ses}
0 \rightarrow \pd_{\tilde{V}}^{\pi'}(\C[\epsilon]) \rightarrow \pd_{\tilde{V}}(\C[\epsilon]) \rightarrow \Gamma(V, R^2\pi'_*\C_{\tilde{V}}),
\end{equation}
where $\pd_{\tilde{V}}^{\pi'}$ is the subfunctor of $\pd_{\tilde{V}}$ given by deformations which are trivial on the preimage of a sufficiently small analytic neighborhood of every point in $V$. These types of subfunctors will be studied more extensively in the next two sections. By the argument in \cite[Lemma 3]{NaMori} (using the Leray spectral sequence), the map $\pd_{\tilde{V}}^{\pi'}(\C[\epsilon]) \rightarrow \pd_{\tilde{V}}(\C[\epsilon])$ is given by the pullback map in cohomology $H^2(V, \C) \rightarrow H^2(\tilde{V}, \C)$. Additionally, as explained in the proof of \Cref{Hgen}, $\Gamma(V, R^2\pi'_*\C_{\tilde{V}})$ has a basis by the irreducible components of the exceptional divisor of $\pi'$ (see the first sentence of \cite[p. 78 (ii)]{NaPD1}). Thus, $H^2(\tilde{V}, \C) \rightarrow \Gamma(V, R^2\pi'_*\C_{\tilde{V}})$ is surjective and \eqref{tildV:ses} is a short exact sequence. 

Next we apply the short exact sequence \eqref{ses:lt} to $V$, obtaining
\begin{equation} \label{Vlt}
0 \rightarrow \pd_{lt, V}(\C[\epsilon]) \rightarrow \pd_V(\C[\epsilon]) \rightarrow \Gamma(V, \cH).
\end{equation}
For the next proof we want to show that the sequence above is also short exact, in particular that the final map is surjective. We can apply \cite[Proposition 4.2]{NaPD1} to $\tilde{V} \rightarrow V$ to obtain $\dim \Gamma(V, \cH) = \dim \Gamma(V, R^2\pi'_*\C_{\tilde{V}})$ and \Cref{ltpd} to obtain $\dim \pd_{lt, V}(\C[\epsilon]) = \dim H^2(V, \C)$. Thus, we have 
\begin{align*}
\dim \pd_V(\C[\epsilon]) &\leq \dim \pd_{lt, V}(\C[\epsilon]) + \dim \Gamma(V, \cH) \\
&=  \dim H^2(V, \C) + \dim \Gamma(V, R^2\pi'_*\C_{\tilde{V}}) \\
&= \dim H^2(\tilde{V}, \C) = \dim \pd_{\tilde{V}}(\C[\epsilon]).
\end{align*}

Since all of the singularities of $V$ are of the form $(S, 0) \times (\C^{2n -2}, 0)$ and $\pi'|_{\tilde{V}}: \tilde{V} \rightarrow V$ is a symplectic resolution, the argument in part $(iii)$ of \cite[Theorem 5.1]{NaPD1} implies that the pushforward map
\begin{equation*}
(\pi'|_{\tilde{V}})_*: \Spec(R_{\tilde{V}}) \rightarrow \Spec(R_V)
\end{equation*}
has finite closed fiber, where $R_{\tilde{V}}$ and $R_V$ are the complete local rings that prorepresent $\pd_{\tilde{V}}$ and $\pd_V$ respectively. Letting $\mathfrak{m}_{\tilde{V}}$ and $\mathfrak{m}_{V}$ be the maximal ideals of $R_{\tilde{V}}$ and $R_V$ respectively, this means that $\dim R_{\tilde{V}}/\mathfrak{m}_{V}R_{\tilde{V}} = 0$. Thus, we have 
\begin{equation*}
\dim R_{\tilde{V}} \leq \dim R_{V} + \dim R_{\tilde{V}}/\mathfrak{m}_{V}R_{\tilde{V}} = \dim R_{V}.
\end{equation*}
We also know that 
$$
\dim R_{V} \leq \dim_{\C} \mathfrak{m}_{V}/(\mathfrak{m}_{V})^2 \hspace{10pt} \text{and} \hspace{10pt} \dim R_{\tilde{V}} = \dim_{\C} \mathfrak{m}_{\tilde{V}}/(\mathfrak{m}_{\tilde{V}})^2,
$$ 
since $\pd_{\tilde{V}}$ is unobstructed and thus $R_{\tilde{V}}$ is regular. Combining the above equations we obtain
\begin{equation*}
\dim_{\C} \mathfrak{m}_{\tilde{V}}/(\mathfrak{m}_{\tilde{V}})^2 \leq \dim R_{V} \leq \dim_{\C} \mathfrak{m}_{V}/(\mathfrak{m}_{V})^2.
\end{equation*}
However these two tangent spaces are given by $\pd_{\tilde{V}}(\C[\epsilon])$ and $\pd_{V}(\C[\epsilon])$ respectively. Since we proved above that $\dim \pd_V(\C[\epsilon]) \leq \dim \pd_{\tilde{V}}(\C[\epsilon])$, the inequalities above must be equalities. Since $\dim R_{V} = \dim_{\C}\mathfrak{m}_{V}/(\mathfrak{m}_{V})^2$, $R_{V}$ is a regular local ring and $\pd_{V} = \Hom(R_{V}, -)$ is unobstructed. This equality also implies that the sequence \eqref{Vlt} is a short exact sequence.

Note that the only features of the resolution $\pi'|_{\tilde{V}}: \tilde{V} \rightarrow V$ that we used were that all of the singularities of $V$ are of the form $(S, 0) \times (\C^{2n -2}, 0)$, that $\pi'|_{\tilde{V}}: \tilde{V} \rightarrow V$ is a symplectic resolution, and that $H^i(V, \cO_V) = 0$ for $i = 1, 2$ (to apply \Cref{ltpd}). These all apply to $\pi'|_{\tilde{U}}: \tilde{U} \rightarrow U'$, and $\pi|_{\tilde{U}}: \tilde{U} \rightarrow U$ as well. In particular, 
\begin{equation} \label{H12U'}
H^i(U, \cO_{U}) = H^i(U', \cO_{U'}) = H^i(\tilde{U}, \cO_{\tilde{U}}) = 0 \hspace{10pt} \text{for} \hspace{10pt} i = 1, 2,
\end{equation}
since $U$ and $U'$ have rational singularities. Thus, the same argument applies to show that $\pd_{U'}$ and $\pd_{U}$ are unobstructed and that the dimensions of the prorepresenting rings are the same as that of $\pd_{\tilde{U}}$.

We only have left to prove that the Poisson deformation functors of $Y$, $X'$, and $X$ are unobstructed. For this we use the fact that $Y - Y_{\reg}$,  $X' - V$, and $X - U$ all have codimension at least four. We will explain the details for $X'$ as the other two cases are identical. For any deformation $X'_n \in \pd_{X'}(S_n)$ we can restrict it to a deformation $V_n \in \pd_{V}(S_n)$. Applying the reasoning of \cite[Proposition 13]{NaFlops}, since $X' - V$ has codimension at least four, we have 
\begin{equation*}
\pd_{X'}(X'_n/T_n, T_n[\epsilon]) \cong \pd_{V}(V_n/T_n, T_n[\epsilon]).
\end{equation*}
Thus, $\pd_{X'}$ satisfies the $T^1$-lifting criterion if and only if $\pd_V$ does. Thus, $\pd_{X'}$ is unobstructed. Finally, since $\pd_{X'}$ and $\pd_V$ are both unobstructed and $r_V$ induces an isomorphism on first order deformations, it is an equivalence. 

As stated above, the same argument shows that $r_{Y_{\reg}}$ and $r_U$ are equivalences and in particular that $\pd_Y$ and $\pd_X$ are unobstructed.
\end{proof}

In the above arguments we have shown that all of the horizontal maps in diagram \eqref{AllSubs} are equivalences, except for $r_{U'}$. We show that below:

\begin{prop} \label[prop]{rU'equiv}

The map $r_{U'}$ is an equivalence.

\end{prop}

\begin{proof}

We have shown that $\pd_V$ and $\pd_{U'}$ are both prorepresentable and unobstructed, so we only need to show that $r_{U'}$ induces an isomorphism on first order deformations. We first study $r_{\tilde{U}}$. 

The construction of sequence \eqref{tildV:ses} is functorial under pullback from the right column to the left column in:
\begin{equation*}
\begin{tikzcd}
	{\tilde{U}} && {\tilde{V}} \\
	\\
	{U'} && V
	\arrow["{i_{\tilde{U}}}", hook, from=1-1, to=1-3]
	\arrow["{\pi'|_{\tilde{U}}}"', from=1-1, to=3-1]
	\arrow["{\pi'|_{\tilde{V}}}"', from=1-3, to=3-3]
	\arrow["{i_U}", hook, from=3-1, to=3-3]
\end{tikzcd}
\end{equation*}

Additionally, we showed in the previous argument that the sequences we obtain are short exact. Thus, we have the diagram:

\begin{equation} \label{Leray}
\begin{tikzcd}
	0 & {H^2(V, \C)} & {H^2(\tilde{V}, \C)} & {\Gamma(V, R^2\pi'_*\C_{\tilde{V}})} & 0 \\
	\\
	0 & {H^2(U', \C)} & {H^2(\tilde{U}, \C)} & {\Gamma(U', R^2\pi'_*\C_{\tilde{V}}) } & 0
	\arrow["{i_{\tilde{U}}^*}"', from=1-3, to=3-3]
	\arrow[from=1-3, to=1-4]
	\arrow[from=3-3, to=3-4]
	\arrow[from=1-4, to=3-4]
	\arrow[from=1-4, to=1-5]
	\arrow[from=3-4, to=3-5]
	\arrow["{(\pi'|_{\tilde{V}})^*}", from=1-2, to=1-3]
	\arrow[from=1-1, to=1-2]
	\arrow["{(\pi'|_{\tilde{U}})^*}", from=3-2, to=3-3]
	\arrow[from=3-1, to=3-2]
	\arrow["{i_{U'}^*}"', from=1-2, to=3-2]
\end{tikzcd}
\end{equation}

We know from the proof of \Cref{TopRow} that $i_{\tilde{U}}^*$ is an isomorphism. This implies that the right column of diagram \eqref{Leray} is surjective. Next we prove it is injective. It is shown in the proof of \cite[Proposition 2.14]{Lo} that the codimension of $Y - \tilde{U}$ in $Y$ is at least two. Thus, $\tilde{V} - \tilde{U}$ has codimension at least two in $\tilde{V}$ as well. This implies that the singular locus $\Sigma_{V - U'}$ of $V - U'$ must be codimension at least three, because otherwise its preimage would be a subset of $\tilde{V} - \tilde{U}$ of codimension less than two. Finally, this implies that $\Sigma_{V - U'}$ has codimension at least one in the singular locus $\Sigma_V$ of $V$. In particular, this means that the singular locus $\Sigma_{U'}$ of $U'$ intersects every connected component of $\Sigma_V$. Thus since $R^2\pi'_*\C_{\tilde{V}}$ is a local system on $\Sigma_V$, the restriction map in the right column of diagram \eqref{Leray} is injective. Thus, it is an isomorphism, and by the commutativity of diagram \eqref{Leray}, $i_{U'}^*$ is an isomorphism. 

We now study the map $\pd_V(\C[\epsilon]) \rightarrow \pd_{U'}(\C[\epsilon])$ induced by $r_{U'}$ directly. We apply the short exact sequence \eqref{ses:lt} to $U'$ and $V$, noting that it is functorial under restriction. Furthermore, in the proof of \Cref{AllSmooth} we saw that the sequences below are also short exact. Thus we have the diagram:

\begin{equation} \label{PartialsFunct}
\begin{tikzcd}
	0 & {\pd_{lt, V}(\C[\epsilon])} & {\pd_{V}(\C[\epsilon])} & {\Gamma(V, \cH)} & 0 \\
	\\
	0 & {\pd_{lt, U'}(\C[\epsilon])} & {\pd_{U'}(\C[\epsilon])} & {\Gamma(U', \cH) } & 0
	\arrow[from=1-3, to=3-3]
	\arrow[from=1-3, to=1-4]
	\arrow[from=3-3, to=3-4]
	\arrow[from=1-4, to=3-4]
	\arrow[from=1-2, to=1-3]
	\arrow[from=1-1, to=1-2]
	\arrow[from=3-2, to=3-3]
	\arrow[from=3-1, to=3-2]
	\arrow[from=1-2, to=3-2]
	\arrow[from=1-4, to=1-5]
	\arrow[from=3-4, to=3-5]
\end{tikzcd}
\end{equation}

It follows from the proof of \cite[Proposition 1.11]{NaPD1} (as applied to \Cref{ltpd}), that $\pd_{lt, V}(\C[\epsilon])$ is given by $H^2(V, \tilde{\Omega}_{V}^{\geq 1})$, where $(\tilde{\Omega}_{V}^{\geq 1}, d)$ is the truncated de Rham complex of $V$ as a $V$-manifold. Like the de Rham complex of a smooth variety, it sits in a distinguished triangle 
\begin{equation*}
\tilde{\Omega}_{V}^{\geq 1} \rightarrow \tilde{\Omega}_V \rightarrow \cO_V \rightarrow \tilde{\Omega}_{V}^{\geq 1}[1].
\end{equation*}

This is also true for $\pd_{lt, U'}(\C[\epsilon])$, and furthermore the above short exact sequence is functorial under restriction, so we obtain the diagram
\begin{equation} \label{ltFunct}
\begin{tikzcd}
	{H^1(V, \cO_V)} & {\pd_{lt, V}(\C[\epsilon])} & {H^2(V, \C)} & {H^2(V, \cO_V)} \\
	\\
	{H^1(U', \cO_{U'})} & {\pd_{lt, U'}(\C[\epsilon])} & {H^2(U', \C)} & {H^2(U', \cO_{U'})}
	\arrow[from=1-2, to=3-2]
	\arrow[from=1-1, to=1-2]
	\arrow[from=1-1, to=3-1]
	\arrow[from=3-1, to=3-2]
	\arrow[from=1-2, to=1-3]
	\arrow[from=3-2, to=3-3]
	\arrow[from=1-3, to=3-3]
	\arrow[from=1-3, to=1-4]
	\arrow[from=1-4, to=3-4]
	\arrow[from=3-3, to=3-4]
\end{tikzcd}
\end{equation}

However, we showed above that $H^i(V, \cO_V) = H^i(U', \cO_{U'}) = 0$ for $i = 1, 2$. Thus, by diagram \eqref{ltFunct}, the restriction map $\pd_{lt, V}(\C[\epsilon]) \rightarrow \pd_{lt, U'}(\C[\epsilon])$ is given by the restriction map in cohomology
\begin{equation*}
i_{U'}^*: H^2(V, \C) \rightarrow H^2(U', \C),
\end{equation*}
which we proved was an isomorphism. 

Additionally, applying \cite[Proposition 4.2]{NaPD1} to both $\tilde{V} \rightarrow V$ and $\tilde{U} \rightarrow U'$ and using the isomorphism in the right column of diagram \eqref{Leray} we have 
\begin{equation*}
\dim \Gamma(V, \cH) = \dim \Gamma(V, R^2\pi'_*\C_{\tilde{V}}) = \dim \Gamma(U', R^2\pi'_*\C_{\tilde{V}}) = \dim \Gamma(U', \cH).
\end{equation*}
Since $\cH$ is also a local system on the singular locus $\Sigma_V$ of $V$, the argument from before that $\Sigma_{U'}$ intersects every connected component of $\Sigma_V$ implies that the restriction map $\Gamma(V, \cH) \rightarrow \Gamma(U', \cH)$ is injective. Since both spaces have the same dimension, it follows that it is an isomorphism. Thus, by the five lemma, diagram \eqref{PartialsFunct} shows that 
\begin{equation*}
r_{U'}(\C[\epsilon]): \pd_{V}(\C[\epsilon]) \rightarrow \pd_{U'}(\C[\epsilon])
\end{equation*}
is an isomorphism, so $r_{U'}$ is an equivalence.
\end{proof}

\begin{cor}

The maps 
\begin{equation*}
R_X \xrightarrow{\rho^{\sharp}} R_{X'} \xrightarrow{\pi'^{\sharp}} R_Y
\end{equation*}
induced by 
\begin{equation*}
\pd_Y \xrightarrow{\pi'_*} \pd_{X'} \xrightarrow{\rho_*} \pd_X
\end{equation*}
are injections.

\end{cor}

\begin{proof}

By \Cref{AllSmooth}, $R_X$, $R_{X'}$, and $R_Y$ all have the same dimension. Additionally, the argument in part $(iii)$ of \cite[Theorem 5.1]{NaPD1} shows that the maps $R_X \rightarrow R_{X'} \rightarrow R_Y$ all have closed finite fiber. The map $R_X/\ker(\rho^{\sharp}) \rightarrow R_{X'}$ induced by $\rho^{\sharp}$ is injective. This map has closed finite fiber since $\rho^{\sharp}$ does, so we have
\begin{equation*}
\dim R_{X'} \leq \dim R_X/\ker(\rho^{\sharp}) \leq \dim R_X.
\end{equation*}
Since $\dim R_{X'} = \dim R_X$, these must be equalities, so $\ker(\rho^{\sharp})$ is trivial. Thus, $\rho^{\sharp}$ is injective, and the same argument shows the map $R_{X'} \rightarrow R_Y$ is injective. 
\end{proof}

\subsection{Algebraization} \label[subsection]{Algebraization}

The goal of this section is to construct the universal graded deformation of $X'$, and show is it compatible with the universal graded deformations of $X$ and $Y$ constructed in \cite[Section 5.4]{NaPD1}. The work in the previous two subsections implies that $\pd_{X'}$ has a universal formal deformation that is compatible with the universal formal deformations of $X$ and $Y$. Additionally, $X$ has a $\C^*$ action that lifts compatibly to $X'$ and $Y$, as we will show. Using the work in \cite[Section 4]{NaFlops}, \cite[Appendix A]{NaFlops}, \cite[Section 5.4]{NaPD1}, and \cite[Proposition 2.6]{Lo}, we can use this $\C^*$ action to algebraicize the universal formal deformations of $X$, $X'$ and $Y$ in a compatible way. Thus, we will obtain compatible universal graded deformations over full affine spaces of the same dimension as $\pd_{Y}(\C[\epsilon])$. In order to replace the bases of the universal formal deformations with affine spaces, we replace the formal rings that prorepresent $\pd_X$, $\pd_{X'}$, and $\pd_Y$ with the subrings generated by $\C^*$-eigenvectors. Throughout the proof of \Cref{UnivConic} below, we will carefully extract the algebraic data that encodes the universal formal deformations, and we will replace the rings that appear with their subrings generated by $\C^*$-eigenvectors, keeping track of the compatibilities between different spaces as we proceed. 

Since $\pd_X$, $\pd_{X'}$, and $\pd_Y$ are all prorepresentable and unobstructed, $X$, $X'$ and $Y$ have compatible universal formal deformations $\hat{\cX}$, $\hat{\cX'}$, and $\hat{\cY}$ over the formal neighborhoods of zero in their respective tangent spaces, endowed with the natural scheme structure as affine spaces. In particular, the functors $\pd_Y \rightarrow \pd_{X'} \rightarrow \pd_X$ induce maps $\hat{\cY} \rightarrow \hat{\cX'} \rightarrow \hat{\cX}$. We let $\A_Y$, $\A_{X'}$, and $\A_X$ be the affine spaces
\begin{align*}
&\A_{Y} := \Spec \Sym(\pd_{Y}(\C[\epsilon])^*), \\ 
&\A_{X'} := \Spec \Sym(\pd_{X'}(\C[\epsilon])^*), \\
&\A_{X} := \Spec \Sym(\pd_{X}(\C[\epsilon])^*)
\end{align*}
that are naturally associated with $\pd_{Y}(\C[\epsilon])$, $\pd_{X'}(\C[\epsilon])$, and $\pd_{X}(\C[\epsilon])$.

Thus we have the following commutative diagram, where $\hat{V}$ denotes the formal completion of an affine space $V$ at zero.
\begin{equation} \label{FormalUniv}
\begin{tikzcd}[row sep=scriptsize]
	Y && {\hat{\cY}} && {\widehat{\A_Y}} \\
	\\
	{X'} && {\hat{\cX'}} && {\widehat{\A_{X'}}} \\
	\\
	X && {\hat{\cX}} && {\widehat{\A_X}}
	\arrow["{i_{Y}}", hook, from=1-1, to=1-3]
	\arrow["{\pi'}"', from=1-1, to=3-1]
	\arrow[from=1-3, to=1-5]
	\arrow["{\hat{\pi}'}"', from=1-3, to=3-3]
	\arrow["{\pi'_*}"', from=1-5, to=3-5]
	\arrow["{i_{X'}}", hook, from=3-1, to=3-3]
	\arrow["\rho"', from=3-1, to=5-1]
	\arrow[from=3-3, to=3-5]
	\arrow["{\hat{\rho}}"', from=3-3, to=5-3]
	\arrow["{\rho_*}"', from=3-5, to=5-5]
	\arrow["{i_X}", hook, from=5-1, to=5-3]
	\arrow[from=5-3, to=5-5]
\end{tikzcd}
\end{equation}

We recall that $X$ has a $\C^*$-action satisfying \Cref{ConicDef}, which lifts uniquely to $Y$ by \Cref{QFac}. Modifying the argument in \cite[Proposition A.7]{NaFlops}, we can show that the action lifts compatibly to $X'$ as well.

\begin{lemma}

There is a $\C^*$-action on $X'$ such that $\pi': Y \to X$ and $\rho: X' \to X$ are $\C^*$-equivariant.

\end{lemma}

\begin{proof}

As before, we let $U \subset X$ be the union of the open and codimension two symplectic leaves, and we let $U' = \rho^{-1}(U)$ and $\tilde{U} = \pi^{-1}(U)$. In particular, $U'$ is a partial resolution of $U$ given by contracting some of the exceptional divisors of $\pi|_{\tilde{U}}: \tilde{U} \to U$, which is the minimal resolution of $U$ (see \cite[p. 52]{NaPD1}, \cite[Section 3.7]{MMY25}). As stated in step 1 of \cite[Proposition A.7]{NaFlops}, the $\C^*$-action on $X$ restricts to an action on $U$, which lifts to $\tilde{U}$ because it is the minimal resolution. Since the $\C^*$-action preserves the exceptional divisors, the action lifts compatibly to $U'$ as well. Finally, $X' - U'$ has codimension at least two in $X'$, so the argument in \cite[Proposition A.7, Step 1]{NaFlops} carries over to show that the $\C^*$-action on $U'$ extends uniquely to $X'$. 
\end{proof}

Before stating the main theorem we recall the definition of a universal graded deformation, as described in \cite[p. 8]{Lo}. 

\begin{rem}
A similar construction, defined as a universal conic deformation, is described in \cite[Section 2.3]{MN}. 
\end{rem}

\begin{df} \label[df]{GradedDef}

Let $Z$ be a Poisson variety with a $\C^*$ action that rescales the Poisson bracket with weight $-d$ for $d \in \Z_{>0}$. Then a graded deformation of $Z$ over a positively graded ring $B$ is a pair $(Z_B, i)$ where $Z_B$ is a flat Poisson scheme over $B$ equipped with a $\C^*$-action and $i: Z_B \times_{\Spec(B)} 0 \xrightarrow{\sim} Z$ is an isomorphism of the central fiber of $Z_B$ with $Z$, satisfying: 

\begin{enumerate}

\item The $\C^*$-action on $Z_B$ transforms the relative Poisson bracket with weight $-d$.
\item The $\C^*$-action induced from the grading on $\Spec(B)$ is compatible with the $\C^*$ action on $Z_B$.
\item $i$ intertwines the $\C^*$ actions on $Z_B$ and $Z$. 

\end{enumerate}

Furthermore, a graded deformation $\cZ \rightarrow V$ is universal if given any graded deformation $Z_B \rightarrow \Spec(B)$ of $Z$, there is a $\C^*$-equivariant map $ \Spec(B) \rightarrow V$ and an isomorphism of graded deformations $Z_B \xrightarrow{\sim} \Spec(B) \times_{V} \cZ$.

\end{df}

\begin{thm} \label[thm]{UnivConic}

Diagram \eqref{FormalUniv} algebraicizes to the diagram of universal graded deformations 
\begin{equation*}
\begin{tikzcd}[row sep=scriptsize]
	Y && \cY && {\A_Y} \\
	\\
	{X'} && {\cX'} && {\A_{X'}} \\
	\\
	X && \cX && {\A_X}
	\arrow["{i_{Y}}", hook, from=1-1, to=1-3]
	\arrow["{\pi'}"', from=1-1, to=3-1]
	\arrow[from=1-3, to=1-5]
	\arrow["{\tilde{\pi}'}"', from=1-3, to=3-3]
	\arrow["{f'}"', from=1-5, to=3-5]
	\arrow["{i_{X'}}", hook, from=3-1, to=3-3]
	\arrow["\rho"', from=3-1, to=5-1]
	\arrow[from=3-3, to=3-5]
	\arrow["{\tilde{\rho}}"', from=3-3, to=5-3]
	\arrow["g"', from=3-5, to=5-5]
	\arrow["{i_X}", hook, from=5-1, to=5-3]
	\arrow[from=5-3, to=5-5]
\end{tikzcd}
\end{equation*}

\end{thm}

\begin{proof} 

We divide the proof below into four steps. 

\textbf{1} (\textit{Replacing formal spaces with schemes}): Let $A_0 := \limit \Gamma(\hat{\cX}_n, \cO_{\hat{\cX}_n})$ as a topological ring, so that $\hat{\cX} \cong \Spf(A_0)$. Similarly, let 
$$
A'_0 := \limit \Gamma(\hat{\cX}'_n, \cO_{\hat{\cX}'_n}) \hspace{5pt} \text{and} \hspace{5pt} B_0 := \limit \Gamma(\hat{\cY}_n, \cO_{\hat{\cY}_n}),
$$ 
so that we have natural maps $\hat{\cX'} \rightarrow \Spf(A'_0)$ and $\hat{\cY} \rightarrow \Spf(B_0)$.

We can extend diagram \eqref{FormalUniv} to the following diagram of formal spaces.
\begin{equation*}
\begin{tikzcd}[column sep=scriptsize]
	{\hat{\cY}} & {\hat{\cX}'_{B_0} := \hat{\cX'} \times_{\Spf(A'_0)} \Spf(B_0)} & {\Spf(B_0)} & {\widehat{\A_Y} \cong \Spf(R_Y)} \\
	\\
	& {\hat{\cX'}} & {\Spf(A'_0)} & {\widehat{\A_{X'}} \cong \Spf(R_{X'})} \\
	\\
	&& {\hat{\cX} \cong \Spf(A_0)} & {\widehat{\A_X} \cong \Spf(R_X)}
	\arrow["{\hat{\pi}'_{B_0}}", from=1-1, to=1-2]
	\arrow["{\hat{\pi}'}"', from=1-1, to=3-2]
	\arrow["{\hat{\rho}_{B_0}}", from=1-2, to=1-3]
	\arrow[from=1-2, to=3-2]
	\arrow[from=1-3, to=1-4]
	\arrow[from=1-3, to=3-3]
	\arrow["{\pi'_*}"', from=1-4, to=3-4]
	\arrow["{\hat{\rho}_{A'_0}}", from=3-2, to=3-3]
	\arrow["{\hat{\rho}}"', from=3-2, to=5-3]
	\arrow[from=3-3, to=3-4]
	\arrow[from=3-3, to=5-3]
	\arrow["{\rho_*}"', from=3-4, to=5-4]
	\arrow[from=5-3, to=5-4]
\end{tikzcd}
\end{equation*}
We also define $\hat{\pi}_{B_0} := \hat{\rho}_{B_0} \circ \hat{\pi}'_{B_0}$. By \cite[Th\'eor\`eme 5.4.5]{EGA III}, we can algebraicize all of the maps above that are projective maps to the formal spectrum of a ring. Thus, we can algebraicize the maps 
$$
\hat{\rho}_{A'_0}: \hat{\cX'} \rightarrow \Spf(A'_0), \hspace{5pt} \hat{\rho}_{B_0}: \hat{\cX}'_{B_0} \rightarrow \Spf(B_0), \hspace{5pt} \text{and} \hspace{5pt} \hat{\pi}_{B_0}: \hat{\cY} \rightarrow \Spf(B_0).
$$
Let the algebraizations of these maps be 
$$
\cX'_{A'_0} \rightarrow \Spec(A'_0), \hspace{5pt} \cX'_{B_0} := \cX'_{A'_0} \times_{\Spec(A'_0)} \Spec(B_0) \rightarrow \Spec(B_0), \hspace{5pt} \text{and} \hspace{5pt} \cY_{B_0} \rightarrow \Spec(B_0)
$$
respectively. By \cite[Th\'eor\`eme 5.4.1]{EGA III}, the map 
\begin{equation*}
\Hom_{\Spec(B_0)}(\cY_{B_0}, \cX'_{B_0}) \rightarrow \Hom_{\Spf(B_0)}(\hat{\cY}, \hat{\cX}'_{B_0})
\end{equation*}
given by mapping a morphism to its formal completion is a bijection. Thus, we can algebraicize the map $\hat{\cY} \rightarrow \hat{\cX}'_{B_0}$ to a map $\cY_{B_0} \rightarrow \cX'_{B_0}$. Also, base changing the map $\Spec(B_0) \rightarrow \Spec(A'_0)$ clearly algebraicizes $\hat{\cX}'_{B_0} \rightarrow \hat{\cX}'$ to a map $\cX'_{B_0} \rightarrow \cX'_{A'_0}$. Thus, we obtain the following diagram of schemes, where we also define $\pi_{B_0} := \pi'_{B_0} \circ \rho_{B_0}$.
\begin{equation} \label{QuarterAlg}
\begin{tikzcd}[row sep=scriptsize]
	{\cY_{B_0}} & {\cX'_{B_0} := \cX'_{A'_0} \times_{\Spec(A'_0)} \Spec(B_0)} & {\Spec(B_0)} & {\Spec(R_Y)} \\
	\\
	& {\cX'_{A'_0}} & {\Spec(A'_0)} & {\Spec(R_{X'})} \\
	\\
	&& {\Spec(A_0)} & {\Spec(R_X)}
	\arrow["{\pi'_{B_0}}", from=1-1, to=1-2]
	\arrow[from=1-1, to=3-2]
	\arrow["{\rho_{B_0}}", from=1-2, to=1-3]
	\arrow[from=1-2, to=3-2]
	\arrow[from=1-3, to=1-4]
	\arrow[from=1-3, to=3-3]
	\arrow[from=1-4, to=3-4]
	\arrow["{\rho_{A'_0}}", from=3-2, to=3-3]
	\arrow[from=3-2, to=5-3]
	\arrow[from=3-3, to=3-4]
	\arrow[from=3-3, to=5-3]
	\arrow[from=3-4, to=5-4]
	\arrow[from=5-3, to=5-4]
\end{tikzcd}
\end{equation}

\textbf{2} (\textit{Replacing the bases with affine spaces}): The next step is to use the tools from \cite[Appendix A]{NaFlops} to further algebraicize the above diagram to a new diagram where we have replaced the rings with the subrings generated by $\C^*$-eigenvectors. As a technical point, in order to use the tools in \cite[Appendix A]{NaFlops}, we have to replace $A_0$, $A'_0$, and $B_0$ with the local complete $\C$-algebras given by taking completions at the origin. We call these rings $A$, $A'$, and $B$ respectively, and they all have residue field $\C$, and good $\C^*$-actions in the sense of \cite[Appendix A]{NaFlops}. We base change the rest of \eqref{QuarterAlg} to obtain diagram \eqref{HalfAlg} below. When we construct an algebraization of \eqref{HalfAlg}, it will also be an algebraization of \eqref{QuarterAlg} (see e.g. the proof of \cite[Lemma 22]{NaFlops}).
\begin{equation} \label{HalfAlg}
\begin{tikzcd}[row sep=scriptsize]
	{\cY_B} & {\cX'_B := \cX'_{A'} \times_{\Spec(A')} \Spec(B)} & {\Spec(B)} & {\Spec(R_Y)} \\
	\\
	& {\cX'_{A'}} & {\Spec (A')} & {\Spec(R_{X'})} \\
	\\
	&& {\Spec(A)} & {\Spec(R_X)}
	\arrow[from=1-3, to=3-3]
	\arrow[from=3-3, to=5-3]
	\arrow[from=1-3, to=1-4]
	\arrow[from=1-4, to=3-4]
	\arrow[from=3-4, to=5-4]
	\arrow[from=3-3, to=3-4]
	\arrow[from=5-3, to=5-4]
	\arrow[from=3-2, to=5-3]
	\arrow["{\rho_{A'}}", from=3-2, to=3-3]
	\arrow["{\rho_B}", from=1-2, to=1-3]
	\arrow[from=1-2, to=3-2]
	\arrow["{\pi'_B}", from=1-1, to=1-2]
	\arrow[from=1-1, to=3-2]
\end{tikzcd}
\end{equation}

 We recall from the proof of \cite[Lemma 22]{NaFlops} that $R_Y \cong \C[[t_i, \dots , t_m]]$, where each of the $t_i$ are $\C^*$-eigenvectors of the same weight $d$ as the symplectic form on $X$. Thus we can apply  \cite[Lemma A.2]{NaFlops} to $R_Y$, as well as to $R_X$ and $R_{X'}$, since they are subrings of $R_Y$. In particular, since $R_Y$ is the ring of functions of the formal neighborhood of zero in $\A_Y$, \cite[Lemma A.2]{NaFlops} lets us recover $\A_Y$ as the spectrum of the subring of $R_Y$ generated by $\C^*$-eigenvectors. Applying the same argument to $R_X$ and $R_{X'}$, we algebraicize the right column of diagram \eqref{HalfAlg} to a sequence of maps
\begin{equation*}
\A_Y \to \A_{X'} \to \A_X,
\end{equation*}
which are also finite and surjective by \cite[Lemma A.4]{NaFlops}.

Similarly, we can algebraicize the second column of diagram \eqref{HalfAlg} on the right by replacing $B$, $A'$, and $A$ by their subrings generated by $\C^*$-eigenvectors, which we will denote $T$, $S'$, and $S$ respectively. Since the maps $\pi_B$, $\rho_B$, and $\rho_{A'}$ are projective maps to affine schemes, by \cite[Proposition A.5]{NaFlops} we can algebraicize them to projective maps
\begin{equation*}
\pi_T: \cY \to \Spec(T), \hspace{10pt} \rho_T: \cX'_T \to \Spec(T), \hspace{10pt} \rho_{S'}: \cX' \to \Spec(S').
\end{equation*}
We note that since $\cX'_B = \cX'_{A'} \times_{\Spec(A')} \Spec(B)$, it follows from the construction in \cite[Proposition A.5]{NaFlops} that
$$
\cX'_T \cong \cX'_{S'} \times_{\Spec(S')} \Spec(T),
$$
and that $\rho_T$ is the natural projection map to $\Spec(T)$. Lastly, we need to algebraicize the map $\pi'_B: \cY_B \rightarrow \cX'_B$. We do this by showing that $\cY$ and $\cX'_T$ have compatible line bundles, as we explain below. 

\textbf{3} (\textit{Algebraicizing $\cY_B \rightarrow \cX'_B$}): Since the map $\rho_B: \cX'_B \rightarrow \Spec(B)$ is projective, we pick a $\rho_B$-ample line bundle $\cL_{\cX'_B}$. Thus, we have  
$$
\cX'_B = \PROJ_B \oplus_{i \geq 0} B'_i, \hspace{10pt} B'_i := \Gamma(\cX'_B, \cL_{\cX'_B}^{\otimes i}).
$$
We can pick $\cL_{\cX'_B}$ so that $B'_* := \oplus_{i \geq 0} B'_i$ is generated in degree 1. Note that by the construction in \cite[Proposition A.5]{NaFlops}, $\cX'_T$ is given by 
$$
\cX'_T := \PROJ_B \oplus_{i \geq 0} B'_{i, T},
$$
where $B'_{i, T}$ is the subspace of $B'_i$ generated by $\C^*$-eigenvectors. Similarly, we pick a $\pi_B$-ample line bundle $\cL_{\cY_B}$, so that
$$
\cY_B = \PROJ_B \oplus_{i \geq 0} B_i, \hspace{10pt} B_i := \Gamma(\cY_B, \cL_{\cY_B}^{\otimes i}).
$$ 
Additionally, we require that $B_* := \oplus_{i \geq 0} B_i$ is generated in degree 1, and we have
$$
\cY_T := \PROJ_B \oplus_{i \geq 0} B_{i, T},
$$
where $B_{i, T}$ is the subspace of $B_i$ generated by $\C^*$-eigenvectors. 

By \cite[\href{https://stacks.math.columbia.edu/tag/01O4}{Lemma 01O4}]{Stacks}, the map 
$$
\pi'_B: \cY_B \to \cX'_B
$$
is determined uniquely by the line bundle $\pi'^*_B \cL_{\cX'_B} \in \Pic(\cY_B)$ and the natural surjection 
\begin{equation} \label{F1B}
\pi^*_B B'_1 = \pi^*_B \Gamma(\cX'_B, \cL_{\cX'_B}) \to \pi'^*_B \cL_{\cX'_B},
\end{equation}
which is given by pulling back the counit map $\rho^*_B B'_1 \to \cL_{\cX'_B}$ along $\pi'^*_B$. Note that this map is $\C^*$-equivariant by construction. Additionally, $\pi^*_B B'_1$ is the sheaf associated to the graded $B_*$-module $\oplus_{i \geq 0} \ B_i \otimes_B B'_1$. By \cite[Proposition A.6]{NaFlops}, there is a coherent sheaf $\cF_T$ on $\cY_T$ such that $\cF_T \otimes_T B \cong \pi^*_B B'_1$. The construction in \cite[Proposition A.6]{NaFlops} replaces each $B_i \otimes_B B'_1$ with the subspace $(B_i \otimes_B B'_1)_T$ generated by $\C^*$-eigenvectors. We will now explain why the natural map
\begin{equation} \label{C^*eigen}
B_{i, T} \otimes_T B'_{1, T} \to (B_i \otimes_B B'_1)_T 
\end{equation}
is an isomorphism. By \cite[Theorem 9.3A]{Ha}, the functor $ - \otimes_T B$ is exact. Additionally, applying this functor to \eqref{C^*eigen} yields an isomorphism since 
\begin{align*}
(B_{i, T} \otimes_T B'_{1, T}) \otimes_T B &\cong B_{i, T} \otimes_T B \otimes_B B'_{1, T} \otimes_T B \\
&\cong (B_{i, T} \otimes_T B) \otimes_B (B'_{1, T} \otimes_T B) \\
&\cong B_i \otimes_B B'_1 \cong (B_i \otimes_B B'_1)_T \otimes_T B.
\end{align*}
Thus the cokernel $M$ of \eqref{C^*eigen} satisfies $M \otimes_T B = 0$. As explained in the proof of \cite[Proposition A.5]{NaFlops}, this means that the support of $M$ is a closed subset of $\Spec(T)$ that does not contain the origin. Since the $\C^*$-action is good, this means that $M$ is trivial and \eqref{C^*eigen} is surjective. Since $ - \otimes_T B$ is exact, we can apply the same argument to the kernel of \eqref{C^*eigen}, and we obtain that the map is indeed an isomorphism. 

Note that the sheaf on $\cY_T$ given as a graded module by $\oplus_{i \geq 0} B_{i, T} \otimes_T B'_{1, T}$ is precisely $\pi^*_T B'_{1, T}$. By \cite[Proposition A.6]{NaFlops}, we also obtain a line bundle $\ol{\cL} \in \Pic(\cY_T)$ such that 
$$
\ol{\cL} \otimes_T B \cong \pi'^*_B \cL_{\cX'_B}.
$$
Because \eqref{F1B} is $\C^*$-equivariant, $\C^*$-eigenvectors will map to $\C^*$-eigenvectors, and thus it follows from the construction in \cite[Proposition A.6]{NaFlops} that we obtain a $\C^*$-linearized map 
\begin{equation} \label{F1T}
\pi^*_T B'_{1, T} \to \ol{\cL}.
\end{equation}
This map becomes surjective after applying $ - \otimes_T B$. Thus, the cokernel $\cM$ of this map satisfies $\cM \otimes_T B = 0$. Using the same argument as above, this implies that $\cM$ is the trivial sheaf and \eqref{F1T} is surjective. By \cite[\href{https://stacks.math.columbia.edu/tag/01O4}{Lemma 01O4}]{Stacks} again, this produces a map 
$$
\pi'_T: \cY \to \cX'_T
$$
over $\Spec(T)$, which algebraicizes $\pi'_B$. Finally, since $\pi_T$ and $\rho_T$ are projective, $\pi'_T$ is projective as well. Recapping, we have algebraicized $\cY_B \rightarrow \cX'_B \rightarrow \Spec(B)$ to a sequence of projective and birational maps
\begin{equation*}
\cY \rightarrow \cX'_T \rightarrow \Spec(T).
\end{equation*}
\textbf{4} (\textit{Final points}): Combining these algebraizations we obtain a final diagram 
\begin{equation} \label{FinalAlg}
\begin{tikzcd}[column sep=scriptsize]
	\cY & {\cX'_T = \cX' \times_{\Spec(S')} \Spec(T)} & {\Spec(T)} & {\A_Y} \\
	\\
	& {\cX'} & {\Spec (S')} & {\A_{X'}} \\
	\\
	&& {\cX \cong \Spec(S)} & {\A_{X}}
	\arrow[from=1-1, to=1-2]
	\arrow[from=1-1, to=3-2]
	\arrow[from=1-2, to=1-3]
	\arrow[from=1-2, to=3-2]
	\arrow[from=1-3, to=1-4]
	\arrow[from=1-3, to=3-3]
	\arrow[from=1-4, to=3-4]
	\arrow[from=3-2, to=3-3]
	\arrow[from=3-2, to=5-3]
	\arrow[from=3-3, to=3-4]
	\arrow[from=3-3, to=5-3]
	\arrow[from=3-4, to=5-4]
	\arrow[from=5-3, to=5-4]
\end{tikzcd}
\end{equation}
where the map $\cX'_T \rightarrow \cX'$ is defined by base change and the diagonal maps are simply defined by composition. Finally, the method in \cite[Lemma 22]{NaFlops} extends to give $\cY$, $\cX'$, and $\cX$ relative Poisson structures over $\A_Y$, $\A_{X'}$, and $\A_X$ respectively. In particular, it is clear that $B$ has a Poisson structure over $R_Y$, and this induces a Poisson structure of $T$ over $\C[\A_Y]$. Similarly, we obtain Poisson structures of $\Spec(S')$ over $\A_{X'}$ and $\cX \cong \Spec(S)$ over $\A_X$. Because the maps $\cY \rightarrow \cX'_T \rightarrow \Spec(T)$ and $\cX' \rightarrow \Spec(S')$ are birational, we can use the relative Poisson structures on $\Spec(T)$ and $\Spec(S')$ to induce relative Poisson structures on $\cY$, $\cX'_T$, and $\cX'$.

Restricting to the right column and the diagonal sequence in diagram \eqref{FinalAlg} we have the desired algebraization. All that is left is to note that it is shown in \cite[Proposition 2.6]{Lo} and \cite[Proposition 2.12]{Lo} that $\cY \to \A_Y$ and $\cX \rightarrow \A_X$ are universal graded deformations, and the argument in \cite[Proposition 2.6]{Lo} carries over to show that $\cX' \rightarrow \A_{X'}$ is a universal graded deformation as well.
\end{proof}

\section{Kleinian Singularities} \label[section]{Klein}

\subsection{Structure of Partial Resolutions} \label[subsection]{Klein.1}

Let $\pi: \tilde{S} \rightarrow S$ be the unique crepant resolution of a Kleinian singularity $S$ corresponding to the Lie algebra $\mathfrak{g}$ with Cartan subalgebra $\mathfrak{h}$ and Weyl group $W$.  Recall that by \cite[Proposition 3.1]{NaPD1} we have universal deformations of $S$ and $\tilde{S}$, obtained from the Grothendieck-Springer resolution of $\mathfrak{g}$: 
\begin{equation} \label{fullKlein}
\begin{tikzcd}[row sep=scriptsize]
	{\tilde{S}} && {\tilde{\cS}} && {\A_{\tilde{S}} \cong \mathfrak{h}} \\
	\\
	S && \cS && {\A_{S} \cong \mathfrak{h}/W}
	\arrow[hook, from=1-1, to=1-3]
	\arrow[from=1-1, to=3-1]
	\arrow[from=1-3, to=1-5]
	\arrow[from=1-3, to=3-3]
	\arrow[from=1-5, to=3-5]
	\arrow[hook, from=3-1, to=3-3]
	\arrow[from=3-3, to=3-5]
\end{tikzcd}
\end{equation}
where $\mathfrak{h} \rightarrow \mathfrak{h}/W$ is the standard quotient map. In particular, $\cS$ can be constructed explicitly by taking a Slodowy slice to the subregular nilpotent orbit of $\mathfrak{g}$ and giving it the relative Poisson structure over $\mathfrak{h}/W$ induced from the Kostant-Kirillov 2-form. Furthermore, we can construct $\tilde{\cS}$ by taking the preimage of $\cS$ under the Grothendieck-Springer map $\tilde{\mathfrak{g}} \rightarrow \mathfrak{g}$. As explained in \Cref{Algebraization}, we can also equip $\mathfrak{h}$ and $\mathfrak{h}/W$ with appropriate $\C^*$-actions, making the above diagram a diagram of universal graded deformations. 

In particular, since $\A_{\tilde{S}} \cong H^2(\tilde{S}, \C)$ is the base of the universal graded deformation of $\tilde{S}$, we have isomorphisms $\mathfrak{h} \cong \Pic(\tilde{S}) \otimes \C \cong H^2(\tilde{S}, \C)$. In order to translate the base of our universal deformation to $H^2(\tilde{S}, \C)$ we use the period map
\begin{equation} \label{PeriodMap}
\mathfrak{h} \rightarrow \Pic(\tilde{S}) \otimes \C \cong H^2(\tilde{S}, \C),
\end{equation}
which we will describe explicitly below. This map is described in \cite[Proposition 3.2]{NaPD1} and induces an action of $W$ on $H^2(\tilde{S}, \C)$. Although we do not strictly need this, it is shown in \cite{Ya} that this $W$-action agrees with the standard monodromy action on $H^2(G/B, \C) \cong H^2(\tilde{S}, \C)$, where $G$ is a Lie group of the same type as the Kleinian singularity. 

Recall that the exceptional fiber $\pi^{-1}(0)$ of $\pi$ is a union of projective lines in the shape of the Dynkin diagram $D_S$ corresponding to $S$. Thus the group $H_2(\tilde{S}, \Z)$ is generated by the fundamental classes $[C_i]$ of the projective lines comprising the exceptional locus. Furthermore, the lattice
\begin{equation*}
H_2(\tilde{S}, \Z) = \{\sum a_i[C_i] \ | \ a_i \in \Z \}
\end{equation*}
has a symmetric bilinear form $( -, - )$ given by the negative of the intersection product. This product endows $H_2(\tilde{S}, \Z)$ with the structure of a root system of the same type as the Dynkin diagram $D_S$. Since the $[C_i]$ are in bijection with the nodes of $D$, we have an isomorphism 
\begin{equation*}
\mathfrak{h} \rightarrow H_2(\tilde{S}, \C), \hspace{10pt} \alpha^{\vee}_i \mapsto [C_i]
\end{equation*}
where $\alpha^{\vee}_i$ is the simple coroot corresponding to the node $i$. Thus, we obtain a natural isomorphism 
\begin{equation} \label{dualPeriod}
\mathfrak{h}^* \xrightarrow{\sim} H^2(\tilde{S}, \C) \cong (H_2(\tilde{S}, \C))^*.
\end{equation}

In particular, this isomorphism identifies the dual fundamental classes $[C_i]^{\vee} \in H^2(\tilde{S}, \C)$ with the fundamental weights in $\mathfrak{h}^*$ corresponding to the same node of the Dynkin diagram. Thus, the dual fundamental classes $[C_i]^{\vee}$ generate the standard dominant Weyl chamber given by elements of $\mathfrak{h}^*$ that have non-negative pairing with each simple coroot. We note that this map is obtained from the map in \cite[Proposition 4.7]{LMM} by base change to $\C$.

Additionally, via the natural isomorphism $H^2(G/B, \C) \rightarrow H^2(\tilde{S}, \C)$ for Lie algebras with simply laced Dynkin diagrams, the isomorphism above is that standard map 
\begin{equation*}
\mathfrak{h}^* \xrightarrow{\sim} \Pic(G/B) \otimes \C \cong H^2(G/B, \C),
\end{equation*}
described for example in \cite[(P.3)]{NaBir} and \cite[Section 6.1]{CG}.

It is shown in the proof of \cite[Proposition 3.2]{NaPD1} that the period map \eqref{PeriodMap} is given by 
\begin{equation*}
\mathfrak{h} \rightarrow \mathfrak{h}^* \rightarrow H^2(\tilde{S}, \C),
\end{equation*}
where the first map is the isomorphism induced by the Killing form, and the second is given by \eqref{dualPeriod}. In particular, this isomorphism identifies the dual fundamental classes $[C_i]^{\vee} \in H^2(\tilde{S}, \C)$ with the corresponding fundamental coweights in $\mathfrak{h}$. This description immediately gives us the following result

\begin{cor} \label[cor]{subsetInv}

Let $J \subset I$ be a subset of the set $I$ of nodes of the Dynkin diagram $D_S$. Then let $V$ be the subspace of $H^2(\tilde{S}, \C)$ generated by $\{[C_j]^{\vee} \ | \ j \in J\}$ and $W'$ be the subgroup of $W$ generated by the simple generators corresponding to $I - J$. Then $W'$ is the stabilizer of $V$ in $H^2(\tilde{S}, \C)$ and $V \cong H^2(\tilde{S}, \C)^{W'}$.

\end{cor}

Given a symplectic singularity $Z$ and a subset $E \subset Z$ we let $(Z, E)$ denote the complex analytic germ of $E$ in $Z$. Additionally, deformations of germs $(Z, E)$ will be taken in the complex analytic topology. In particular, let $(S, 0)$ be the germ of $S$ at zero and let $(\tilde{S}, F)$ be the germ of $\tilde{S}$ at $F := \pi^{-1}(0)$. It is proven in \cite[Proposition 3.1]{NaPD1} that the universal deformations of $(S, 0)$ and $(\tilde{S}, F)$ are given by the universal deformations of $S$ and $\tilde{S}$ respectively, which are described above. In particular, this means that we have the following commutative diagram, where the horizontal maps are given by restriction. 
\begin{equation} \label{LocDiagram}
\begin{tikzcd}
	{\pd_{\tilde{S}}} && {\pd_{(\tilde{S}, F)}} \\
	\\
	{\pd_{S}} && {\pd_{(S, 0)}}
	\arrow["\sim", from=1-1, to=1-3]
	\arrow["{\pi_*}"', from=1-1, to=3-1]
	\arrow["{(\pi|_{(\tilde{S}, F)})_*}"', from=1-3, to=3-3]
	\arrow["\sim", from=3-1, to=3-3]
\end{tikzcd}
\end{equation}
The above diagram shows that all the data of the Poisson deformations is determined locally near the singularity. This is useful for the study of partial resolutions of Kleinian singularities because they all have finitely many singular points and locally around each singular point they are isomorphic to Kleinian singularities corresponding to subsets of the Dynkin diagram of the original Kleinian singularity. We now study these in more detail.

Let $\rho: S' \rightarrow S$ be a partial resolution of $S$. Since $\tilde{S}$ is the unique crepant resolution of $S$, we have a map $\pi': \tilde{S} \rightarrow S'$. Let $D_{S'} \subset D_S$ be the subset of the Dynkin diagram $D_S$ corresponding to the projective lines contracted by $\pi'$.

By \cite[p. 16-18]{Nak} (see also \cite{Slo}) the singular locus $\Sigma$ of $S'$ is a set of points in bijection with the connected components $D_i$ of $D_{S'}$. Furthermore, the singularity at a point $x_i \in \Sigma$ corresponding to a component $D_i$ is a Kleinian singularity $S_i$ of the same type as $D_i$. For the rest of the section let $\mathfrak{h}_i$ and $W_i$ be the Cartan and Weyl group associated to this Kleinian singularity and let $T'_i$ denote the complex analytic germ $(S', x_i)$. Thus, we have a restriction map 
\begin{equation} \label{PartGerm}
\prod_{x_i \in \Sigma} r_{T'_i}: \pd_{S'} \rightarrow \pd_{(S', \Sigma)} = \prod_{x_i \in \Sigma}\pd_{T'_i} \cong \prod_{x_i \in \Sigma}\pd_{(S_i, 0)}.
\end{equation}
Additionally if we let $\tilde{F_i} := \pi'^{-1}(x_i)$ and let $\tilde{T}_i$ be the complex analytic germ $(\tilde{S}, \tilde{F_i})$, we obtain the sequence 
\begin{equation} \label{FullGerm}
\prod_{x_i \in \Sigma} r_{\tilde{T}_i}: \pd_{\tilde{S}} \rightarrow \pd_{(\tilde{S}, \pi'^{-1}(\Sigma))} = \prod_{x_i \in \Sigma}\pd_{\tilde{T}_i} \cong \prod_{x_i \in \Sigma}\pd_{(\tilde{S}_i, F_i)}
\end{equation}
where $F_i$ is the exceptional fiber of the crepant resolution $\tilde{S}_i$ of $S_i$. 

The kernel of the restriction functor in diagram \eqref{PartGerm} is the subfunctor $\pd_{lt, S'}$ of $\pd_{S'}$ given by locally trivial deformations, while the kernel of restriction functor in diagram \eqref{FullGerm} is given by the subfunctor $\pd_{\tilde{S}}^{\pi'}$ of $\pd_{\tilde{S}}$. This is defined in \cite[p. 348]{NaMori}, and we reproduce the definition below:

\begin{df} \label[df]{LocSub}

Given a sequence of partial resolutions $X' \xrightarrow{\mu} X'' \rightarrow X$ of a conical affine symplectic singularity $X$, the subfunctor $\pd_{X'}^{\mu} \subset \pd_{X'}$ is given by letting $\pd_{X'}^{\mu}(A)$ be the set of Poisson deformations over $A$ that are trivial on the preimage of a sufficiently small analytic neighborhood of every $x \in X''$, for any $A \in \Art_{\C}$.

\end{df}

In particular, the proof of \cite[Lemma 3]{NaMori} also tells us the following about this functor.

\begin{lemma}[{\cite[p. 348]{NaMori}}] \label[lemma]{PDPartial}

Using the notation of \Cref{LocSub}, if $X'$ is a symplectic resolution of $X$, then the functor $\pd_{X'}^{\mu}$ is prorepresentable and unobstructed, with tangent space $\pd_{X'}^{\mu}(\C[\epsilon]) \cong H^2(X'', \C)$.

\end{lemma}

Assembling the above information we obtain the following diagram of functors:
\begin{equation} \label{PartKlein}
\begin{tikzcd} [row sep=scriptsize]
	{\pd_{\tilde{S}}^{\pi'}} && {\pd_{\tilde{S}}} && {\displaystyle \prod_{x_i \in \Sigma}\pd_{\tilde{T}_i}} \\
	\\
	{\pd_{lt, S'}} && {\pd_{S'}} && {\displaystyle \prod_{x_i \in \Sigma}\pd_{T'_i}}
	\arrow["{\prod r_{\tilde{T}_i}}", from=1-3, to=1-5]
	\arrow["{\pi_*'}"', from=1-3, to=3-3]
	\arrow["{ \prod (\pi'|_{\tilde{T}_i})_*}"', from=1-5, to=3-5]
	\arrow["{ \prod r_{T'_i}}", from=3-3, to=3-5]
	\arrow[hook, from=3-1, to=3-3]
	\arrow[hook, from=1-1, to=1-3]
	\arrow["{\pi_{*, lt}'}"', from=1-1, to=3-1]
\end{tikzcd}
\end{equation}
Recall that for every $A \in \Art_{\C}$, each $\tilde{\cS} \in \pd_{\tilde{S}}^{\pi'}(A)$ is trivial on the preimage of some neighborhood of each $x_i$. Thus, $\pi_*'(A)(\tilde{\cS})$ is trivial on some neighborhood of each $x_i$, so it is an element of $\pd_{lt, S'}(A)$. This implies that for every $A \in \Art_{\C}$, $\pi_*'(A)$ maps $\pd_{\tilde{S}}^{\pi'}(A)$ to $\pd_{lt, S'}(A)$, inducing the functor $\pi_{*, lt}'$ in the above diagram. We need to prove that this map induces an isomorphism on first order deformations. We prove a more general statement below that will also be needed in the next section. 

\begin{prop} \label[prop]{KerIsoPartial}

Let $Y \xrightarrow{\pi'} X' \xrightarrow{\rho} X''$ be a sequence of crepant partial resolutions of a conical affine symplectic singularity, with $\pi := \rho \circ \pi'$. Furthermore, let $U \subset X''$ be the union of the open and codimension two symplectic leaves and let $U' := \rho^{-1}(U)$ and $\tilde{U} := \pi^{-1}(U)$. Then the pushforward maps 
\begin{equation*}
\pd_{\tilde{U}}^{\pi}(\C[\epsilon]) \xrightarrow{\pi'_*(\C[\epsilon])} \pd_{U'}^{\rho}(\C[\epsilon]) \xrightarrow{\rho_*(\C[\epsilon])} \pd_{lt, U}(\C[\epsilon])
\end{equation*}
are both isomorphisms. 

\end{prop}

\begin{proof}

Let $i: U_0 \subset U$ be the inclusion of the smooth locus of $U$, which we also treat as a subset of $U'$ via $\rho$. Since $\rho$ is an isomorphism over $U_0$ we have the following commutative diagram
\begin{equation} \label{rtoU0}
\begin{tikzcd}[row sep=scriptsize]
	{\pd_{U'}^{\rho}(\C[\epsilon])} \\
	&& {\pd_{U_0}(\C[\epsilon])} \\
	{\pd_{lt, U}(\C[\epsilon])}
	\arrow["{r^{U'}_{U_0}}", from=1-1, to=2-3]
	\arrow["{\rho_*(\C[\epsilon])}"', from=1-1, to=3-1]
	\arrow["{r^{U}_{U_0}}"', from=3-1, to=2-3]
\end{tikzcd}
\end{equation}
where the maps from left to right are restriction maps. We shall show that $r^{U'}_{U_0}$ is injective. We recall from the argument in \Cref{rU'equiv} (see \eqref{ltFunct}) that we have a commutative diagram
\begin{equation} \label{U'toU0}
\begin{tikzcd}
	{H^1(U', \cO_{U'})} & {\pd_{lt, U'}(\C[\epsilon])} & {H^2(U', \C)} & {H^2(U', \cO_{U'})} \\
	\\
	{H^1(U_0, \cO_{U_0})} & {\pd_{U_0}(\C[\epsilon])} & {H^2(U_0, \C)} & {H^2(U_0, \cO_{U_0})}
	\arrow[from=1-1, to=1-2]
	\arrow[from=1-1, to=3-1]
	\arrow[from=1-2, to=1-3]
	\arrow["{r^{U'}_{U_0}}"', from=1-2, to=3-2]
	\arrow[from=1-3, to=1-4]
	\arrow["{i^*}"', from=1-3, to=3-3]
	\arrow[from=1-4, to=3-4]
	\arrow[from=3-1, to=3-2]
	\arrow[from=3-2, to=3-3]
	\arrow[from=3-3, to=3-4]
\end{tikzcd}
\end{equation}
Additionally, we recall from \eqref{H12U'} that $H^1(U', \cO_{U'}) = H^2(U', \cO_{U'}) = 0$. By the argument in \cite[Lemma 3]{NaMori}, the Leray spectral sequence shows that the inclusion map 
$$
\pd_{U'}^{\rho}(\C[\epsilon]) \hookrightarrow \pd_{lt, U'}(\C[\epsilon])
$$
is given by the pullback in cohomology $\rho^*: H^2(U, \C) \to H^2(U', \C)$. Combining this with the middle columns of \eqref{U'toU0}, we obtain the commutative diagram 
\begin{equation} \label{PDtoH2}
\begin{tikzcd}
	{\pd_{U'}^{\rho}(\C[\epsilon])} & {H^2(U, \C)} \\
	\\
	{\pd_{lt, U'}(\C[\epsilon])} & {H^2(U', \C)} \\
	\\
	{\pd_{U_0}(\C[\epsilon])} & {H^2(U_0, \C)}
	\arrow["\sim", from=1-1, to=1-2]
	\arrow[hook, from=1-1, to=3-1]
	\arrow["{\rho^*}"', from=1-2, to=3-2]
	\arrow["\sim", from=3-1, to=3-2]
	\arrow["{r^{U'}_{U_0}}"', from=3-1, to=5-1]
	\arrow["{i^*}"', from=3-2, to=5-2]
	\arrow[from=5-1, to=5-2]
\end{tikzcd}
\end{equation}
Both of the maps on the right column of \eqref{PDtoH2} are pullback maps in cohomology, and the composition $\rho \circ i$ is the inclusion of $U_0$ into $U$. Thus, the composition of the maps on the right is the pullback $H^2(U, \C) \to H^2(U_0, \C)$ in cohomology, which is an isomorphism by the argument in \cite[Lemma 2.8]{Lo}. In particular, the composition $\pd_{U'}^{\rho}(\C[\epsilon]) \to H^2(U_0, \C)$ is an isomorphism, so by the commutativity of \eqref{PDtoH2} the composition
$$
r^{U'}_{U_0}: \pd_{U'}^{\rho}(\C[\epsilon]) \to \pd_{U_0}(\C[\epsilon])
$$ 
is an injection. Thus, it follows from the commutativity of \eqref{rtoU0} that 
$$
\rho_*(\C[\epsilon]): \pd_{U'}^{\rho}(\C[\epsilon]) \to \pd_{lt, U}(\C[\epsilon])
$$ 
is an injection as well. By considering the case where $\pi'$ is the identity, the same argument shows that 
$$
\pi_*(\C[\epsilon]): \pd_{\tilde{U}}^{\pi}(\C[\epsilon]) \to \pd_{lt, U}(\C[\epsilon])
$$
is injective. Thus, we obtain a sequence of injections
\begin{equation*}
\pd_{\tilde{U}}^{\pi}(\C[\epsilon]) \xrightarrow{\pi'_*(\C[\epsilon])} \pd_{U'}^{\rho}(\C[\epsilon]) \xrightarrow{\rho_*(\C[\epsilon])} \pd_{lt, U}(\C[\epsilon])
\end{equation*}
Applying \Cref{ltpd} and the argument for \Cref{PDPartial} in \cite[p. 348]{NaMori}, we have 
\begin{equation*}
\pd_{\tilde{U}}^{\pi}(\C[\epsilon]) \cong \pd_{lt, U}(\C[\epsilon]) \cong H^2(U, \C),
\end{equation*}
so the above injections must both be isomorphisms.
\end{proof}

\begin{cor} \label[cor]{LeftProSmooth}

Using the notation from \Cref{KerIsoPartial}, all of the functors in the sequence
\begin{equation*}
\pd_{\tilde{U}}^{\pi} \xrightarrow{\pi'_*} \pd_{U'}^{\rho} \xrightarrow{\rho_*} \pd_{lt, U}
\end{equation*}
are prorepresentable and unobstructed and the maps between them are equivalences.

\end{cor}

\begin{proof}

The prorepresentability of these functors follows from the work in \Cref{ProDef}. Additionally, the argument for \Cref{PDPartial} proves that $\pd_{\tilde{U}}^{\pi}$ is unobstructed. Finally, \Cref{KerIsoPartial} shows that the maps $\pi'_*$ and $\rho_*$ induce isomorphisms on first order deformations. This is only possible if all of the prorepresenting rings are regular and the induced maps between them are isomorphisms. 
\end{proof}

Next we are going to study the rows of diagram \eqref{PartKlein}. In particular, the following result will show that they induce short exact sequences on first order deformations.

\begin{lemma} \label[lemma]{KleinSurj}

The maps $\prod r_{\tilde{T}_i}$ and $\prod r_{T'_i}$ in diagram \eqref{PartKlein} are both surjective. 

\end{lemma}

\begin{proof}

We begin with $\prod r_{T'_i}$. Considering the bottom row of diagram \eqref{PartKlein} and taking first order deformations we obtain the exact sequence
\begin{equation*}
0 \rightarrow \pd_{lt, S'}(\C[\epsilon]) \rightarrow \pd_{S'}(\C[\epsilon]) \rightarrow \prod_{x_i \in \Sigma} \pd_{T'_i}(\C[\epsilon]).
\end{equation*}
By \Cref{AllSmooth} we have $\dim \pd_{S'}(\C[\epsilon]) = \dim H^2(\tilde{S}, \C)$. Furthermore, by \Cref{Partlt} we have $\dim \pd_{lt, S'}(\C[\epsilon]) = \dim H^2(S', \C)$. Finally by \eqref{PartGerm}, \eqref{LocDiagram}, and \eqref{fullKlein}, the dimension of $\pd_{T'_i}(\C[\epsilon])$ is equal to the number of projective lines contracted to $x_i$ by $\pi'$. The sum of these dimensions over $x_i \in \Sigma$ is the total number of projective lines contracted by $\pi'$, which is precisely $\dim H^2(\tilde{S}, \C) - \dim H^2(S', \C)$. Thus, the above sequence must be short exact. 

Moving on to $\prod r_{\tilde{T}_i}$, we consider the top row of diagram \eqref{PartKlein} and take first order deformations to obtain
\begin{equation*}
0 \rightarrow \pd_{\tilde{S}}^{\pi'}(\C[\epsilon]) \rightarrow \pd_{\tilde{S}}(\C[\epsilon]) \rightarrow \prod \pd_{\tilde{T}_i}(\C[\epsilon]).
\end{equation*}
Similarly to above, we have $\dim \pd_{\tilde{S}}(\C[\epsilon]) = \dim H^2(\tilde{S}, \C)$ by diagram \eqref{fullKlein} and $\dim \pd_{\tilde{S}}^{\pi'}(\C[\epsilon]) = \dim H^2(S', \C)$ by \Cref{PDPartial}. Finally, we can also apply \eqref{PartGerm}, \eqref{LocDiagram}, and \eqref{fullKlein} to show that $\sum_{x_i \in \Sigma} \dim \pd_{\tilde{T}_i}(\C[\epsilon])$ is the number of projective line contracted by $\pi'$, which is  $\dim H^2(\tilde{S}, \C) - \dim H^2(S', \C)$. Thus, the above sequence is also short exact. 

Since all of the functors in \eqref{PartKlein} are prorepresentable and unobstructed, the fact that $\prod r_{\tilde{T}_i}$ and $\prod r_{T'_i}$ induce surjections on first order deformation implies that the functors are surjective when evaluated on any Artinian ring.
\end{proof}

In the next section we are going to need an analogue of diagram \eqref{LocDiagram} and \cite[Proposition 3.1]{NaPD1} for partial resolutions of Kleinian singularities, and we now have all the tools to prove it below.

\begin{prop} \label[prop]{PKLocSys}

Let $\cS' \rightarrow \A_{S'}$ be the universal graded deformation of $S'$. Then it is also a universal deformation of the germ $T' := (S', \rho^{-1}(0))$. Additionally, $\cS' \times \C^{2n - 2} \rightarrow \A_{S'}$ is a universal deformation of the germ $V' := T' \times (\C^{2n - 2}, 0)$.

\end{prop}

\begin{proof}

First we note that we can apply the arguments in \Cref{ProDef} in the complex analytic topology to show that $\pd_{T'}$ and $\pd_{V'}$ are prorepresentable. 

We first prove the proposition for $\pd_{T'}$. Clearly $\cS'$ is a deformation of the germ $T'$, and the induced map 
\begin{equation} \label{PartKleinUniv}
T_0(\A_{S'}) \cong \pd_{S'}(\C[\epsilon]) \rightarrow \pd_{T'}(\C[\epsilon])
\end{equation}
is simply the restriction map on first order deformations. The kernel of this map is by definition $\pd_{S'}^{\rho}(\C[\epsilon])$, which by \Cref{PDPartial} is isomorphic to $H^2(S, \C) = 0$, so \eqref{PartKleinUniv} is injective. Thus, we only need to check that $\dim \pd_{T'}(\C[\epsilon]) = \dim \pd_{S'}(\C[\epsilon])$. Letting $T'_i$ be the germ of each singular point in $S'$ as above, we have a short exact sequence
\begin{equation} \label{PartKleinLoc}
0 \rightarrow \pd_{lt, T'}(\C[\epsilon]) \rightarrow \pd_{T'}(\C[\epsilon]) \rightarrow \prod_{x_i \in \Sigma} \pd_{T'_i}(\C[\epsilon]) \cong \prod_{x_i \in \Sigma}\pd_{(S_i, 0)}.
\end{equation}

We want to show that the restriction map from $S'$ to $T'$ induces an isomorphism 
\begin{equation} \label{S'toT'}
r_{lt, T'}: \pd_{lt, S'}(\C[\epsilon]) \rightarrow \pd_{lt, T'}(\C[\epsilon]),
\end{equation}
which is equivalent to showing that every locally trivial deformation of $T'$ extends to $S'$, since we already know the map is injective. To do this let $\tilde{T}$ be the complex analytic germ $(\tilde{S}, \pi^{-1}(0))$, and observe that we have a pushforward map 
\begin{equation*}
(\pi'|_{\tilde{T}})_*: \pd_{\tilde{T}}(\C[\epsilon]) \rightarrow \pd_{T'}(\C[\epsilon]).
\end{equation*}

Given an element $\cT' \in \pd_{lt, T'}(\C[\epsilon])$ we can construct a lift $\tilde{\cT}$ to $\pd_{\tilde{T}}(\C[\epsilon])$ by hand as follows. Since $\cT'$ is a complex analytic deformation of $T'$, it is a deformation $\cU$ of an open neighborhood $U$ of $\rho^{-1}(0)$ in the complex analytic topology. In particular, $\cT'$ is determined by the data deformations $\cU_1, \dots , \cU_n$ of disjoint open neighborhoods of each singular point $x_i \in U$, a deformation $\cU_0$ of the smooth locus $U_0$ of $U$, and identifications $\Phi_i$ between $\cU_0$ and  $\cU_i$ on $U_0 \cap U_i$, for $i = 1, \dots, n$. By picking $U_i$ sufficiently small for $i = 1, \dots, n$, we can ensure $\cU_i$ is trivial for $i = 1, \dots , n$, since $\cT'$ is locally trivial. Since $\pi'$ is an isomorphism over $U_0$ we can treat $U_0$ as a subset of $\tilde{S}$. Additionally, we can use the $\Phi_i$ to glue $\cU_0$ to trivial deformations of each $\pi'^{-1}(U_i)$. Thus, we obtain a complex analytic deformation $\tilde{\cU}$ of $\pi'^{-1}(U)$ that pushes forward to $\cU$. In particular, treating $\tilde{\cU}$ as a complex analytic deformation of the germ $\tilde{T}$, it provides a lift of $\cT'$. 

Furthermore, by diagram \eqref{LocDiagram}, $\tilde{\cT}$ extends to an element $\tilde{\cS}$ of $\pd_{\tilde{S}}(\C[\epsilon])$ and thus $\pi'_*(\tilde{\cS})$ is an extension of $\cT'$ by construction. Thus $r_{lt, T'}$ is an isomorphism and 
\begin{equation} \label{ltTiH2S'}
\pd_{lt, T'}(\C[\epsilon]) \cong \pd_{lt, S'}(\C[\epsilon]) \cong H^2(S', \C),
\end{equation}
with the second isomorphism follows from \Cref{ltpd}. 

Returning to diagram \eqref{PartKleinLoc}, we see that
\begin{equation*}
\dim \pd_{T'}(\C[\epsilon]) \leq \dim H^2(S', \C) + \sum_{x_i \in \Sigma} \dim \pd_{(S_i, 0)} = \dim \pd_{S'}(\C[\epsilon])
\end{equation*}
so indeed the map \eqref{PartKleinUniv} is an isomorphism and $\cS'$ is a universal deformation of $T'$.

We now move on to consider the deformation $\cS' \times \C^{2n - 2} \rightarrow \A_{S'}$ of $V'$. Observe that $\cS' \times \C^{2n - 2}$ is a deformation of $V'$ and the induced map 
\begin{equation} \label{PartVUniv}
T_0(\A_{S'}) \cong \pd_{S'}(\C[\epsilon]) \rightarrow \pd_{V'}(\C[\epsilon])
\end{equation} 
is given by $\cS' \mapsto r_{T'}(\cS') \times (\C^{2n - 2}, 0)$. This map is clearly injective by the previous part, so we just need to check that $\dim \pd_{V'}(\C[\epsilon]) = \dim \pd_{S'}(\C[\epsilon])$. If we let $V'_i$ be the germ $(S' \times \C^{2n - 2}, (x_i, 0)) \cong (S_i, 0) \times (\C^{2n - 2}, 0)$, then we have the short exact sequence
\begin{equation} \label{PartVLoc}
0 \rightarrow \pd_{lt, V'}(\C[\epsilon]) \rightarrow \pd_{V'}(\C[\epsilon]) \rightarrow \prod_{x_i \in \Sigma} \pd_{V'_i}(\C[\epsilon]) \cong \prod_{x_i \in \Sigma} \pd_{(S_i, 0)},
\end{equation}
where the last isomorphism follows from \cite[Proposition 3.1]{NaPD1}. Finally, we want to show that the map $\pd_{lt, T'}(\C[\epsilon]) \rightarrow \pd_{lt, V'}(\C[\epsilon])$ given by $\cT' \mapsto \cT' \times (\C^{2n - 2}, 0)$ is an isomorphism. It is clearly injective by construction. To show that it is surjective note that as before we can construct a lift of any element $\cV' \in \pd_{lt, V'}(\C[\epsilon])$ to $\pd_{\tilde{V}}(\C[\epsilon])$, where $\tilde{V}$ is $\tilde{T} \times (\C^{2n - 2}, 0)$. However by \cite[Proposition 3.1]{NaPD1} every element of $\pd_{\tilde{V}}(\C[\epsilon])$ is of the form $\tilde{\cT} \times (\C^{2n - 2}, 0)$ for some $\tilde{\cT} \in \pd_{\tilde{T}}(\C[\epsilon])$. Thus, $\cV' = (\pi'|_{\tilde{T}})_*(\tilde{\cT}) \times (\C^{2n - 2}, 0)$ and indeed $\pd_{lt, T'}(\C[\epsilon]) \cong \pd_{lt, V'}(\C[\epsilon])$. Using equation \eqref{ltTiH2S'} we obtain $\pd_{lt, T'}(\C[\epsilon]) \cong H^2(S', \C)$, and thus by equation \eqref{PartVLoc} we have 
\begin{equation*}
\dim \pd_{V'}(\C[\epsilon]) \leq \dim H^2(S', \C) + \sum_{x_i \in \Sigma} \dim \pd_{(S_i, 0)} = \dim \pd_{S'}(\C[\epsilon]).
\end{equation*}
Thus, \eqref{PartVUniv} is an isomorphism and we are done.
\end{proof}

\subsection{Identifying $W_{S'}$}

In this subsection we will prove the main theorem, \Cref{MainThm}, for Kleinian singularities. In particular, we will show that the map $f': \A_{\tilde{S}} \to \A_{S'}$ from \Cref{UnivConic} is of the form $\mathfrak{h} \to \mathfrak{h}/W_{S'}$, where $W_{S'}$ is the Namikawa Weyl group of a given partial resolution $S'$ of a Kleinian singularity $S$. In order to do this, we first algebraicize the right square of \eqref{PartKlein}, to obtain a diagram of universal deformations. Next we show that this diagram is a pullback square. We have seen that every map $(\pi'|_{\tilde{T}_i})_*: \pd_{\tilde{T}_i} \to \pd_{T'_i}$ is a Galois cover with Galois group $W_i$, the Weyl group of the Kleinian singularity $T'_i$. This implies that $\pi'_*: \pd_{\tilde{S}} \to \pd_{S'}$ is a Galois cover whose Galois group is the product of the $W_i$. Finally, with the help of \Cref{subsetInv}, we identify this group as a subgroup of $W$, the Weyl group of $S$, with the Namikawa Weyl group $W_{S'}$ of $S'$. 

We are going to use some techniques from \cite[Theorem 1.1]{NaPD2} to algebraicize and then further study the right square of diagram \eqref{PartKlein}. Recall that for any algebraic variety (or germ of an algebraic variety) $Z$, such that $\pd_Z$ is prorepresentable and unobstructed, it has a universal formal deformation $\hat{\cZ} \rightarrow \hat{\A_Z}$ over the formal completion of the affine space
\begin{equation*}
\A_{Z} := \Spec \Sym(\pd_{Z}(\C[\epsilon])^*).
\end{equation*}
If we can algebraicize this deformation, then we obtain a deformation $\cZ$ over $\A_Z$ that is formally universal at zero. In this case we let $\PDef(Z)$ be a small neighborhood of $0 \in \A_Z$, so that $\cZ|_{\PDef(Z)}$ is also a formally universal deformation. Namikawa refers to this as a Kuranishi space for the Poisson deformations of $Z$.

Since all of the functors in diagram \eqref{PartKlein} are prorepresentable and unobstructed, we can take the associated formal spaces and maps. Doing so for the right square we obtain
\begin{equation} \label{FormalKlein}
\begin{tikzcd}[row sep=.42cm]
	{\widehat{\PDef(\tilde{S})}} && {\displaystyle\prod_{x_i \in \Sigma} \widehat{\PDef(\tilde{T}_i)}} \\
	\\
	{\widehat{\PDef(S')}} && {\displaystyle\prod_{x_i \in \Sigma} \widehat{\PDef(T'_i)}}
	\arrow["{\prod \hat{\phi}_i}", from=3-1, to=3-3]
	\arrow["{\hat{f}}"', from=1-1, to=3-1]
	\arrow["{\prod \hat{\varphi}_i}", from=1-1, to=1-3]
	\arrow["{\prod \hat{f}_i}"', from=1-3, to=3-3]
\end{tikzcd}
\end{equation}
By \Cref{UnivConic}, we can algebraicize the vertical maps in diagram \eqref{FormalKlein} above. In order to algebraicize each of the $\hat{\varphi}_i$ we use the period map construction that Namikawa describes in \cite[Theorem 1.1, Step 2]{NaPD2}. This construction works identically here and for each $x_i \in \Sigma$ it produces a map $\varphi_i: \PDef(\tilde{S}) \rightarrow \PDef(\tilde{T}_i)$ that algebraicizes $\hat{\varphi}_i$ and fits into the following diagram
\begin{equation} \label{varphi}
\begin{tikzcd}[row sep=.6cm]
	{\PDef(\tilde{S})} && {H^2(\tilde{S}, \C)} \\
	\\
	{\PDef(\tilde{T}_i)} && {H^2(\tilde{T}_i, \C) \cong H^2(\tilde{S}_i, \C)}
	\arrow["{p_{\tilde{S}}}", from=1-1, to=1-3]
	\arrow["{\varphi_i}"', from=1-1, to=3-1]
	\arrow["{r_{\tilde{S}_i}}"', from=1-3, to=3-3]
	\arrow["{p_{\tilde{T}_i}}", from=3-1, to=3-3]
\end{tikzcd}
\end{equation}
where $p_{\tilde{S}}$ and $p_{\tilde{T}_i}$ are open embeddings and $r_{\tilde{S}_i}$ is the pullback in cohomology. 

Thus, we have algebraicized every map in diagram \eqref{FormalKlein} except for the bottom map. Using Namikawa's argument in \cite[Theorem 1.1, Step 4(i)]{NaPD2}, these algebraizations imply that the bottom map algebraicizes as well. In particular, we obtain a diagram:
\begin{equation} \label{KleinAlg}
\begin{tikzcd}[row sep=scriptsize]
	{\PDef(\tilde{S})} && {\displaystyle\prod_{x_i \in \Sigma} \PDef(\tilde{T}_i)} \\
	\\
	{\PDef(S')} && {\displaystyle\prod_{x_i \in \Sigma} \PDef(T'_i)}
	\arrow["{\prod \phi'_i}", from=3-1, to=3-3]
	\arrow["f"', from=1-1, to=3-1]
	\arrow["{\prod \varphi_i}", from=1-1, to=1-3]
	\arrow["{\prod f_i}"', from=1-3, to=3-3]
\end{tikzcd}
\end{equation}
Note that each of the maps $f_i: \PDef(\tilde{T}_i) \rightarrow \PDef(T'_i)$ is the restriction to an open neighborhood of zero of the map $\mathfrak{h}_i \rightarrow \mathfrak{h}_i/W_i$ from diagram \eqref{fullKlein}, where $\mathfrak{h}_i$ and $W_i$ are the Cartan subalgebra and Weyl group associated to the Kleinian singularity $S_i$.

Following the strategy in \cite[Theorem 1.1, Step 4(iii)]{NaPD2} we obtain:

\begin{prop} \label[prop]{KleinPull}

Diagram \eqref{KleinAlg} is a pullback square.

\end{prop}

\begin{proof}

First we show that 
\begin{equation*}
\PDef(S') \times_{\prod \PDef(T'_i)} \prod \PDef(\tilde{T}_i)
\end{equation*}
is smooth. The map $\prod \phi'_i: \PDef(S') \rightarrow \prod \PDef(T'_i)$ is a map between smooth varieties and by \Cref{KleinSurj} the kernel of the tangent map
\begin{equation*}
T_0(\prod \phi'_i) = \prod r_{T'_i}: \pd_{S'}(\C[\epsilon]) \rightarrow \prod_{x_i \in \Sigma} \pd_{T'_i}(\C[\epsilon])
\end{equation*}
has dimension $\dim \pd_{lt, S'}(\C[\epsilon]) = \dim \PDef(S') - \dim \prod \PDef(T'_i)$. This implies that $\prod \phi'_i$ is a smooth map. Thus, by base change, the map 
\begin{equation*}
\PDef(S') \times_{\prod \PDef(T'_i)} \prod \PDef(\tilde{T}_i) \rightarrow \prod_{x_i \in \Sigma} \PDef(\tilde{T}_i)
\end{equation*}
is also smooth. In particular, since $\prod_{x_i \in \Sigma} \PDef(\tilde{T}_i)$ is smooth, $\PDef(S') \times_{\prod \PDef(T'_i)} \prod \PDef(\tilde{T}_i)$ is as well.

Next we show that the induced map 
\begin{equation*}
\PDef(\tilde{S}) \rightarrow \PDef(S') \times_{\prod \PDef(T'_i)} \prod \PDef(\tilde{T}_i)
\end{equation*}
is an isomorphism. Since both spaces are smooth it is enough to show the derivative 
\begin{equation} \label{InducedPull}
T_0(\PDef(\tilde{S})) \cong \pd_{\tilde{S}}(\C[\epsilon]) \rightarrow T_0(\PDef(S') \times_{\prod \PDef(T'_i)} \prod \PDef(\tilde{T}_i))
\end{equation}
is an isomorphism of vector spaces.

The map $\PDef(\tilde{T}_i) \rightarrow \PDef(T'_i)$ is the restriction to an open neighborhood of zero of the map $\mathfrak{h}_i \rightarrow \mathfrak{h}_i/W_i$ from diagram \eqref{fullKlein}. Thus the map $T_0(\PDef(\tilde{T}_i)) \rightarrow T_0(\PDef(T'_i))$ is the zero map between vector spaces of the same dimension. Thus, 
\begin{align*}
\dim T_0(\PDef(S') \times_{\prod \PDef(T'_i)} \prod \PDef(\tilde{T}_i)) &= \sum_{x_i \in \Sigma} \dim \pd_{\tilde{T}_i}(\C[\epsilon]) + \dim \pd_{lt, S'}(\C[\epsilon]) \\
&= \sum_{x_i \in \Sigma} \dim \pd_{\tilde{T}_i}(\C[\epsilon]) + \dim \pd_{\tilde{S}}^{\pi'}(\C[\epsilon]) \\
&= \dim \pd_{\tilde{S}}(\C[\epsilon]).
\end{align*}

Therefore, we only need to show that \eqref{InducedPull} is injective. Assume that we have some $\tilde{\cS}$ in the kernel of the map. Then because it is in the kernel of $\prod r_{\tilde{T}_i}(\C[\epsilon])$ it must be an element of $\pd_{\tilde{S}}^{\pi'}(\C[\epsilon])$. However, $\tilde{\cS}$ is also in the kernel of $\pi'_*(\C[\epsilon])$, and by \Cref{KerIsoPartial} this map is an isomorphism when restricted to $\pd_{\tilde{S}}^{\pi'}(\C[\epsilon])$. Thus $\tilde{\cS} = 0$, and the map \eqref{InducedPull} is injective.
\end{proof}

As described above, each of the maps $f_i: \PDef(\tilde{T}_i) \rightarrow \PDef(T'_i)$ are Galois covering with Galois group $W_i$. Thus, the map $\prod f_i$ in diagram \eqref{KleinAlg} is a Galois covering with Galois group $\prod_{x_i \in \Sigma} W_i$. Finally, by \Cref{KleinPull}, the map $f: \PDef(\tilde{S}) \rightarrow \PDef(S')$ is a Galois covering with the same Galois group. 

The map $f: \PDef(\tilde{S}) \rightarrow \PDef(S')$ can also be extended to a sequence 
\begin{equation*}
\PDef(\tilde{S}) \xrightarrow{f} \PDef(S') \xrightarrow{j} \PDef(S).
\end{equation*}
Since these are all graded deformations that are universal at zero, the above sequence is given by the restrictions to neighborhoods of zero in the sequence
\begin{equation*}
\A_{\tilde{S}} \cong \mathfrak{h} \xrightarrow{f'} \A_{S'} \xrightarrow{g} \A_{S} \cong \mathfrak{h}/W
\end{equation*}
that algebraicizes $\pd_{\tilde{S}} \rightarrow \pd_{S'} \rightarrow \pd_S$ to a diagram of universal graded deformations, from \Cref{UnivConic}.

In particular, because $\mathfrak{h} \rightarrow \mathfrak{h}/W$ is a Galois covering and $\A_{S'}$ is normal, the map $f'$ must be a quotient by some $W' \subset W$. However, we have already shown that $f: \PDef(\tilde{S}) \rightarrow \PDef(S')$ is a Galois covering with Galois group $\prod_{x_i \in \Sigma} W_i$, so we must have $W' \cong \prod_{x_i \in \Sigma} W_i$. 

All that is left to do is to identify $W'$ as a subgroup of $W$. 

\begin{prop}

The subgroup $W' \subset W$ is the Namikawa Weyl group $W_{S'} \subset W$ of the partial resolution $\rho: S' \rightarrow S$.

\end{prop}

\begin{proof}

We observe that for the symplectic resolution $\tilde{S} \rightarrow S$, the essential curve classes are simply the components of the exceptional fiber. In particular, $W' \cong \prod_{x_i \in \Sigma} W_i$ is isomorphic to the subgroup of $W$ generated by the simple generators corresponding to contracted essential curve classes. This subgroup is $W_{S'}$, as defined in \Cref{PartWeyl}. This implies that $W' \cong W_{S'}$ as groups, but we still need to show that they are the same subgroup of $W = W_{\tilde{S}}$. 

We observe that it follows from \Cref{KleinPull} that $\PDef(\tilde{S})^{W'} = \ker(\prod \varphi_i)$. We explain this below in more detail. Note that following \eqref{varphi} and \eqref{KleinAlg}, by $\ker(\prod \varphi_i)$ we mean the intersection of $\PDef(\tilde{S}) \subset \A_{\tilde{S}} \cong H^2(\tilde{S}, \C)$ with the kernel of the restriction map $H^2(\tilde{S}, \C) \to \prod H^2(\tilde{T}_i, \C)$ in cohomology. 

Since $f: \PDef(\tilde{S}) \rightarrow \PDef(S')$ is a quotient by $W'$, $x \in \PDef(\tilde{S})$ is an element of $\PDef(\tilde{S})^{W'}$ is and only if $f^{-1}(f(x))$ consists of a single point. By \Cref{KleinPull}, the fiber $f^{-1}(f(x))$ has the same size as the fiber 
$$
(\prod f_i)^{-1}(((\prod \phi'_i) \circ f)(x)).
$$
However, zero is the only point in $\prod \PDef(T'_i)$ whose fiber over $\prod f_i$ is a single point. Thus, $x \in \PDef(\tilde{S})^{W'}$ if and only if $((\prod \phi'_i) \circ f)(x) = 0$. By the commutativity of the \eqref{KleinAlg}, this is equivalent to $((\prod \varphi_i) \circ (\prod f_i))(x) = 0$, which is equivalent to $\prod \varphi_i(x) = 0$. In particular, $x \in \PDef(\tilde{S})^{W'}$ if and only if $\prod \varphi_i(x) = 0 \in \prod \PDef(\tilde{T}_i)$, which implies that 
\begin{equation} \label{W'inv}
\PDef(\tilde{S})^{W'} = \ker(\prod \varphi_i).
\end{equation}

On the other hand, we recall from diagram \eqref{varphi} that each $\varphi_i$ is given by restriction in cohomology, so $\ker(\prod \varphi_i)$ is the intersection of $\PDef(\tilde{S})$ with the subspace of $H^2(\tilde{S}, \C)$ generated by the dual fundamental classes $[C_i]^{\vee}$ of the exceptional curves that are not contracted by $\pi'$. Note that this subspace is $\pi'^*(H^2(S', \C))$. 

We recall that the period map $\mathfrak{h} \xrightarrow{\sim} H^2(\tilde{S}, \C)$ identifies these $[C_i]^{\vee}$ with the fundamental coweights in $\mathfrak{h}$ corresponding to the same nodes of the Dynkin diagram, and in particular by \Cref{subsetInv} we have
\begin{equation*}
\ker(\prod \varphi_i) = \PDef(\tilde{S})^{W_{S'}} \hspace{10pt} \text{and} \hspace{10pt} \Stab_W(\ker(\prod \varphi_i)) = W_{S'}.
\end{equation*}
Combining this with equation \eqref{W'inv}, we obtain that $W' \subset \Stab_W(\ker(\prod \varphi_i)) = W_{S'}$. However, since $W' \cong W_{S'}$, they must be the same subgroup of $W$.
\end{proof}

Thus, $f'$ is a quotient by the subgroup $W_{S'} \subset W$, and it is explicitly of the form 

\begin{equation*} 
\mathfrak{h} \cong H^2(S',  \C) \oplus \prod\limits_{x_i \in \Sigma} \mathfrak{h}_i \longrightarrow \mathfrak{h}/W_{S'} \cong H^2(S',  \C) \oplus \prod\limits_{x_i \in \Sigma} (\mathfrak{h}_i/W_i).
\end{equation*}

We can expand diagram \eqref{fullKlein} from the beginning of the section to the diagram below. In particular, this proves the main theorem of the paper, \Cref{MainThm}, for Kleinian singularities.
\begin{equation} \label{KleinResult}
\begin{tikzcd}[row sep=scriptsize]
	{\tilde{S}} && {\tilde{\cS}} && {\A_{\tilde{S}} \cong H^2(\tilde{S}, \C) \cong \mathfrak{h}} \\
	\\
	{S'} && {\cS'} && {\A_{S'} \cong \mathfrak{h}/W_{S'}} \\
	\\
	S && \cS && {\A_{S} \cong \mathfrak{h}/W}
	\arrow[hook, from=1-1, to=1-3]
	\arrow[from=1-1, to=3-1]
	\arrow[from=1-3, to=1-5]
	\arrow[from=1-3, to=3-3]
	\arrow[from=1-5, to=3-5]
	\arrow[hook, from=3-1, to=3-3]
	\arrow[from=3-1, to=5-1]
	\arrow[from=3-3, to=3-5]
	\arrow[two heads, from=3-3, to=5-3]
	\arrow[from=3-5, to=5-5]
	\arrow[hook, from=5-1, to=5-3]
	\arrow[from=5-3, to=5-5]
\end{tikzcd}
\end{equation}

The following lemma is not used anywhere else and follows largely from unpacking definitions, but we include it for completeness.

\begin{lemma}

Let $\rho:X' \rightarrow X$ be a partial resolution of a conical affine symplectic singularity, and let $\pi': Y \rightarrow X'$ be a $\Q$-factorial terminalization. Then the Namikawa Weyl group $W_{X'}$ of the partial resolution is the subgroup of $W_X$ that fixes the face $F \subset H^2(Y^{\reg}, \R) \cong \Pic(Y) \otimes \R$ corresponding to $\rho$.

\end{lemma}

\begin{proof}

Let $s$ be a simple generator of the Namikawa Weyl group $W_X$, which corresponds to a node $i$ of the folded Dynkin diagram $D_{\cL}$ corresponding to a codimension two leaf $\cL$ of $X$. Recall that the Weyl group $W_{\cL}$ of $D_{\cL}$ is a component of $W_X$. For any point $x \in \cL$, it is shown in \cite[Theorem 1.1]{NaPD2}, as mentioned in \cite[(1)]{BPW}, that the restriction map 
\begin{equation} \label{HyperRest}
\Pic(Y) \otimes \R  \cong H^2(Y^{\reg}, \R) \rightarrow H^2(\pi^{-1}(x), \R)^{\pi_1(\cL, x)}
\end{equation}
is $W_{\cL}$-equivariant. This is also explained in \Cref{Main}. Note that the action of $W_{\cL} = \hat{W}_{\cL}^{\pi_1(\cL, x)}$ on $H^2(\pi^{-1}(x), \R)^{\pi_1(\cL, x)}$ is given by restriction of the action of $\hat{W}_{\cL}$ on $H^2(\pi^{-1}(x), \R)$, where $\hat{W}_{\cL}$ is the Weyl group corresponding to the Kleinian singularity given by taking a transverse slice $S_{\cL}$ to $\cL$. Let $\hat{D}_{\cL}$ be the simply laced Dynkin diagram corresponding to $S_{\cL}$, and let $J$ be set of nodes of $\hat{D}_{\cL}$ in the $\pi_1(\cL, x)$-orbit corresponding to the node $i$ of $D_{\cL}$.

The hyperplane $H_s$ corresponding to $s$ is given by the set of classes that pair to zero with the essential curve class $[C] \in \cN_1(Y/X)$ corresponding to the same node $i$ of $D_{\cL}$ as $s$. In particular, $H_s$ pulls back under \eqref{HyperRest} to the intersection of $H^2(\pi^{-1}(x), \R)^{\pi_1(\cL_i, x)}$ with the subspace generated by the set of classes in $H^2(\pi^{-1}(x), \R) \cong \Pic(\pi^{-1}(x)) \otimes \R$ that pair to zero with every projective line in $\pi^{-1}(x)$ in the locus of $[C]$. This is the subspace of $H^2(\pi^{-1}(x), \R)$ generated by all the dual cohomology classes $[C_k]^{\vee}$ corresponding to nodes of $\hat{D}_{\cL}$ that are not in $J$. By \Cref{subsetInv}, this is $H^2(\pi^{-1}(x), \R)^{W_J}$, where $W_J$ is the parabolic subgroup of $\hat{W}_{\cL}$ corresponding to $J$. Thus, $H_s$ pulls back under \eqref{HyperRest} to
\begin{equation} \label{Hs}
H^2(\pi^{-1}(x), \R)^{\pi_1(\cL, x)} \cap H^2(\pi^{-1}(x), \R)^{W_J}.
\end{equation}

Furthermore, this subspace \eqref{Hs} is generated by the sums $\sum_{k \in i'} [C_k]^{\vee}$ for every $\pi_1(\cL, x)$-orbit $i'$ in $\hat{D}_{\cL}$ except the one corresponding to $i$. This subspace clearly has codimension one in $H^2(\pi^{-1}(x), \R)^{\pi_1(\cL, x)}$, and it is invariant under the reflection $s \in W_{\cL} \subset W_X$, since $s$ is a product of elements of $W_J$. Thus, this is the subspace fixed by $s$. Since \eqref{HyperRest} is $W_{\cL}$-equivariant, the hyperplane $H_s \subset \Pic(Y) \otimes \R$ is the hyperplane fixed by $s$.

Given a face $F$ corresponding to a partial resolution $\rho:X' \rightarrow X$, $F$ is open in an intersection $H_1 \cap \dots \cap H_k$ of hyperplanes. Out of these, let $H_{s_1}, \dots, H_{s_k}$ be the subset of these hyperplanes that correspond to simple generators of $W_X$. Then, by the above work, the stabilizer of $F$ is the subgroup of $W_X$ generated by $s_1, \dots, s_k$. This is precisely the Namikawa Weyl group $W_{X'}$ of $\rho$, as defined in \Cref{PartWeyl}.
\end{proof}

\section{Main Theorem} \label[section]{Main}

In this section we adapt the strategy from \cite[Theorem 1.1]{NaPD2} using results from the previous sections, particularly the previous section, to prove the main result. 

Let $Y \xrightarrow{\pi'} X' \xrightarrow{\rho} X$ be such that $X$ is an affine symplectic singularity, $X'$ is a partial resolution of $X$, and $Y$ is a $\Q$-factorial terminalization of $X$ and $X'$. Furthermore, let $\pi = \rho \circ \pi'$. Let $\{\cL_i\}_{i \in I}$ be the set of codimension two symplectic leaves of $X$ and let $U$ be the union of the open leaf and the codimension two leaves, i.e. the complement of the symplectic leaves of codimension $\geq$ 4.  Recall that by \cite[p. 62]{NaPD1}, we have $\pd_X(\C[\epsilon]) \cong \pd_U(\C[\epsilon])$. Let $\Sigma$ be the singular locus of $U$, i.e. the union of the $\cL_i$. For each $i \in I$, pick a point $x_i \in \cL_i$, and let $S_i$ be the Kleinian singularity given by taking a transverse slice to $\cL_i$. Note that $S_i$ does not depend on the choice of $x_i$.

In order to study the Poisson deformations of $X'$, we have to analyze their behavior over each of the codimension two strata of $X$. For this we will use some results of Namikawa that essentially tell us that the local behavior is determined by the corresponding $S_i$.

\begin{lemma}[{\cite[Lemma 1.3, Lemma 1.6]{NaPD1}}]  \label[lemma]{LocSys}

For each $x_i$, there is an analytic neighborhood $U_i$ of $x_i$ such that $U_i$ is isomorphic to an open subset of $S_i \times \C^{2n - 2}$ as a Poisson variety. 

Furthermore, every first order deformation $\cV \in \pd_{(X, x_i)}$ of the germ 
\begin{equation*}
(X, x_i) \cong (S_i, 0) \times (\C^{2n - 2}, 0)
\end{equation*}
is a product of deformations of $(S_i, 0)$ and $(\C^{2n - 2}, 0)$. In particular it is of the form 
\begin{equation*}
\cV \cong (\cS_i, 0) \times (\C^{2n - 2}, 0)
\end{equation*}
for some deformation $\cS_i \in \pd_{S_i}(\C[\epsilon])$.
\end{lemma}

For each codimension two leaf $\cL_i$, we set $F'_i := \rho^{-1}(x_i)$ and $\tilde{F}_i := \pi^{-1}(x_i)$. Additionally, we define the germs $T_i := (X, x_i)$, $T'_i := (X', F'_i)$, and $\tilde{T}_i := (Y, \tilde{F}_i)$. Then we have the following diagram, where the isomorphisms in the middle column are proved in \Cref{AllSmoothSec}:
\begin{equation} \label{Partial}
\begin{tikzcd} [row sep=.49cm]
	{\pd_{\tilde{U}}^{\pi}} && {\pd_{\tilde{U}} \cong \pd_Y} && {\displaystyle\prod_{i \in I} \pd_{\tilde{T}_i}} \\
	\\
	{\pd_{U'}^{\rho}} && {\pd_{U'} \cong \pd_{X'}} && {\displaystyle \prod_{i \in I} \pd_{T'_i}} \\
	\\
	{\pd_{lt, U}} && {\pd_{U} \cong \pd_X} && {\displaystyle\prod_{i \in I} \pd_{T_i}}
	\arrow["{\prod r_{T'_i}}", from=3-3, to=3-5]
	\arrow["{\rho_*}"', from=3-3, to=5-3]
	\arrow[from=3-5, to=5-5]
	\arrow["{\prod r_{T_i}}", from=5-3, to=5-5]
	\arrow[hook, from=5-1, to=5-3]
	\arrow[hook, from=3-1, to=3-3]
	\arrow[from=3-1, to=5-1]
	\arrow[from=1-1, to=3-1]
	\arrow[hook, from=1-1, to=1-3]
	\arrow["{\pi'_*}"', from=1-3, to=3-3]
	\arrow["{\prod r_{\tilde{T}_i}}", from=1-3, to=1-5]
	\arrow[from=1-5, to=3-5]
\end{tikzcd}
\end{equation}

Note that by \Cref{LocSys}, $T_i \cong (S_i, 0) \times (\C^{2n - 2}, 0)$ and this implies that 
\begin{equation} \label{tTi&T'i}
\tilde{T}_i \cong (\tilde{S}_i, \tilde{F}_i) \times (\C^{2n - 2}, 0) \hspace{10pt} \text{and} \hspace{10pt} T'_i \cong (S'_i, F'_i) \times (\C^{2n - 2}, 0),
\end{equation}
where $\tilde{S}_i$ is the crepant resolution of $S_i$ and $S'_i$ is the partial resolution given by the restriction of $X'$. Furthermore, the second part of \Cref{LocSys} and \Cref{PKLocSys} imply that the sequence $\pd_{\tilde{T}_i} \rightarrow \pd_{T'_i} \rightarrow \pd_{T_i}$ given by the natural pushforward maps is isomorphic to $\pd_{\tilde{S}_i} \rightarrow \pd_{S'_i} \rightarrow \pd_{S_i}$. This sequence has compatible universal graded deformations, as described in \eqref{KleinResult}. In particular, the sequence given by taking the bases of the universal deformations is of the form 
\begin{equation} \label{AlgRightColumn}
\hat{\mathfrak{h}_i} \rightarrow \hat{\mathfrak{h}_i}/W_{S'_i} \rightarrow \hat{\mathfrak{h}_i}/W_{S_i},
\end{equation}
where $\hat{\mathfrak{h}_i}$ is the Cartan subalgebra of the Lie algebra corresponding to $S_i$.

Next we study the rows of diagram \eqref{Partial}.

\begin{lemma} 
The sequences
\begin{gather*}
\pd_{\tilde{U}}^{\pi}(\C[\epsilon]) \rightarrow \pd_{\tilde{U}}(\C[\epsilon]) \rightarrow \prod_{i \in I} \pd_{\tilde{T}_i}(\C[\epsilon]) \\
\pd_{U'}^{\rho}(\C[\epsilon]) \rightarrow \pd_{U'}(\C[\epsilon]) \rightarrow \prod_{i \in I} \pd_{T'_i}(\C[\epsilon]) \\
\pd_{lt, U}(\C[\epsilon]) \rightarrow \pd_{U}(\C[\epsilon]) \rightarrow \prod_{i \in I} \pd_{T_i}(\C[\epsilon]) 
\end{gather*}
induced by the rows of diagram \eqref{Partial} are each exact at the middle term. 
\end{lemma}

\begin{proof}

We shall prove the result for the middle row in the statement of the lemma. The other two cases are boundary cases of this one, and can be deduced directly by the same method. 

By definition, $\cU' \in \pd_{U'}^{\rho}(\C[\epsilon])$ restricts to a Poisson deformation $\cU'|_{T'_i}$ of each germ $T'_i$ whose underlying flat deformation is trivial. However, by \Cref{PKLocSys}, this implies that $\cU'|_{T'_i}$ is also trivial as a Poisson deformation. Thus, $\pd_{U'}^{\rho}(\C[\epsilon])$ is contained in the kernel of
$$
\prod  r_{T'_i}: \pd_{U'}(\C[\epsilon]) \to \prod_{i \in I} \pd_{T'_i}(\C[\epsilon]).
$$

In order to prove the converse, let $\cU' \in \ker(\prod  r_{T'_i})$. We want to show that $\cU' \in \pd_{U'}^{\rho}(\C[\epsilon])$. In particular, we want to show that for any $x \in U$, there is some analytic neighborhood $V_x$ of $x$ such that $\cU'$ restricts to the trivial deformation of $V$. This is immediate if $x$ is in the smooth locus of $U$, so we can pick $x$ to be in the singular locus of $U$, which means it is in some codimension two leaf $\cL_i$. Let $C$ be a closed path from $x$ to $x_i$ in $\cL_i$. Given any point $y \in C$, by \eqref{tTi&T'i} and \Cref{PKLocSys}, we can find a neighborhood $V_y$ of $y$ such that there is an open embedding $i_y: V_y \hookrightarrow S'_i \times \C^{2n - 2}$ and $\cU'|_{V_y} \cong (\cS'_y \times \C^{2n - 2})|_{V_y}$ for some $\cS'_y \in \pd_{S'}(\C[\epsilon])$. Since $C$ is compact, we can take a finite subcover $V_{x_i} = V_{y_0}, V_{y_1}, \dots, V_{y_k} = V_{x}$. Since $\cU' \in \pd_{U'}^{\rho}(\C[\epsilon])$, $\cU'|_{V_{y_0}}$ must be the trivial deformation. Thus, we can pick some $y' \in V_{y_0} \cap V_{y_1}$ and it follows that $\cU'$ restricts to the trivial deformation of the germ $(U', \rho^{-1}(y'))$. By the construction of the $V_y$, this implies that $\cU'|_{V_{y_1}}$ is trivial as well. Applying this argument inductively we deduce that each $\cU'|_{V_{y_j}}$ is trivial. In particular, $\cU'|_{V_{x}}$ is trivial, so we are done. 
\end{proof}

\begin{thm}[Main Theorem] \label[thm]{MainThm}

We have the following diagram of universal graded deformations
\begin{equation*}
\begin{tikzcd}
	Y && \cY & {\A_{Y}} \\
	\\
	{X'} && {\cX'} & {\A_{X'} \cong \A_{Y}/W_{X'}} \\
	\\
	X && \cX & {\A_{X} \cong \A_{Y}/W_{X}}
	\arrow["{i_{Y}}", hook, from=1-1, to=1-3]
	\arrow["{\pi'}"', from=1-1, to=3-1]
	\arrow[from=1-3, to=1-4]
	\arrow["{\tilde{\pi}'}"', from=1-3, to=3-3]
	\arrow["{f'}"', from=1-4, to=3-4]
	\arrow["{i_{X'}}", hook, from=3-1, to=3-3]
	\arrow["\rho"', from=3-1, to=5-1]
	\arrow[from=3-3, to=3-4]
	\arrow["{\tilde{\rho}}"', from=3-3, to=5-3]
	\arrow["g"', from=3-4, to=5-4]
	\arrow["{i_X}", hook, from=5-1, to=5-3]
	\arrow[from=5-3, to=5-4]
\end{tikzcd}
\end{equation*}
with $f'$ and $f := g \circ f'$ quotients by $W_{X'}$ and $W_X$ respectively. 

\end{thm}

\begin{proof}

We divide the proof into three steps. 

\textbf{1} (\textit{Algebraization}): Since all of the functors in the right two columns of diagram \eqref{Partial} are prorepresentable and unobstructed, we can take the associated formal spaces and maps. Thus we obtain:
\begin{equation} \label{FormaltoAlg}
\begin{tikzcd}[row sep=scriptsize]
	{\widehat{\PDef(\tilde{U})} \cong \widehat{\PDef(Y)}} && {\displaystyle\prod_{i \in I} \widehat{\PDef(\tilde{T}_i)}} \\
	\\
	{\widehat{\PDef(U')} \cong \widehat{\PDef(X')}} && {\displaystyle\prod_{i \in I} \widehat{\PDef(T'_i)}} \\
	\\
	{\widehat{\PDef(U)} \cong \widehat{\PDef(X)}} && {\displaystyle\prod_{i \in I} \widehat{\PDef(T_i)}}
	\arrow["{\prod  \hat{\varphi}_i}", from=1-1, to=1-3]
	\arrow["{\hat{f'}}"', from=1-1, to=3-1]
	\arrow["{\prod \hat{f}'_i}"', from=1-3, to=3-3]
	\arrow["{\prod \hat{\phi}'_i}", from=3-1, to=3-3]
	\arrow["{\hat{g}}"', from=3-1, to=5-1]
	\arrow["{\prod \hat{g}_i}"', from=3-3, to=5-3]
	\arrow["{\prod \hat{\phi}_i}", from=5-1, to=5-3]
\end{tikzcd}
\end{equation}

We apply \Cref{UnivConic} directly to algebraicize the left column. It is explained after diagram \eqref{Partial} how we can also use \Cref{UnivConic} to algebraicize the right column as well, and we obtain open neighborhoods of zero in diagram \eqref{AlgRightColumn}.

As before we apply the period map construction from \cite[Theorem 1.1, Step 2]{NaPD2} to algebraicize each formal map
$$
\hat{\varphi}_i: \widehat{\PDef(\tilde{U})} \rightarrow \widehat{\PDef(\tilde{T}_i)}
$$
induced by $r_{\tilde{T}_i}: \pd_{\tilde{U}} \rightarrow \pd_{\tilde{T}_i}$ to a map $\varphi_i: \PDef(\tilde{U}) \rightarrow \PDef(\tilde{T}_i)$, which fits into the diagram

\begin{equation} \label{PeriodAlg}
\begin{tikzcd}
	{\PDef(\tilde{U}) \cong \PDef(Y)} && {H^2(\tilde{U}, \C) \cong H^2(Y^{\reg}, \C)} \\
	\\
	{\PDef(\tilde{T}_i)} && {H^2(\tilde{T}_i, \C) \cong H^2(\tilde{S}_i, \C)}
	\arrow["{p_{\tilde{U}}}", from=1-1, to=1-3]
	\arrow["{\varphi_i}"', from=1-1, to=3-1]
	\arrow["{r_{\tilde{S}_i}}"', from=1-3, to=3-3]
	\arrow["{p_{\tilde{T}_i}}", from=3-1, to=3-3]
\end{tikzcd}
\end{equation}
where $p_{\tilde{U}}$ and $p_{\tilde{T}_i}$ are open embeddings and $r_{\tilde{S}_i}$ is the pullback in cohomology. Finally, we can again apply Namikawa's argument in \cite[Theorem 1.1, Step 4(i)]{NaPD2} to algebraicize the bottom and middle maps. Thus we obtain a diagram:

\begin{equation} \label{MainAlg}
\begin{tikzcd} [row sep=scriptsize]
	{\PDef(Y)} && {\displaystyle\prod_{i \in I} \PDef(\tilde{T}_i)} \\
	\\
	{\PDef(X')} && {\displaystyle\prod_{i \in I} \PDef(T'_i)} \\
	\\
	{\PDef(X)} && {\displaystyle\prod_{i \in I} \PDef(T_i)}
	\arrow["{\prod \phi'_i}", from=3-1, to=3-3]
	\arrow["g"', from=3-1, to=5-1]
	\arrow["{\prod g_i}"', from=3-3, to=5-3]
	\arrow["{\prod \phi_i}", from=5-1, to=5-3]
	\arrow["{f'}"', from=1-1, to=3-1]
	\arrow["{\prod \varphi_i}", from=1-1, to=1-3]
	\arrow["{\prod f'_i}"', from=1-3, to=3-3]
\end{tikzcd}
\end{equation}

Note that we are slightly abusing notation by denoting the maps $f'$, $g$, and $f$ on the bases of the universal graded deformations the same way as their restrictions to open neighborhoods of zero above. 

\textbf{2} (\textit{Image of restriction to germs}): We want to study the image of the left column of the above diagram in the right column. We recall from diagram \eqref{AlgRightColumn} that the right column is given by restricting the maps 
$$
\hat{\mathfrak{h}_i} \rightarrow \hat{\mathfrak{h}_i}/W_{S'_i} \rightarrow \hat{\mathfrak{h}_i}/W_{S_i}
$$
to open neighborhoods of zero. 

Note that the image of $r_{\tilde{S_i}}$ in diagram \eqref{PeriodAlg} and thus the image $V_i$ of $\varphi_i$ is the invariants under the monodromy action of $\pi_1(\cL_i, x_i)$. Let $m$ be the sum of the dimensions of the $V_i$, i.e. the dimension of the image of $\prod \varphi_i$. The $V_i$ are each smooth, and each map $\varphi_i$ is smooth over $V_i$. Thus, $\dim \pd_{\tilde{U}}^{\pi}(\C[\epsilon]) = \dim \PDef(Y) - m$.

Next we study the maps $\prod \phi'_i$ and $\prod \phi_i$. Since the $f'_i$ are finite and surjective we have $m = \dim \prod f'_i(V_i)$. Additionally, since $f'$ is surjective, $\prod f'_i(V_i)$ is also the image of $\prod \phi'_i$. By the short exact sequence 
\begin{equation} \label{sesVsm}
0 \rightarrow \pd_{U'}^{\rho}(\C[\epsilon]) \rightarrow \pd_{X'}(\C[\epsilon]) \rightarrow \prod \pd_{T'_i}(\C[\epsilon]),
\end{equation}
the kernel of the tangent map of $\prod \phi'_i$ at 0 is $\pd_{U'}^{\rho}(\C[\epsilon])$. In particular,
\begin{equation*}
\dim \pd_{U'}^{\rho}(\C[\epsilon]) = \dim \pd_{\tilde{U}}^{\pi}(\C[\epsilon]) = \dim \PDef(Y) - m = \dim \PDef(X') - \dim \prod f'_i(V_i),
\end{equation*}
so $\prod \phi'_i$ is smooth as a map onto its image. Since $\PDef(X')$ is also smooth this implies that $\prod f'_i(V_i)$ is smooth as well. The same argument applies to $\prod \phi_i$.
 
Let $\hat{W}_i$ be the Weyl group of the Lie algebra $\hat{\mathfrak{g}_i}$ corresponding to $S_i$. Then by \cite[Lemma 1.2]{NaPD2}, the subgroup of $\hat{W_i}$ which preserves $V_i$ is the Weyl group $W_i$ of the folded Dynkin diagram $\hat{D}_i/\pi_1(\cL_i, x_i)$. Thus, by \Cref{FoldingResult}, the subgroup of $W_{S'_i}$ which preserves $V_i$ is $W'_i := W_i \cap W_{S'_i}$, which is the parabolic subgroup of $W_i$ generated by the simple reflections corresponding to curve classes contracted by $\pi'$. In particular, the restriction of $f'_i$ to $V_i$ factors as 
\begin{equation*}
V_i \rightarrow V_i/W'_i \rightarrow f'_i(V_i) = \phi'_i(\PDef(X')).
\end{equation*}
We note that the second map above must be birational. This is because for any $w \in W_{S'_i}$ that doesn't preserve $V_i$, the image $w(V_i) \cap V_i$ is at most a hyperplane in $V_i$ and thus the quotient by the action of $w$ is birational. Since  $W_{S'_i}$ is a finite group this implies the second map is birational. However we showed above that $\phi'_i(\PDef(X'))$ is smooth, and $V_i/W'_i$ is also smooth since it is an open neighborhood of a quotient of $\mathfrak{h}_i$ by a finite reflection group. Thus, the induced map $V_i/W'_i \rightarrow f'_i(V_i) = \phi'_i(\PDef(X'))$ is a birational finite map between smooth varieties and thus an isomorphism. 

Clearly the same argument shows that the maps $f_i|_{V_i}: V_i \rightarrow \phi_i(\PDef(X))$ are also given by the restriction to open neighborhoods of zero of the natural quotient $\mathfrak{h}_i \rightarrow \mathfrak{h}_i/W_i$. By \Cref{Stein}, the $W_i$ act generically freely on the $\mathfrak{h}_i$, which implies the quotients by $W_i$ and $W'_i$ above are Galois covers with Galois groups $W_i$ and $W'_i$ respectively.

Finally, we observe that the product $\prod_{i \in I} W_i$ is the product of Weyl groups of the folded Dynkin diagrams $\hat{D}_i/\pi_1(\cL_i, x_i)$, which is the Namikawa Weyl group $W_X$. Furthermore, $\prod_{i \in I} W'_i$ is the parabolic subgroup generated by simple generators which correspond to essential curve classes contracted by $\pi'$, which is precisely the Namikawa Weyl group $W_{X'}$, as defined in \Cref{PartWeyl}. Thus, restricting the right column of diagram \eqref{MainAlg} to the image of the left column we obtain the following diagram, where all the maps are surjective. 

\begin{equation} \label{PullSquares}
\begin{tikzcd}
	{\PDef(Y)} && {\displaystyle\prod_{i \in I} V_i} \\
	\\
	{\PDef(X')} && {\displaystyle\prod_{i \in I} f'_i(V_i)\cong \prod_{i \in I} V_i/W'_i \cong (\prod_{i \in I} V_i)/W_{X'} } \\
	\\
	{\PDef(X)} && {\displaystyle\prod_{i \in I} f_i(V_i)\cong \prod_{i \in I} V_i/W_i \cong (\prod_{i \in I} V_i)/W_{X}}
	\arrow["{\prod \phi'_i}", from=3-1, to=3-3]
	\arrow["g"', from=3-1, to=5-1]
	\arrow["{\prod g_i}"', from=3-3, to=5-3]
	\arrow["{\prod \phi_i}", from=5-1, to=5-3]
	\arrow["{f'}"', from=1-1, to=3-1]
	\arrow["{\prod \varphi_i}", from=1-1, to=1-3]
	\arrow["{\prod f'_i}"', from=1-3, to=3-3]
\end{tikzcd}
\end{equation}

\textbf{3} (\textit{Pullback squares}): Next we are going to show that the squares in diagram \eqref{PullSquares} are pullback squares. We are going to show this for the bottom square, as the proof for the top square and the big square are identical. 

We begin by showing that 
\begin{equation*}
\PDef(X) \times_{\prod f_i(V_i)} \prod f'_i(V_i)
\end{equation*}
is smooth. As explained above, $\prod \phi_i$ is smooth, which by base change further implies that the map 
\begin{equation*}
\PDef(X) \times_{\prod f_i(V_i)} \prod f'_i(V_i) \rightarrow \prod f'_i(V_i)
\end{equation*}
is smooth. Since $\prod f'_i(V_i)$ is smooth, we have shown that $\PDef(X) \times_{\prod f_i(V_i)} \prod f'_i(V_i)$ is smooth.

Then we show that the induced map 
\begin{equation*}
\PDef(X') \rightarrow \PDef(X) \times_{\prod f_i(V_i)} \prod f'_i(V_i)
\end{equation*}
is an isomorphism. Since both spaces are smooth it is enough to show the derivative 
\begin{equation} \label{FullIndPull}
T_0(\PDef(X')) \cong \pd_{X'}(\C[\epsilon]) \rightarrow T_0(\PDef(X) \times_{\prod f_i(V_i)} \prod f'_i(V_i))
\end{equation}
is an isomorphism of vector spaces.

Because the linear map $T_0(\PDef(X)) \rightarrow T_0(\prod f_i(V_i))$ induced by $\prod \phi_i$ is surjective, we have 
\begin{gather*}
\dim T_0(\PDef(X) \times_{\prod f_i(V_i)} \prod f'_i(V_i)) \\
= \dim T_0(\PDef(X)) - \dim T_0(\prod f_i(V_i)) + \dim T_0(\prod f'_i(V_i)) = \dim T_0(\PDef(X))
\end{gather*}

Therefore, we only need to show that the map \eqref{FullIndPull} is injective. Assume that we have some $\cX'$ in the kernel of the map. Then it is in the kernel of $\prod r_{T'_i}(\C[\epsilon])$ from diagram \eqref{Partial}, so it must be an element of $\pd_{X'}^{\rho}(\C[\epsilon])$. However, $\cX'$ is also in the kernel of $\rho_*(\C[\epsilon])$, and by \Cref{KerIsoPartial} this map is an isomorphism when restricted to $\pd_{X'}^{\rho}(\C[\epsilon])$. Thus $\cX' = 0$, and the map \eqref{FullIndPull} is injective. This completes the proof that the squares in diagram \eqref{PullSquares} are pullback squares. 

Thus, we can apply the results about the right column of diagram \eqref{PullSquares} to deduce that the left column is of the form 
\begin{equation} \label{ResultNbhd}
\PDef(Y) \rightarrow \PDef(X') \cong \PDef(Y)/W_{X'} \rightarrow \PDef(X) \cong \PDef(Y)/W_X,
\end{equation}
where the first map and the composition are Galois covers with Galois groups $W_{X'}$ and $W_X$ respectively.

Finally, we recall from \Cref{UnivConic} that we have compatible universal graded deformations of $X$, $X'$, and $Y$. We have shown that if we take the induced maps
\begin{equation*}
\A_Y \xrightarrow{f'} \A_{X'} \xrightarrow{g} \A_{X}
\end{equation*}
from \Cref{UnivConic} and restrict them to neighborhoods of zero we obtain diagram \eqref{ResultNbhd}. Since the universal deformations are conical, this completes the proof.
\end{proof}

\section{Springer Theory} \label[section]{Springer}

Let $\rho: X' \rightarrow X$ be a partial resolution of a conical affine symplectic singularity $X$. For this section, we will require that $\pi: Y \rightarrow X$ is a symplectic resolution (not simply a $\Q$-factorial terminalization) that covers $\rho$, with induced map $\pi': Y \rightarrow X'$. Note that it follows from \cite[Corollary 5.6]{NaPD1} that since $X$ has a symplectic resolution, every $\Q$-terminalization of $X$ is a symplectic resolution.  

We reproduce the diagram from \Cref{MainThm} below but with the bases of the universal graded deformations of $Y$, $X'$, and $X$ denoted $B_Y$, $B_{X'}$, and $B_X$, following \cite{MN}.
\[\begin{tikzcd}
	Y && \cY & {B_Y} \\
	\\
	{X'} && {\cX'} & {B_{X'} \cong B_Y/W_{X'}} \\
	\\
	X && \cX & {B_{X} \cong B_Y/W_{X}}
	\arrow["{i_{Y}}", hook, from=1-1, to=1-3]
	\arrow[from=1-3, to=1-4]
	\arrow["{\pi'}"', from=1-1, to=3-1]
	\arrow["{\tilde{\pi}'}"', from=1-3, to=3-3]
	\arrow["{i_{X'}}", hook, from=3-1, to=3-3]
	\arrow[from=1-4, to=3-4]
	\arrow[from=3-3, to=3-4]
	\arrow["\rho"', from=3-1, to=5-1]
	\arrow["{\tilde{\rho}}"', from=3-3, to=5-3]
	\arrow["{i_X}", hook, from=5-1, to=5-3]
	\arrow[from=5-3, to=5-4]
	\arrow[from=3-4, to=5-4]
\end{tikzcd}\]
We further let $\tilde{\pi} := \tilde{\pi}' \circ \tilde{\rho}$ and $\pi := \pi' \circ \rho$. Following \cite{MN}, we define the symplectic Harish-Chandra sheaves $\HC' \in \Perv(\cX')$ and $\HC \in \Perv(\cX)$ by 
\begin{equation*}
\HC' := \tilde{\pi}'_*{\C_{\cY}}[\dim(\cY)], \hspace{20pt} \HC := \tilde{\pi}_*{\C_{\cY}}[\dim(\cY)] = \tilde{\rho}_*{\HC'}.
\end{equation*}
Similarly, we define the symplectic Springer sheaves $\Spr' \in \Perv(X')$ and $\Spr \in \Perv(X)$ by 
\begin{equation*}
\Spr' := \pi'_*{\C_{Y}}[\dim(Y)], \hspace{20pt} \Spr := \pi_*{\C_{Y}}[\dim(Y)] = \rho_*{\Spr'}.
\end{equation*}
Finally, we define $\Perv_{\sympl}(X')$ (resp. $\Perv_{\sympl}(X))$ to be the category of perverse sheaves on $X'$ (resp. $X$) that are constructible with respect to the stratification by symplectic leaves. We recall the definition of a small map below. 

\begin{df}[{\cite[Definition 3.8.5]{PA}}]

Let $f: Y \to X$ be a morphism of varieties with $X$ irreducible and $U \subset X$ an open and dense subvariety. Then $f$ is small with respect to $U$ if the following hold:

\begin{enumerate}

\item For each $x \in U$, $f^{-1}(x)$ is a finite set.
\item There exists a stratification of $X$ such that $U$ is a union of strata, and for every stratum $S \subset X$ and point $x \in S$, we have 
$$
\dim f^{-1}(x) < \frac{1}{2}(\dim Y - \dim S).
$$

\end{enumerate}

\end{df}

Let $B_X^{\reg}$ be the complement of the discriminant locus and $\cX^{\reg} := \cX \times_{B_X} B_X^{\reg}$. Then, the reasoning of \cite[Proposition 2.15]{MN} implies the following lemma:

\begin{lemma} \label[lemma]{small}

The maps $\tilde{\rho}$ and $\tilde{\pi}'$ are both small with respect to the dense open sets $\cX^{\reg}$ and $\tilde{\rho}^{-1}(\cX^{\reg})$ respectively. 

\end{lemma}

Furthermore, the restriction $\tilde{\pi}'|_{\tilde{\pi}^{-1}(\cX^{\reg})}: \tilde{\pi}^{-1}(\cX^{\reg}) \rightarrow \tilde{\rho}^{-1}(\cX^{\reg})$ is an unramified Galois $W_{X'}$-cover and the map $\tilde{\pi}|_{\tilde{\pi}^{-1}:(\cX^{\reg})}: \tilde{\pi}^{-1}(\cX^{\reg}) \rightarrow \cX^{\reg}$ is an unramified Galois $W_X$ cover, such that the composition
\begin{equation} \label{Gal}
\tilde{\pi}'^{-1}(\cX^{\reg}) \rightarrow \tilde{\rho}^{-1}(\cX^{\reg}) \rightarrow \cX^{\reg}
\end{equation}
is given by 
\begin{equation*}
\tilde{\pi}^{-1}(\cX^{\reg}) \rightarrow \tilde{\pi}^{-1}(\cX^{\reg})/W_{X'} \rightarrow \tilde{\pi}^{-1}(\cX^{\reg})/W_X.
\end{equation*}

In \cite{MN}, McGerty and Nevins proved a number of results about $\HC$ and $\Spr$ which together established the viability of symplectic Springer theory. The proofs carry over to produce the analogous statements for $\HC'$ and $\Spr'$, which are collected below. 

\begin{thm} \label[thm]{mn} Let $W_{X'}$ be the Namikawa Weyl group of $X'$. Then:

\begin{enumerate}

\item $\HC'$ and $\Spr'$ are semisimple perverse sheaves, and $\HC'$ is given by taking the IC extension of the local system on ${\cX'}^{\reg} := \tilde{\rho}^{-1}(\cX^{\reg})$ given by the regular representation.
\item $\Spr' \in \Perv_{\sympl}(X')$.
\item $i_{X'}^*(\HC') = \Spr'$.
\item There is a natural isomorphism $\C[W_{X'}] \cong \End(\HC')$.
\item There is an adjunction 
$$( - \otimes_{\C[W']} \Spr'): \C[W_{X'}] - \Mod \leftrightarrows \Perv_{\sympl}(X'): \Hom(\Spr', - ).$$

\end{enumerate}
\end{thm}

\begin{rem} \label[rem]{NoRes}

The fact that $\Spr \in \Perv_{\sympl}(X)$ is proven in \cite[Section 3.2]{MN} using a nearby cycles construction over the universal deformation $\cX \rightarrow B_X$. In particular, this implies that $\Spr$ does not depend on a choice of a symplectic resolution $\pi: Y \rightarrow X$. Similarly, $\Spr'$ can be obtained by the analogous nearby cycles construction over the universal deformation $\cX' \rightarrow B_{X'}$. This proves the second point above, and shows that $\Spr'$ doesn't depend on the choice of a symplectic resolution $\pi': Y \rightarrow X'$.

\end{rem}

Additionally, we can consider how the Springer and Harish-Chandra sheaves of different partial resolutions interact. In particular, since star and shriek pushforward agree for proper maps, proper base change implies that the following diagram commutes:
\begin{equation} \label{Spr}
\begin{tikzcd} [row sep=scriptsize]
	{\C_Y} && {\C_{\cY}} \\
	\\
	{\Spr'} && {\HC'} \\
	\\
	\Spr && \HC
	\arrow["{\pi_*'}"', maps to, from=1-1, to=3-1]
	\arrow["{\tilde{\pi}'_*}"', maps to, from=1-3, to=3-3]
	\arrow["{i_Y^*}"', maps to, from=1-3, to=1-1]
	\arrow["{i_{X'}^*}"', maps to, from=3-3, to=3-1]
	\arrow["{\tilde{\rho}_*}"', maps to, from=3-3, to=5-3]
	\arrow["{\rho_*}"', maps to, from=3-1, to=5-1]
	\arrow["{i_X^*}"', maps to, from=5-3, to=5-1]
\end{tikzcd}
\end{equation}
We can further apply $\End(-)$ to the above diagram, and using the functoriality of $(-)_!$ and $(-)^*$, we obtain:
\begin{equation} \label{End}
\begin{tikzcd} 
	\C && \C \\
	\\
	{\End(\Spr')} && {\C[W_{X'}]} \\
	\\
	{\End(\Spr)} && {\C[W_X]}
	\arrow["{i_Y^*}"', from=1-3, to=1-1]
	\arrow["{\pi_*'}"', hook, from=1-1, to=3-1]
	\arrow["{\tilde{\pi}'_*}"', hook, from=1-3, to=3-3]
	\arrow["{i_{X'}^*}"', from=3-3, to=3-1]
	\arrow["{\tilde{\rho}_*}"', hook, from=3-3, to=5-3]
	\arrow["{\rho_*}"', hook, from=3-1, to=5-1]
	\arrow["{i_X^*}"', from=5-3, to=5-1]
\end{tikzcd}
\end{equation}

\begin{rem}

If we consider partial resolutions as a poset with the partial ordering given by covering, the above reasoning shows that the constructions of $\Spr'$ and $\HC'$ are functorial. In particular, we could extend diagrams \eqref{Spr} and \eqref{End} to consider a longer chain of partial resolutions, or even the entire poset. For simplicity, we shall limit ourselves to the above diagrams for the remainder of the section. 

\end{rem}

Recall that if we let $Z := Y \times_X Y$ and $Z':= Y \times_{X'} Y$, then the Borel-Moore homology groups $H^{\BM}_{\bullet}(Z)$ and $H^{\BM}_{\bullet}(Z')$ have natural convolution structures, which are were studied extensively in \cite{CG}. Additionally, by \cite[Lemma 2.11]{Kal}, $\pi$ and $\pi'$ are semi-small, because they are symplectic resolutions. This has many consequences for the structure of $H^{\BM}_{\bullet}(Z)$ and $H^{\BM}_{\bullet}(Z')$ (see \cite[Section 8.9]{CG}), and in particular we have the following description of the endomorphism algebras in \eqref{End}. 

\begin{lemma}[{\cite[Theorem 5.4]{VG}}] \label[lemma]{SteinSpr}

The top homology groups $H^{\BM}_{top}(Z) \subset H^{\BM}_{\bullet}(Z)$ and $H^{\BM}_{top}(Z') \subset H^{\BM}_{\bullet}(Z')$ are subalgebras under convolution. Additionally, we have algebra isomorphisms $\End(\Spr) \cong H^{\BM}_{top}(Z)$ and $\End(\Spr') \cong H^{\BM}_{top}(Z')$.

\end{lemma}

As well as the map $\rho_*$ in diagram \eqref{End}, proper pushforward in Borel-Moore homology along the inclusion map $i: Z' \hookrightarrow Z$ provides us with a map
\begin{equation*}
i_*: H^{\BM}_{top}(Z') \rightarrow H^{\BM}_{top}(Z).
\end{equation*}
Both of these maps, along with isomorphisms witnessing \Cref{SteinSpr}, can be described in terms of units and counits of various maps in the following diagram below. Note that the $\pi_i$ and $\pi'_i$ are the natural projection maps and both the inner and outer squares are Cartesian.
\begin{equation*}
\begin{tikzcd}
	&&& Y \\
	\\
	Z && {Z'} && {X'} && X \\
	\\
	&&& Y
	\arrow["{\pi'}"', from=1-4, to=3-5]
	\arrow["\pi", from=1-4, to=3-7]
	\arrow["{\pi_1}", from=3-1, to=1-4]
	\arrow["{\pi_2}", from=3-1, to=5-4]
	\arrow["{\pi_1'}"', from=3-3, to=1-4]
	\arrow["i"', hook', from=3-3, to=3-1]
	\arrow["{\pi_2'}", from=3-3, to=5-4]
	\arrow["\rho", from=3-5, to=3-7]
	\arrow["{\pi'}", from=5-4, to=3-5]
	\arrow["\pi", from=5-4, to=3-7]
\end{tikzcd}
\end{equation*}
Thus we have the following lemma, whose details we leave to the reader to verify:
\begin{lemma}
The proper pushforward $i_*$ in Borel-Moore homology is isomorphic to the map $\rho_*$ in diagram \eqref{End}.
\end{lemma}

In particular, $H^{\BM}_{top}(Z)$ has a basis by the fundamental classes of the irreducible components of $Z$ (see \cite[Proposition 2.6.14]{CG}), and $i_*$ is the inclusion of the subspace generated by the irreducible components of $Z'$. Thus, $\End(\Spr')$ has the structure of a subalgebra of $\End(\Spr)$. As explained in \Cref{NoRes}, $\Spr'$ and $\Spr$ do not depend on the choice of a symplectic resolution $Y$, so the above construction assigns to each partial resolution of $X$ a canonical subalgebra of $\End(\Spr)$.

By \cite[Proposition 3.8.7]{PA}, \Cref{small} implies that the sheaves on the right column of diagram \eqref{Spr} are given by the intermediate extension of their restrictions to diagram \eqref{Gal}. Additionally, by \cite[Lemma 3.3]{MN}, these intermediate extensions induce isomorphisms on $\End(-)$. After restricting to \eqref{Gal}, the pushforward maps are given by:
\[\begin{tikzcd}
	{\Loc(\tilde{\pi}^{-1}(\cX^{\reg}))} \\
	\\
	{\Loc(\tilde{\rho}^{-1}(\cX^{\reg})) \cong W_{X'}-\Mod(\Loc(\tilde{\pi}^{-1}(\cX^{\reg})))} \\
	\\
	{\Loc(\cX^{\reg}) \cong W_X-\Mod(\Loc(\tilde{\pi}^{-1}(\cX^{\reg}))) }
	\arrow["{-  \otimes \C[W_{X'}]}", from=1-1, to=3-1]
	\arrow["{-  \otimes_{\C[W_{X'}]} \C[W_X]}", from=3-1, to=5-1]
\end{tikzcd}\]

In particular, the maps $\tilde{\pi}'_*$ and $\tilde{\rho}_*$ in diagram \eqref{End} are the standard embeddings. Thus, if we use $i_{X'}^*$ and $i_X^*$ in diagram \eqref{End} to induce actions of $W_{X'}$ and $W_X$ on $\Spr'$ and $\Spr$ respectively, the bottom square of diagram \eqref{Spr} is $W_{X'}$-equivariant.  Taking $W_{X'}$-invariants of diagram \eqref{Spr} we obtain 
\begin{equation} \label{inv}
\begin{tikzcd}
	{\Spr'^{W_{X'}}} && {\HC'^{W_{X'}}} \\
	\\
	{\Spr^{W_{X'}}} && {\HC^{W_{X'}}}
	\arrow["{i_{X'}^*}"', from=1-3, to=1-1]
	\arrow["{\tilde{\rho}_*}"', maps to, from=1-3, to=3-3]
	\arrow["{\rho_*}"', maps to, from=1-1, to=3-1]
	\arrow["{i_X^*}"', maps to, from=3-3, to=3-1]
\end{tikzcd}
\end{equation}

Since $\HC'$ is the intermediate extension of $\C[W_{X'}]$ on $\tilde{\rho}^{-1}(\cX^{\reg})$, $\HC'^{W_{X'}}$ is given by $\IC(\tilde{\rho}^{-1}(\cX^{\reg}), \C_{\tilde{\rho}^{-1}(\cX^{\reg})})$. Additionally, $\HC^{W_{X'}}$ is the intermediate extension of $\C[W]^{W_{X'}}$ on $\cX^{\reg}$. In particular,
\begin{equation*}
\End(\HC^{W_{X'}}) = \End_{W_{X}-\Mod(\Loc(\tilde{\pi}^{-1}(\cX^{\reg})))}(\C[W_X]^{W_{X'}}) \cong \cH(W_X, W_{X'}),
\end{equation*}
where $\cH(W_X, W_{X'})$ is the Hecke algebra of a finite group and a subgroup, given by the set of bi-invariant functions $\C[W_{X'} \backslash W_X/W_{X'}]$ under convolution. Thus, the induced map 
\begin{equation*}
i_X^*: \End(\HC^{W_{X'}}) \rightarrow \End(\Spr^{W_{X'}})
\end{equation*}
is an explicit construction of the universal action of the Hecke algebra on $W_{X'}$-invariants. 

The following corollary is inspired by \cite[Sections 2.6-2.8]{BM}, but requires the additional condition that $\cX'$ is rationally smooth. 

\begin{cor} \label[cor]{SymplFibersInv}

If $\cX'$ is rationally smooth then we have $\HC^{W_{X'}} \cong \tilde{\rho}_*(\C_{\cX'})$ and $\Spr^{W_{X'}} \cong \rho_*(\C_{X'})$.  In particular, given $x \in X$, $H^*(\rho^{-1}(x), \C) \cong H^*(\pi^{-1}(x), \C)^{W_{X'}}$.

\end{cor}

\begin{proof}

If $\cX'$ is rationally smooth, then by \cite[Proposition 8.2.21]{HTT}, 
\begin{equation*}
\HC'^{W_{X'}} \cong \IC(\tilde{\rho}^{-1}(\cX^{\reg}), \C_{\tilde{\rho}^{-1}(\cX^{\reg})}) \cong \C_{\cX'}. 
\end{equation*}
Thus by diagram \eqref{inv} we have $\HC^{W_{X'}} \cong \tilde{\rho}_*(\C_{\cX'})$.  Additionally, 
$$
\Spr'^{W_{X'}} \cong i_{X'}^*(\C_{\cX'}) \cong \C_{X'},
$$
and thus $\Spr^{W_{X'}} \cong \rho_*{\C_{X'}}$. By proper base change, the stalks of $\Spr$ and $\rho_*{\C_{X'}}$ are $H^*(\pi^{-1}(x), \C)$ and $H^*(\rho^{-1}(x), \C)$. Thus, taking the stalk at $x$ of the isomorphism $\Spr^{W_{X'}} \cong \rho_*{\C_{X'}}$ yields $H^*(\rho^{-1}(x), \C) \cong H^*(\pi^{-1}(x), \C)^{W_{X'}}$.
\end{proof}

\begin{rem}

In \cite{BM}, Borho and MacPherson study partial resolutions of the nilpotent cone. In that case the universal deformations $\cX'$ of the partial resolutions are all given by the smooth varieties $\tilde{\mathfrak{g}}^P \cong G \times_P \mathfrak{p}$, for different parabolic subgroups $P \subset G$ of a reductive lie group $G$. 

\end{rem}

\begin{Ex}

\Cref{SymplFibersInv} also applies to all the partial resolutions of Kleinian singularities, as we briefly explain below. Note that in this case the conclusion of \Cref{SymplFibersInv} is not a difficult computation to do by hand. 

Given a Kleinian singularity $S$, we recall from \cite[Proposition 3.1]{NaPD1} that the universal deformation of $S$ is given by taking a Slodowy slice $\cS$ to the subregular nilpotent orbit of the corresponding Lie algebra $\mathfrak{g}$. In particular, $\cS$ is an affine space, so it is smooth. Given a partial resolution $\rho: S' \to S$, let $I$ be the subset of the simple roots of $\mathfrak{g}$ that corresponds to projective lines contracted by the induced map $\pi': \tilde{S} \to S'$. Additionally, let $\mathcal{P}$ be the conjugacy class of parabolic subalgebras of $\mathfrak{g}$ that corresponds to $I$. The universal deformation $\cS'$ of $S'$ is given by taking the preimage of $\cS$ under the canonical map $\tilde{\eta}: \tilde{\mathfrak{g}}^\mathcal{P} \to \mathfrak{g}$. Note that $\cS$ is transverse to every adjoint orbit in $ \mathfrak{g}$ (see \cite[Section 2.2]{GG}), and this implies that it is transverse to the map $\tilde{\eta}$. In particular, this implies that $\cS' = \tilde{\eta}^{-1}(\cS)$ is smooth, so \Cref{SymplFibersInv} applies. 

\end{Ex}

The above result is false in general because there are many examples of symplectic singularities with nontrivial symplectic resolutions that have trivial Namikawa Weyl group. For example, the symplectic resolution $T^*\Proj^2 \rightarrow \overline{\mathbb{O}_{\text{subreg}}} \subset \mathcal{N}_{\mathfrak{sl}_3}$ has trivial Weyl group, but the fiber at zero is the zero section. However, in the cases where the corollary applies, it shows a strong relationship between the geometry and the combinatorics of a given partial resolution.

\end{document}